\documentclass{article}%
\usepackage{amsfonts}
\usepackage{amsmath}
\usepackage{amssymb}
\usepackage{amsxtra}
\usepackage{graphicx}
\usepackage{geometry}
\usepackage{color}
\usepackage{colortbl}
\usepackage{caption}
\usepackage[draft]{varioref}%
\providecommand{\U}[1]{\protect\rule{.1in}{.1in}}
\newtheorem{theorem}{Theorem}

\newtheorem{corollary}[theorem]{Corollary}

\newtheorem{lemma}[theorem]{Lemma}

\newtheorem{proposition}[theorem]{Proposition}
\newtheorem{remark}[theorem]{Remark}

\numberwithin{equation}{section}
\begin{document}

\title{Scalar boundary value problems on junctions \\of thin rods and plates.\\I. Asymptotic analysis and error estimates}
\author{R. Bunoiu\\Universit\'{e} de Lorraine, Institut Elie Cartan de Lorraine, \\UMR 7502, Metz, F-57045, France.\\renata.bunoiu@univ-lorraine.fr
\and G.Cardone\\Universit\`{a} del Sannio - Dipartimento di Ingegneria \\Piazza Roma, 21 - 84100 Benevento, Italy\\email: gcardone@unisannio.it
\and S.A.Nazarov\\Mathematics and Mechanics Faculty, St. Petersburg State University\\198504, Universitetsky pr., 28, Stary Peterhof, Russia.\\email: srgnazarov@yahoo.co.uk}
\maketitle

\begin{abstract}
We derive asymptotic formulas for {the} solutions of the mixed boundary value
problem for the Poisson equation on the union of a thin cylindrical plate and
several thin cylindrical rods. One of {the} ends of each rod is set into a
hole in the plate and the other one is supplied with the Dirichlet condition.
The Neumann conditions are imposed on the whole remaining part of the
boundary. Elements of the junction are assumed to have contrasting properties
so that the small parameter, i.e. the relative thickness, appears in the
differential equation, too, while the asymptotic structures crucially depend
on the contrastness ratio. Asymptotic error estimates are derived in
anisotropic weighted Sobolev norms.

\medskip


\medskip


\end{abstract}

\section{Introduction\label{sect1}}

\subsection{Formulation of the problem\label{sect1.1}}

Let $\mathbb{\omega}_{0}$ and $\mathbb{\omega}_{j\text{ }}$be domains in the
plane $\mathbb{R}^{2}$ bounded by smooth simple closed contours $\partial
\mathbb{\omega}_{p};$ here and everywhere in the paper $j=1,....J$ and
$p=0,....J$, while $J\in\mathbb{N}=\{  1,2,3,...\}  .$
Moreover, the summation over $j=1,...,J$ will be further denoted by $\sum_j$.
The closures
$\overline{\mathbb{\omega}}_{p}=\mathbb{\omega}_{p\text{ }}\cup\partial
\mathbb{\omega}_{p}$ are compact and the origin of the Cartesian coordinates
$y=(  y_{1},y_{2})  $ belongs to $\mathbb{\omega}_{j}.$ We fix some
points $P^{1},....,P^{{J}}$ inside $\mathbb{\omega}_{0}$, $P^{j}\neq P^{k}$
for $j\neq k,$ and introduce the thin plate
\begin{equation}
\Omega_{0}(  h)  =\{  x=(  y,z)  \in\mathbb{R}%
^{3}:y\in\mathbb{\omega}_{0},\text{ }\zeta:=h^{-1}z\in(  0,1)
\}  \label{1.1}%
\end{equation}
and the thin rods%
\begin{equation}
\Omega_{j}(  h)  =\{  x:\eta^{j}=:h^{-1}(  y-P^{j}%
)  \in\mathbb{\omega}_{j},\text{ \ }z\in(  0,l_{j})
\}  \label{1.2}%
\end{equation}
where $h\in(  0,h_{0}]  $ is a small parameter and $l_{1}%
,...,l_{{J}}$ are fixed positive numbers. The bound $h_{0}>0$ is chosen such
that the closures of the small sets
$\mathbb{\omega}_{j}^{h}=\{  y:\eta^{j}\in\mathbb{\omega}_{j}\}$
are contained in the domain $\mathbb{\omega}_{0}$ for all $h\in(
0,h_{0}]  $. In the sequel, if necessary, we may diminish this bound but
always keep the notation $h_{0}.$

\begin{figure}
[ptb]
\begin{center}
\includegraphics[
height=1.5688in,
width=2.9922in
]%
{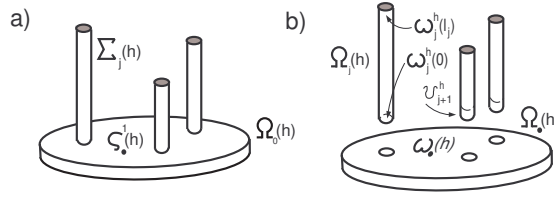}%
\caption{The junction (a) and its elements (b). The Dirichlet zones are
shaded.}%
\label{f1}%
\end{center}
\end{figure}

The junction
\begin{equation}
\Xi(  h)  =\Omega_{\bullet}(  h)  \cup\Omega_{1}(
h)  \cup....\cup\Omega_{J}(  h)  \label{1.4}%
\end{equation}
of the plate and rods, see fig. \ref{f1},a, involves the intact rods
(\ref{1.2}) but the plate $\Omega_{\bullet}(  h)  $ with
cylindrical holes (column sockets), fig. \ref{f1},b,
\begin{equation}
\Omega_{\bullet}(  h)  =\mathbb{\omega}_{\bullet}(  h)
\times(  0,h)  , \quad \mathbb{\omega}_{\bullet}(  h)  =\mathbb{\omega}_{0}%
\mathbb{\diagdown}(  \overline{\mathbb{\omega}_{1}^{h}}\cup
.....\cup\overline{\mathbb{\omega}_{{J}}^{h}})  .\label{1.5}%
\end{equation}
%TCIMACRO{\FRAME{ftbpFU}{2.9014in}{1.1096in}{0pt}{\Qcb{The junction (a) and its
%elements (b). The Dirichlet zones are shaded.}}{\Qlb{f1}}{junction1new.eps}%
%width 2.9014in;  height 1.1096in;  depth 0pt;  original-width 7.1676in;
%original-height 2.7242in;  cropleft "0";  croptop "1";  cropright "1";
%cropbottom "0";  filename 'Junction1new.eps';file-properties "XNPEU";}} }%
%BeginExpansion
%EndExpansion
The bases of the intact (\ref{1.1}) and perforated (\ref{1.5}) plates are
denoted respectively by
\begin{equation}
\varsigma_{0}^{i}(  h)  =\{  x:y\in\mathbb{\omega}%
_{0},\ z=ih\}  ,\ \ \ \varsigma_{\bullet}^{i}(  h)  =\{
x:y\in\mathbb{\omega}_{\bullet},\ z=ih\}  ,\ \ i=0,1, \label{1.7}%
\end{equation}
and the common lateral side by $\upsilon_{0}(  h)  =\partial
\mathbb{\omega}_{0}\times(  0,h)  .$ The ends of the rod
(\ref{1.2}) are
\begin{equation}
\mathbb{\omega}_{j}^{h}(  0)  =\mathbb{\omega}_{j}^{h}%
\times\{  0\}  ,\ \ \mathbb{\omega}_{j}^{h}(  l_{j})
=\mathbb{\omega}_{j}^{h}\times\{  l_{j}\}  \label{1.8}%
\end{equation}
while the lateral side of the rod is divided into two parts%
\begin{equation}
\Sigma_{j}(  h)  =\partial\mathbb{\omega}_{j}^{h}\times(
h,l_{j})  ,\ \ \mathbb{\upsilon}_{j}^{h}=\partial\mathbb{\omega}_{j}%
^{h}\times(  0,h)  , \label{1.9}%
\end{equation}
the latter being the contact zone of $\Omega_{\bullet}(  h)  $ and
$\Omega_{j}(  h)  .$ These sets are indicated in fig. \ref{f1}.

In the junction (\ref{1.4}) we consider the Poisson equations%
\begin{align}
-\Delta_{x}u_{0}(  h,x)   &  =f_{0}(  h,x)
,\ \ \ x\in\Omega_{\bullet}(  h)  ,\label{1.10}\\
-\gamma_{j}(  h)  \Delta_{x}u_{j}(  h,x)   &
=f_{j}(  h,x)  ,\ \ \ x\in\Omega_{j}(  h)  ,
\label{1.11}%
\end{align}
equipped with the Neumann and Dirichlet boundary conditions%
\begin{align}
\partial_{\nu}u_{0}(  h,x)   &  =0,\ \ \ \ x\in\Sigma_{\bullet
}(  h)  =\varsigma_{\bullet}^{0}(  h)  \cup
\varsigma_{\bullet}^{1}(  h)  \cup\mathbb{\upsilon}_{0}(
h) \label{1.12}\\
\gamma_{j}(  h)  \partial_{\nu}u_{j}(  h,x)   &
=0,\ \ \ x\in\Sigma_{j}(  h)  \cup\mathbb{\omega}_{j}^{h}(
0)  , \label{1.13}\\
u_{j}(  h,x)  &=0,\ \ \ x\in\mathbb{\omega}_{j}^{h}(
l_{j})  , \label{1.14}%
\end{align}
together with the transmission conditions%
\begin{align}
u_{0}(  h,x)   &  =u_{j}(  h,x)  ,\ \ \ \ x\in
\mathbb{\upsilon}_{j}^{h},\label{1.15}\\
\partial_{\nu}u_{0}(  h,x)   &  =\gamma_{j}(  h)
\partial_{\nu}u_{j}(  h,x)  ,\ \ \ x\in\mathbb{\upsilon}_{j}^{h}.
\label{1.16}%
\end{align}
Here, $\Delta_{x}$ is the Laplace operator in
$x=(  y,z)  =(  x_{1},x_{2},x_{3})  ,$ $u_{0}$ and
$u_{j}$ are restrictions of the function $u$ on the subdomains $\Omega
_{\bullet}(  h)  $ and $\Omega_{j}(  h)  ,$
$f_{0}=f|_{\Omega_{\bullet}(  h)  }$ and $f_{j}=f|_{\Omega
_{j}(  h)  }$ having the similar meaning, $\partial_{\nu}=\nu
\cdot\nabla_{x}$ is the directional derivative, while $\nabla_{x}%
=\operatorname{grad}$ and $\nu$ is the unit vector of the outward normal on
the surfaces $\partial\overline{\Xi(  h)  }$ and $\partial
\Omega_{j}(  h)  .$ Notice that
$\partial_{\nu}=\partial_{z}=\partial/\partial z\text{ \ \ on }\varsigma
_{0}^{1}(  h)  \text{ \ and }\partial_{\nu}=-\partial
_{z}\text{\ \ on }\varsigma_{0}^{0}(  h)  .
$

The coefficient%
\begin{equation}
\gamma_{j}(  h)  =\gamma_{j}h^{-\alpha}>0 \label{1.17}%
\end{equation}
in (\ref{1.11}), (\ref{1.13}) and (\ref{1.16}) describes contrasting
properties of elements in the junction (\ref{1.4}). Namely, regarding
$u(  h,x)  $ as a stationary temperature field in $\Xi(
h)  ,$ we obtain a homogeneous body in the case $\alpha=0,$ $\gamma
_{j}=1$ but in the case $\alpha>0$ the conductivity of the rods $\Omega
_{j}(  h)  $ is much bigger than of the plate $\Omega_{\bullet
}(  h)  .$ In what follows we deal with two typical cases%
\begin{equation}
\alpha=0\text{ \ and \ }\alpha=1. \label{1.18}%
\end{equation}
As mentioned above, the first case may be attributed to the homogeneous
junction. In the second case the integral conductivity of the plate vertical
segment $(  0,h)  $ and of the rod cross-section $\omega_{j}^{h},$
that are $h$ and $\gamma_{j}(  h)  \operatorname*{mes}_{2}%
\omega_{j}^{h}=h\gamma_{j}\operatorname*{mes}_{2}\omega_{j}$ respectively,
become of the same order in $h.$ The latter complicates both, the asymptotic
ans\"{a}tze for the solution $u(  h,x)  $ of problem and the asymptotic procedure. To construct the
asymptotics as $h\rightarrow+0$ and to justify it by deriving error estimates
are just the main goal of the paper.

The variational formulation of problem (\ref{1.10})-(\ref{1.16}) reads: to
find a function $u\in H_{0}^{1}(  \Xi(  h)  ;\Gamma(
h)  )  $ satisfying the integral identity \cite{Lad}%
\begin{equation}
a(  u,v;\Xi(  h)  )  =(  f,v)  _{\Xi(
h)  }\ \ \ \forall v\in H_{0}^{1}(  \Xi(  h)
;\Gamma(  h)  )  . \label{1.19}%
\end{equation}
Here, $H_{0}^{1}(  \Xi(  h)  ;\Gamma(  h)  )
$ is a subspace of functions in the Sobolev space $H^{1}(  \Xi(
h)  )  $ which meet the Dirichlet condition (\ref{1.14}) and
$(  \ ,\ )  _{\Xi(  h)  }$ is the natural scalar product
in the Lebesgue space $L^{2}(  \Xi(  h)  )  .$
Furthermore,
\begin{equation}
a(  u,v;\Xi(  h)  )  =(  \nabla_{x}u_{0},\nabla
_{x}v_{0})  _{\Omega_{\bullet}(  h)  }+{\textstyle\sum\nolimits_{j}}%
\gamma_{j}(  h)  (  \nabla_{x}u_{j},\nabla_{x}v_{j})
_{\Omega_{j}(  h)  } \label{1.20}%
\end{equation}
and $\Gamma(  h)  $ is the union of the upper ends $\omega_{j}%
^{h}(  l_{j})  $ of the rods. As usual, the integral identity is obtained by multiplying the equations (\ref{1.10}), (\ref{1.11}){ by} test functions $v_{0},$ $v_{j}$ and
integrating by parts in $\Omega_{\bullet}(  h)  ,$ $\Omega
_{j}(  h)  $ while taking into account the boundary (\ref{1.12}),
(\ref{1.13}) and transmission (\ref{1.16}) conditions for $u$ and also the
stable conditions (\ref{1.14}), (\ref{1.15}), absorbed in the space
$H_{0}^{1}(  \Xi(  h)  ;\Gamma(  h)  )  $ and
therefore given to the test function $v.$

\subsection{Asymptotic analysis of junctions\label{sect1.2}}

Bodies made of thin elements are met everywhere in our daily life; one may
think about miscellaneous mechanisms and their details, bridges and wheels
with spokes, chairs and tables, etc. Even flora and fauna give examples of all
kinds of such junctions. Mathematical studies of junctions of domains with
different limit dimensions have provided several approaches diversified by
levels of rigor and ways to formulate results.

%In such structures, at least one of the portion of the whole three-dimensional
%structure has a small thickness or a small diameter, which is proportional to
%a small dimensionless parameter (corresponding to the parameter denoted by $h$
%in our case). In order to investigate such problems for elastic bodies, P.G.
%Ciarlet (see \cite{Ciarlet2, Ciarlet3, Ciarlet4}) (making use of some ideas
%advocated by J.L. Lions (see \cite{Lions2}) proposed the application of the
%asymptotic expansions method to the variational formulation of the problem,
%proved convergence results when the small parameter tends to zero and derived
%the limit problem, which rigorously involved limit junction conditions. This
%problem was a coupled, pluri-dimensional variational problem of a new type;
%the pluri-dimensional character comes from the fact that the problem is
%simultaneously posed over sets of different dimensions.
%
%In the case when at least one portion of the three-dimensional structure is a
%rod,

\begin{figure}
[ptb]
\begin{center}
\includegraphics[
height=0.8311in,
width=2.7121in
]%
{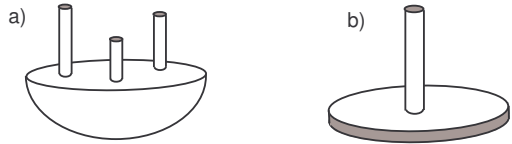}%
\caption{The juncton of a massive body and rods (a). The junction of a
mushroom type (b). The Dirichlet zones are shaded.}%
\label{f2}%
\end{center}
\end{figure}

Junctions of type "massive body/thin rods", fig \ref{f2},a, were
inspected from all sides, most intensely among other types. Asymptotic
formulas together with error estimates for fields of various physical nature
are obtained and hybrid variational-asymptotic models are derived, see, e.g.,
the papers \cite{JT, B, A, KoMaMo1, na220, Gad, na344} for scalar differential
equations and \cite{BG3, BG4, LeSP, COT, na212, KoMaMo2, KoMaMo3, na341, na345, na514}
for elasticity and other systems; vast citation lists can also be found in the
monograph \cite{LD, KoMaMo} and the review paper \cite{na396}.%

%TCIMACRO{\FRAME{ftbpFU}{2.7614in}{1.19in}{0pt}{\Qcb{The juncton of a massive
%body and rods (a). The junction of a mushroom type (b). The Dirichlet zones
%are shaded.}}{\Qlb{f2}}{junction2new.eps}%
%width 2.7614in;  height 1.19in;  depth 0pt;  original-width 7.1252in;
%original-height 3.0571in;  cropleft "0";  croptop "1";  cropright "1";
%cropbottom "0";  filename 'Junction2new.eps';file-properties "XNPEU";}} }%
%BeginExpansion
%EndExpansion

Junction of thin plates and rods have been considered as well, see
\cite{na285, G1, G2} for scalar equations and \cite{Gr1, Gr2, French, BG1,
BG2} for elasticity. However, the known results, cf. \cite{G1, G2}, for scalar
equations concern only the case \thinspace$J=1,$ see fig. \ref{f2},b, and,
more importantly, the lateral side $\upsilon_{0}(  h)  $ of the
plate $\Omega_{0}(  h)  $ used to be endowed with the Dirichlet
condition, that is%
\begin{align}
-\partial_{z}u_{0}(  h,y,0)   &  =\partial_{z}u_{0}(
h,y,h)  =0,\ \ \ y\in\omega_{\bullet}(  h)  ,\label{1.d}\\
u_{0}(  h,y,z)   &  =0,\ \ \ (  y,z)  \in\upsilon
_{0}(  h)  :=\partial\omega_{0}\times(  0,h)  ,\nonumber
\end{align}
instead of (\ref{1.12}). In elasticity the Dirichlet condition means that the
edge of the "mushroom cap" in fig. \ref{f2},b, must be clamped as well as the
lower part of "mushroom leg". Besides, for the stationary heat equation under
consideration, the constant (null) temperature in problem (\ref{1.10}),
(\ref{1.11}), (\ref{1.d}), (\ref{1.13})-(\ref{1.16}) is maintained not only at
the "soles" of the rods but on the lateral side of the plate, too. Evidently,
both cases (\ref{1.12}) and (\ref{1.d}) may occur in practice and in the
sequel we state the scalar boundary value problem on junction (\ref{1.4}) just
in the same way as for the above-mentioned junctions "massive body/thin rods"
in \cite{KoMaMo1, na220, na344}, \cite{KoMaMo} and others.

Another departure from results obtained in \cite{G1, G2} and other
publications consists in the derivation of estimates of the asymptotic
remainders while the previous treatment asserts convergence results without
estimating the convergence rate. For a spatial junction of type "thin
plate/thin rods" formulation of asymptotic formulas is a matter of principle
because one of the limit problem is planar and, as discovered in \cite{Ilin1}
(see also \cite{Ilinbook}, \cite[Ch.2, 4, 5]{MaNaPl}) a singular perturbation
of the boundary may lead to the rational dependence on $\ln h$ of asymptotic terms expansions. This indeed happens in the case $\alpha=1$
(cf. Section \ref{sect3.1}) due to the appearance of logarithm in the
fundamental solution of the Laplacian in $\mathbb{R}^{2}.$ Asymptotic series
in $\vert \ln h\vert ^{-1},$ of course, are available but, being
convergent, they however provide rather rough proximity of order $\vert
\ln h\vert ^{-N}$ while taking into account the rational dependence
brings order $h^{\delta},$ $\delta>0.$

The case $\alpha=0$ generates another effect: asymptotic expansions of the
solution $u(  h,x)  $ of problem (\ref{1.10})-(\ref{1.16}) gain
terms of order $h^{-1}$ even for a smooth uniformly bounded right-hand sides
$f_{p}(  h,x)  $ in the differential equations, see Section
\ref{sect3.3}. Moreover the next terms becomes linear functions in $\ln h,$
i.e., also grow when $h\rightarrow+0.$ These facts make a formulation of the
asymptotic decomposition as a convergent result doubtful, cf. a pour
formulation in Section \ref{sect4.6}.

We emphasize that both the above-mentioned peculiarities of the asymptotic
behavior of the solution $u(  h,x)  $ vanish if the boundary
conditions (\ref{1.12}) are replaced with (\ref{1.d}). In particular all
solutions with logarithmic singularities we use below, disappear from the
analysis of main asymptotic terms and the convergence theorems of \cite{G1,
G2}.

\subsection{The asymptotic method\label{sect1.3}}

%Asymptotic analysis of singularly perturbed elliptic boundary value problems is
%of interest in physics, medicine, material science etc. relatively
%to different type of structures. There are recent mathematical publications
%describing the behavior of thin tubular structures (cf. \cite{CaCoPa2, CaCoPa, CaFaPa, CaPaSi, Panas}) or array structures (cf. \cite{CioSJP, Panas, BeCaGoPa}), the behavior of quantum or acoustic waveguides in the case of perturbation by mixed boundary conditions (cf. \cite{BoCa2, BoCa3, BoBuCa, BoBuCa2, BoBuCa3}) or in the case of an oscillating boundary (cf. \cite{BoCa4}).

Asymptotic expansions of solutions to elliptic boundary value problems in
domains with singularly perturbed boundary can be constructed by means of two,
certainly equivalent methods, namely the method of matched asymptotic
expansions and the method of compound asymptotic expansions; we refer,
respectively, to the monographs \cite{VanDyke, Ilinbook} and \cite{MaNaPl,
KoMaMo}, these lists could be elongated quite much, where a complete
description of the methods is given and, furthermore, the equivalency is
clarified in \cite[Ch.2]{MaNaPl}. Previous studies of junctions of type
"massive body/thin rods" were mainly based on the method of compound expansion
because just this method elucidates boundary layer effects in the vicinity of
the juncture zones. These effects are described by boundary value problems in
the union of a semi-cylinder and a half-space $\mathbb{R}_{+}^{3}$ while the
intrinsic power-law decay of their solutions makes the method of compound
expansions much more preferable.

For the homogeneous ($\alpha=0$) junction $\Xi(  h)  $ of thin
plate and rods the boundary layer terms are solutions of the Neumann problem
in the union $\Xi_{j}$ of the semi-cylinder $Q_{j}$ and the perforated layer
$\Lambda_{j},$ see (\ref{2.35}). Moreover, the contrasting properties of the
junction elements, that is $\alpha=1$ in (\ref{1.17}), split the problem in
$\Xi_{j}$ into two independent problems in $Q_{j}$ and $\Lambda_{j},$ cf.
Sections \ref{sect2.4} and \ref{sect3.3}. Unfortunately, the solutions of the
problems in $Q_{j}$ and $\Lambda_{j}$ no longer get a decay at infinity in the
layer but quite the contrary{ they} gain a logarithmic growth. The latter
brings a complication into the application of compound expansions (cf. \cite[Ch.2]{MaNaPl}) and that is why in the paper we use the method of
matched asymptotic expansions. Namely we construct outer and inner expansions
which, respectively, involve singular solutions of the limit problems in
$\omega_{0}$ and growing solutions of the limit problems in $\Xi_{j}$ and
$\Lambda_{j}.$ The expansions acquire a family of free constants, actually
$A_{0}$ and $A_{1},...,A_{J}$ in Proposition \ref{prop2.1}, but the standard
matching procedure, which equalizes main asymptotic terms of the outer
expansions as $x\rightarrow P^{j}$ and main asymptotic terms of the inner
expansions at infinite, results into a system of linear algebraic equations
{which allow us } to determine the free constants. It should be stressed that
every particularity in the asymptotic behavior of the solution $u(
h,x)  $ originates in the structure of the algebraic system which, in
its turn, is totally predetermined by {the} general attributes of the junction.

We also point out that the choice of the method of matched asymptotic
expansions is prompted by a goal of our proceeding paper, namely to create an
asymptotic variational model of a junction of type "thin plate/thin rods"
which involves certain self-adjoint extensions of differential operators of
the limit problems. Indeed, it is the paper \cite{na239} where an intrinsic
relationship of the method and the technique of self-adjoint extensions was
revealed for elliptic problem in singularly perturbed domains. We emphasize
that in \cite{Lions} an application of self-adjoint extensions for junctions
of domains with different limit dimensions was formulated as an open question;
however, only junction of type "massive body/thin rods" have been examined in
this way, cf. \cite{na344, na341, na345}.

Once asymptotics is clarified, we use a construction of a global approximation
on involving cut-off functions with "overlapping supports" (see \cite[Ch.2]{MaNaPl}) which helps to satisfy all boundary and transmission conditions and
to diminish residuals left in the differential equations (\ref{1.10}) and
(\ref{1.11}). After applying a priori estimates of solutions to the
variational problem (\ref{1.19}), we finally cleanse the used intermediate
approximate solution from those cut-off functions and formulate error
estimates on each element of the junction. The latter explains soundly
convergence theorems and will be used to justify the above-mentioned
asymptotic-variational models.

\subsection{Architecture of the paper\label{sect1.4}}

In Section \ref{sect2} we derive and analyze different limit problems whose
solutions are used in Section \ref{sect3} in order to construct the formal
asymptotics of the solution $u(  h,x)  $ to the original problem in
the junction cases (\ref{1.18}). We also indicate in Sections \ref{sect3.2}
and \ref{sect3.4} significant simplifications occurring due to the restriction
$J=1$ and the Dirichlet boundary condition on the plate lateral side
$\upsilon_{0}(  h)  =\partial\omega_{0}\times(  0,h)  .$
A justification of the general asymptotic procedures in Section \ref{sect3.1}
for $\alpha=1$ and in Section \ref{sect3.3} for $\alpha=0$ is given in Section
\ref{sect4} which starts with derivation of several weighted estimates
summarized in Theorem \ref{th4.1}. After listing some necessary restrictions
on the problem data in Section \ref{sect4.3}, we prove in Section
\ref{sect4.5} Theorem \ref{th4.N3} for $\alpha=1$ and in Section \ref{sect4.6}
Theorem \ref{th4.M3} for $\alpha=0$ which provide sharp estimates of the
remainders. These estimates also allow us to conclude in Corollaries
\ref{cor4.N4} and \ref{cor4.M4} evident assertions on convergence of the
solution components $u_{0}(  h,x)  $ and $u_{j}(  h,x)
$ as $h\rightarrow+0.$

Finally in Section \ref{sect5} we outline available generalizations in our
present formulation of the problem.

\section{Limit problems\label{sect2}}

We here discuss the limit problems whose solutions form the asymptotic
expansions of the solution $u(  h,x)  $ of problem (\ref{1.10}%
)-(\ref{1.16}) in the junction (\ref{1.4}) specified in (\ref{1.17}) and
(\ref{1.18}). The first couple of the limit problems for rods and plates is
rather standard, cf. \cite{MaNaPl, KoMaMo, Panas, Nabook} and others, and we
outline them briefly while paying the most attention to singular solutions of
the Neumann problem in $\omega_{0}$ which play an important role in the outer
expansion in the plate $\Omega_{\bullet}(  h)  .$ Other limit
problems appear in the union of a layer and a semi-cylinder or in a perforated
layer, cf. fig. \ref{f3}, and are intended to describe the inner expansions in
the vicinity of the junction zones. These problems and the decompositions of
their solutions are studied in Section \ref{sect2.3} and \ref{sect2.4} in detail.

\subsection{The limit problems for the rods\label{sect2.1}}

We assume that the right-hand side of the differential equation (\ref{1.11})
in the rod $\Omega_{j}(  h)  $ takes the form%
\begin{equation}
f_{j}(  h,x)  =h^{-1-\alpha}f_{j}^{\bot}(  \eta^{j},z)
+h^{-\alpha}f_{j}^{0}(  z)  +\widetilde{f}_{j}(  h,x)  ,
\label{2.1}%
\end{equation}
where $f_{j}^{\bot}$ and $f_{j}^{0}$ are some functions on $\Omega_{j}(
1)  =\omega_{j}\times(  0,l_{j})  $ and $(
0,l_{j})  ,$ respectively, while
\begin{equation}
\langle {f^\bot}\rangle_j(z):=\int_{\omega_{j}}f_{j}^{\bot}(  \eta,z)  d\eta=0,\ \ \ z\in(
0,l_{j})  . \label{2.2}%
\end{equation}
In this section we suppose that $f_{j}^{\bot}$ and $f_{j}^{0}$ are smooth but
a necessary restriction on their smoothness as well as a smallness of the
remainder $\widetilde{f}_{j}$ will be formulated in Section \ref{sect4.2}.

Taking (\ref{2.1}) and (\ref{1.17}) into account we readily accept the
asymptotic ansatz
\begin{equation}
u_{j}(  h,x)  =U_{j}^{0}(  z)  +hU_{j}^{1}(
\eta^{j},z)  +h^{2}U_{j}^{2}(  \eta^{j},z)  +... \label{2.3}%
\end{equation}
Here and in what follows, dots stand for low order terms, inessential for our
formal asymptotic analysis performed in Section \ref{sect2}. Inserting
(\ref{2.3}) and (\ref{2.1}) into (\ref{1.11}) and (\ref{1.13}) we derive the
recurrent sequence of the Neumann problems in the cross-section $\omega_{j}%
\ni\eta^{j}$ with the parameter $z\in(  0,l_{j})  ,$%
\begin{equation}
-\gamma_{j}\Delta_{\eta}U_{j}^{k}(  \eta^{j},z)  =F_{j}^{k}(
\eta^{j},z)  ,\ \ \ \eta^{j}\in\omega_{j},\ \ \ \ \ \ \ \ \gamma
_{j}\partial_{\nu(  \eta)  }U_{j}^{k}(  \eta^{j},z)
=0,\ \ \ \eta^{j}\in\partial\omega_{j}, \label{2.4}%
\end{equation}
where $k=0,1,2$ and
\[
F_{j}^{0}=0,\ F_{j}^{1}(  \eta^{j},z)  =f_{j}^{\bot}(
\eta^{j},z)  ,\ F_{j}^{2}(  \eta^{j},z)  =f_{j}^{0}(
z)  -\gamma_{j}\partial_{z}^{2}U_{j}^{-1}(  z)  .
\]
Clearly, relations (\ref{2.4}) with $k=0$ are fulfilled. The problem
(\ref{2.4}) with $k=1$ has a solution due to the orthogonality condition
(\ref{2.2}) and this solution can be subject to the orthogonality conditions $\langle {U_{j}^{1}}\rangle_j(z)=0$.
%\begin{equation}
%\int_{\omega_{j}}U_{j}^{1}(  \eta^{j},z)  d\eta^{j}=0. \label{2.5}%
%\end{equation}
Finally, the problem (\ref{2.4}) with $k=2$ becomes solvable provided%
\begin{equation}
-\gamma_{j}\vert \omega_{j}\vert \partial_{z}^{2}U_{j}^{0}(
z)  =\vert \omega_{j}\vert f_{j}^{0}(  z)
,\ \ \ z\in(  0,l_{j})  , \label{2.6}%
\end{equation}
where the factor $\vert \omega_{j}\vert =\operatorname*{mes}%
_{2}\omega_{j},$ the area of the rod cross-section, appears after integration
of $F_{j}^{2}(  \eta^{j},z)  $.

We supply the ordinary differential equation (\ref{2.6}) with the condition%
\begin{equation}
U_{j}^{c}(  l_{j})  =0 \label{2.7}%
\end{equation}
inherited from the Dirichlet condition (\ref{1.14}). The other condition at the point $z=0$ which closes the limit
problem for the rod, depends on the exponent $\alpha$ in (\ref{1.17}) and we
will find it in Section \ref{sect3} only.

\subsection{The limit problem for the plate\label{sect2.2}}

Similarly to (\ref{2.1})-(\ref{2.3}), we set%
\begin{align}
f_{0}(  h,x)   &  =h^{-1}f_{0}^{\bot}(  y,\zeta)
+f_{0}^{0}(  y)  +\widetilde{f}_{0}(  h,y)
,\label{2.8}\\
\langle f^\bot_0\rangle (\zeta):=&\int_{0}^{1}f_{0}^{\bot}(  y,\zeta)  d\zeta   =0,\ \ \ y\in
\omega_{0},\label{2.9}\\
u_{0}(  h,x)  &=U_{0}^{0}(  y)  +hU_{0}^{1}(
y,\zeta)  +h^{2}U_{0}^{2}(  y,\zeta)  +...\label{2.10}
\end{align}
We insert the asymptotic ans\"{a}tze (\ref{2.10}) and (\ref{2.8}) into the equation (\ref{1.10}) in the perforated plate (\ref{1.5}) as well
as into the boundary conditions (\ref{1.12}) at its bases $\varsigma_{\bullet
}^{0}(  h)  $ and $\varsigma_{\bullet}^{1}(  h)  ,$ see
(\ref{1.7}). In this way we obtain a recurrent sequence of differential
equations in the fast variable $\zeta=h^{-1}z\in(  0,1)  $ with the
Neumann conditions at $\zeta=0$ and $\zeta=1.$ Owing to the
orthogonality condition (\ref{2.9}) we solve the problem%
\begin{equation}
-\partial_{\zeta}^{2}U_{0}^{k}(  y,\zeta)  =F_{0}^{k}(
y,\zeta)  ,\ \ \ \zeta\in(  0,1)  ,\ \ \ \ \ \ \partial
_{\zeta}U_{0}^{k}(  y,0)  =\partial_{\zeta}U_{0}^{k}(
y,1)  =0\label{2.11}%
\end{equation}
with the index $k=1$ and the right-hand side $F_{0}^{1}$ from the list%
\[
F_{0}^{0}=0,\ F_{0}^{1}(  y,\zeta)  =f_{0}^{\bot}(
y,\zeta)  ,\ F_{0}^{2}(  y,\zeta)  =f_{0}^{0}(
y)  -\Delta_{y}U_{0}^{-1}(  y)  .
\]
Since this solution is defined up to an additive constant in $\zeta,$ we can
satisfy the natural requirement $\langle U^1_0\rangle (\zeta)=0$.
%\begin{equation}
%\int_{0}^{1}U_{0}^{1}(  y,\zeta)  d\zeta=0.\label{2.12}%
%\end{equation}
Noting that relation (\ref{2.11}) with $k=-1$ is obviously met, we observe
that the problem (\ref{2.11}) with $k=1$ gets a solution if and only if
\begin{equation}
-\Delta_{y}U_{0}^{0}(  y)  =f_{0}^{0}(  y)
,\ \ \ y\in\omega_{0}.\label{2.13}%
\end{equation}
By virtue of the Neumann boundary condition (\ref{1.12}) at the lateral side
$\upsilon_{0}(  h)  $ of the plate, we impose
\begin{equation}
\partial_{\nu}U_{0}^{0}(  y)  =0,\ \ y\in\partial\omega
_{0},\label{2.14}%
\end{equation}
and regard (\ref{2.13}), (\ref{2.14}) as the limit problem for the plate
$\Omega_{\bullet}(  h)  .$

In contrast to the limit problems for the rods $\Omega_{j}(  h)  $
which involve the Dirichlet conditions (\ref{2.7}) and hence are uniquely
solvable, problem (\ref{2.13}),
(\ref{2.14}) has no bounded solution in the case%
\begin{equation}
\langle f_{0}^{0}\rangle_0:=\int_{\omega_{0}}f_{0}^{0}(  y)  dy\neq0. \label{2.15}%
\end{equation}
A reason of this lack is hidden in the formal asymptotic procedure performed
above. Indeed, a priori the three-dimensional Poisson equation (\ref{1.10})
holds true only in the plate $\Omega_{\bullet}(  h)  $ with holes
and the boundary conditions at the perforated bases $\varsigma_{\bullet}^{i}.$
Thus, assuming $U_{0}^{0}$ in (\ref{2.10}) to be smooth in the intact domain
$\omega_{0}$ we somehow extend the differential equation onto the small
domains $\theta_{j}^{h}=\omega_{j}^{h}\times(  0,h)  $ which shrink
to the points $P^{j}$ as $h\rightarrow+0.$ In other words, there is no
intrinsic reason to deal with a smooth function $U_{0}^{0}$ so that $U_{0}%
^{0}$ may have singularities at $P^{1},...,P^{J}.$

A singular solution of the problem (\ref{2.13}), (\ref{2.14}) exists
unconditionally and, to construct such solution, we recall the notion of the
(generalized) Green function, cf. \cite{Smir2}. Let $G(  y,P)  $ be
a distributional solution to the problem%
\begin{equation}
-\Delta_{y}G(  y,P)   =\delta(  y-P)  -\vert
\omega_{0}\vert ^{-1},\ y\in\omega_{0},\ \ \partial_{\nu
}G(  y,P)  =0,\ y\in\partial\omega_{0},\ \ \langle G(\cdot,P)\rangle_0 =0, \label{2.16}%
\end{equation}
where $\vert \omega_{0}\vert =\operatorname*{mes}_{2}\omega_{0}$
and $\delta$ is the Dirac mass. We set $G_{j}(  y)  =G(  y,P^{j})$, and $y^{j}=y-P^{j},$ $r_{j}=\vert y^{j}\vert .$ There hold the
representations%
\begin{equation}
G_{j}(  y)  =\delta_{j,k}(2\pi)^{-1}\ln(1/r_{k})
+G_{jk}+O(  r_{k})  ,\ \ \ r_{k}\rightarrow+0,\ \ k=1,...,J.
\label{2.18}%
\end{equation}
Here, $\delta_{j,k}$ is the Kronecker symbol, $-(  2\pi)  ^{-1}%
\ln\vert y\vert $ is the fundamental solution of the Laplacian in
the plane, and $G_{jk}$ are some constants composing the $J\times J-$matrix%
\begin{equation}
G=(  G_{jk})  _{j,k=1}^{J}. \label{2.19}%
\end{equation}
The following calculation uses relations (\ref{2.16}), (\ref{2.18}) and shows
that the matrix (\ref{2.19}) is symmetric:%
\begin{align}
0  &  =\frac{1}{\vert \omega_{0}\vert }\int_{\omega_{0}}(
G_{j}(  y)  -G_{k}(  y)  )  dy=\lim_{\varrho
\rightarrow+0}\frac{1}{\vert \omega_{0}\vert }\int_{\omega
_{0}\setminus(  \mathbb{B}_{\varrho}(  P^{j})  \cup
\mathbb{B}_{\varrho}(  P^{k})  )  }(  G_{j}(
y)  -G_{k}(  y)  )  dy=\label{2.20}\\
&  =\lim_{\varrho\rightarrow+0}\int_{\omega_{0}\setminus(  \mathbb{B}%
_{\varrho}(  P^{j})  \cup\mathbb{B}_{\varrho}(  P^{k})
)  }(  G_{j}(  y)  \Delta_{y}G_{k}(  y)
-G_{k}(  y)  \Delta_{y}G_{j}(  y)  )
dy=\nonumber\\
&  =\lim_{\varrho\rightarrow+0}\int_{\partial\mathbb{B}_{\varrho}(
P^{j})  \cup\partial\mathbb{B}_{\varrho}(  P^{k})  }(
G_{k}(  y)  \frac{\partial G_{j}}{\partial r}(  y)
-G_{j}(  y)  \frac{\partial G_{k}}{\partial r}(  y)
)  ds_{y}=\nonumber\\
&  =\lim_{\varrho\rightarrow+0}(  G_{k}(  P^{j})
-G_{j}(  P^{k})  )  \varrho\int_{0}^{2\pi}\frac{1}{2\pi}%
\frac{\partial}{\partial r}\left.\ln\frac{1}{r}\right|_{r=\varrho}d\varphi=G_{j}(
P^{k})  -G_{k}(  P^{j})  .\nonumber
\end{align}
Here, $j\neq k,$ $\mathbb{B}_{\varrho}(  P^{j})  =\{
y:r<\varrho\}  ,$ and $\varphi$ is the angular variable.

In what follows it is convenient to restrict the equation (\ref{2.13}) onto
the punctured domain%
\begin{equation}
\omega_{\odot}=\omega_{0}\setminus\{  P^{1},...,P^{J}\}  ,
\label{2.00}%
\end{equation}
that is%
\begin{equation}
-\Delta_{y}U_{0}^{0}(  y)  =f_{0}^{0}(  y)
,\ \ \ y\in\omega_{\odot}. \label{2.21}%
\end{equation}
In this way all Green functions $G_{j}$ become singular solutions of the
problem (\ref{2.21}), (\ref{2.14}) with the constant right-hand side
$f_{0}^{0}(  y)  =-\vert \omega_{0}\vert ^{-1}.$ Notice
that $G_{j}\in L^{2}(  \omega_{0})  $ due to (\ref{2.18}).

\begin{proposition}
\label{prop2.1}

\begin{enumerate}
\item If $f_{0}^{0}\in L^{2}(  \omega_{0})  $ satisfies $\langle f_{0}^{0}\rangle_0=0$, % \label{2.22}%
the problem (\ref{2.13}), (\ref{2.14}) admits a solution in $H^{2}(
\omega_{0})  .$ This solution $U_{0}^{0}$ is defined up to an additive
constant and under the orthogonality condition $\langle U_{0}^{0}\rangle_0=0$,
it becomes unique and meets the estimate $\Vert U_{0}^{0};H^{2}(  \omega_{0})  \Vert \leq
c_{0}\Vert f_{0}^{0};L^{2}(  \omega_{0})  \Vert .$

\item If $f_{0}^{0}\in L^{2}(  \omega_{0})  $ satisfies the orthogonality condition $\langle f_{0}^{0}\rangle_0=0$, the problem (\ref{2.21}), (\ref{2.14}) has a solution $U_{0}%
^{0}\in H_{loc}^{2}(  \overline{\omega_{0}}\setminus\{
P^{1},...,P^{J}\}  )  \cap L^{2}(  \omega_{0})  $ in
the form%
\begin{equation}
U_{0}^{0}(  y)  =U_{\bot}^{0}(  y)  +A_{0}+{\textstyle\sum\nolimits_{j}}A_{j}G_{j}(  y)  \label{2.25}%
\end{equation}
where $A_{0},$ $A_{1},...,$ $A_{J}$ are arbitrary constants subject to the
relation%
\begin{equation}
A_{1}+...+A_{J}=-\langle f_{0}^{0}\rangle_0=0 \label{2.26}%
\end{equation}
and $U_{\bot}^{0}\in H^{2}(  \omega_{0})  $ is the solution of the
problem (\ref{2.13}), (\ref{2.14}) with the orthogonality condition
$\langle U_{0}^{0}\rangle_0=0$ and the new right-hand side%
\begin{equation}
f_{\bot}^{0}(  y)  =f_{0}^{0}(  y)  -\vert
\omega_{0}\vert^{-1}\langle f_{0}^{0}\rangle_0.
\label{2.27}%
\end{equation}
There holds the estimate%
\begin{equation}
\Vert U_{\bot}^{0};H^{2}(  \omega_{0})  \Vert \leq
c\Vert f_{0}^{0};L^{2}(  \omega_{0})  \Vert .
\label{2.bot}%
\end{equation}

\item Any solution $U_{0}^{0}\in L^{2}(  \omega_{0})  \cap
H_{loc}^{2}(  \overline{\omega_{0}}\setminus\{  P^{1},...,P^{J}%
\}  )  $ of the problem (\ref{2.21}), (\ref{2.14}) gets the form
(\ref{2.25}) with coefficients as in (\ref{2.26}).
\end{enumerate}
\end{proposition}

\textbf{Proof}. The first assertion is well known (cf. \cite{LiMa}) and the
second one follows directly from relations (\ref{2.16}) and the fact that the
function (\ref{2.27}) meets the orthogonality condition $\langle f_{0}^{0}\rangle_0=0$. To
confirm the last assertion, we observe that $U_{0}^{0}$ is a distributional
solution of the problem in the intact domain $\omega_{0}$ while the right-hand
side of the Poisson equation is a sum of $f_{0}^{0}\in L^{2}(  \omega
_{0})  $ and a linear combination of the Dirac functions $\delta(
y-P^{j})  $ and their derivatives (the latter is due to the theorem on a
distribution with a point support; cf. \cite[\S 1.2.6]{Vlad}). It remains to
mention that derivatives of the Green function $G(  y,P)  $ in the
second argument live outside $L^{2}(  \omega_{0})  $ but $G$ itself
falls into $L^{2}(  \omega_{0})  .$ For the representation
(\ref{2.25}), we also refer to the paper \cite{na101} where much more general
situation was considered on the base of the Kondratiev theory \cite{Ko}.
$\blacksquare$

\bigskip

Notice that, owing to the orthogonality condition in (\ref{2.16}), we have%
\begin{equation}
U_{\bot}^{0}(  P^{j})  =\int_{\omega_{0}}G_{j}(  y)
f_{0}^{0}(  y)  dy. \label{2.BOT}%
\end{equation}

\subsection{The limit problem in the junction of a layer and a
semi-cylinder\label{sect2.3}}

\begin{figure}
[ptb]
\begin{center}
\includegraphics[
height=1.0265in,
width=2.6662in
]
{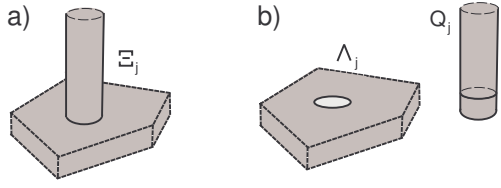}%
\caption{The union of the layer and semi-cylinder (a). The layer and the
semi-cylinder (b). }%
\label{f3}%
\end{center}
\end{figure}

Here, we consider the case $\alpha= 0$. In order to describe the boundary layer phenomenon in the vicinity of the
points $P^{1},...,P^{J}$ we need to examine solutions of the boundary value
problem%
\begin{align}
-\Delta_{\xi}w_{0}^{j}(  \xi)   &  =\mathcal{F}_{0}^{j}(
\xi)  ,\ \ \ \xi\in\Lambda_{j},\label{2.28}\\
-\gamma_{j}\Delta_{\xi}w_{j}^{j}(  \xi)   &  =\mathcal{F}_{j}%
^{j}(  \xi)  ,\ \ \ \xi\in Q_{j},\label{2.29}\\
-\partial_{\zeta}w_{0}^{j}(  \eta,0)   &  =\partial_{\zeta}%
w_{0}^{j}(  \eta,1)  =0,\ \ \ \eta\in\mathbb{R}^{2}\setminus
\overline{\omega_{j}},\label{2.30}\\
-\gamma_{j}\partial_{\zeta}w_{j}^{j}(  \eta,0)   &  =0,\ \ \ \eta
\in\omega_{j},\label{2.31}\\
-\gamma_{j}\partial_{\nu}w_{j}^{j}(  \xi)   &  =0,\ \ \ \xi
\in\partial\omega_{j}\times[  1,+\infty)  ,\label{2.32}\\
w_{0}^{j}(  \xi)   &  =w_{0j}^{j}(  \xi)
,\ \ \ \partial_{\nu}w_{0}^{j}(  \xi)  =\gamma_{j}\partial_{\nu
}w_{j}^{j}(  \xi)  ,\ \ \ \xi\in\upsilon_{j}^{1}=\partial\omega
_{j}\times(  0,1)  \label{2.33}%
\end{align}
in the junction, fig. \ref{f3},a,%
\begin{equation}
\Xi_{j}=\Lambda_{j}\cup Q_{j} \label{2.34}%
\end{equation}
of the perforated layer and semi-cylinder%
\begin{equation}
\Lambda_{j}=\{  \xi=(  \eta,\zeta)  :\eta\in\mathbb{R}%
^{2}\setminus\overline{\omega_{j}},\ \zeta\in(  0,1)  \}
,\ \ \ Q_{j}=\{  \xi=(  \eta,\zeta)  :\eta\in\omega
_{j},\ \zeta>0\}  . \label{2.35}%
\end{equation}%
%TCIMACRO{\FRAME{ftbpFU}{2.412in}{1.4633in}{0pt}{\Qcb{The union of the layer
%and semi-cylinder (a). The layer and the semi-cylinder (b). }}{\Qlb{f3}%
%}{junction3new.eps}{\special{ language "Scientific Word";  type "GRAPHIC";
%width 2.412in;  height 1.4633in;  depth 0pt;  original-width 8.003in;
%original-height 4.8395in;  cropleft "0";  croptop "1";  cropright "1";
%cropbottom "0";  filename 'Junction3new.eps';file-properties "XNPEU";}} }%
%BeginExpansion
%EndExpansion

The set (\ref{2.34}) is obtained from the set (\ref{1.4}) by going over to the
stretched coordinates%
\begin{equation}
\xi^{j}=(  \eta^{j},\zeta)  ,\ \ \ \eta^{j}=(  \eta_{1}%
^{j},\eta_{2}^{j})  =h^{-1}(  y-P^{j})  ,\ \ \ \zeta=h^{-1}z,
\label{2.36}%
\end{equation}
see (\ref{1.2}) and (\ref{1.1}), and putting $h=0$ formally. Note that in
the equations (\ref{2.28})-(\ref{2.33}) we do not mark the variables $\xi$ and
$\eta$ with the superscript $j.$ In the sequel we also omit this superscript
on the functions $w$ and $\mathcal{F}$ whose restrictions on $\Lambda_{j}$ and
$Q_{j}$ are denoted by $w_{0},$ $\mathcal{F}_{0}$ and $w_{j},$ $\mathcal{F}%
_{j},$ respectively. Moreover, let us notice that the domain $Q_{j}$ is a
semi-infinite cylinder; for simplicity, we call it the semi-cylinder.

The equations (\ref{2.28}), (\ref{2.29}) and the boundary
conditions (\ref{2.30})-(\ref{2.32}) are obtained from (\ref{1.10}),
(\ref{1.11}) and (\ref{1.12}), (\ref{1.13}), respectively, by the change
$x\mapsto\xi.$ The transmission conditions (\ref{2.33}) at the surface
$\upsilon_{j}^{1}=\partial\Lambda_{j}\cup\partial Q_{j}$ result from the
original conditions (\ref{1.15}), (\ref{1.16}) with the exponent
$\alpha=0$ in (\ref{1.17}).
%We emphasize that, for $\alpha=1,$ the conditions
%(\ref{1.15}), (\ref{1.16}) decouple and the corresponding limit problems will
%be considered in the next section.
In view of a simple geometry the problem (\ref{2.28})-(\ref{2.33}) can be
investigated by the Fourier method (see, e.g., \cite{Smir}). However, thinking
about possible generalizations, see Section \ref{sect5}, we here realize a
different approach.

The variational formulation of the problem (\ref{2.28})-(\ref{2.33}) appeals
to the integral identity%
\begin{equation}
(  \nabla_{\xi}w_{0},\nabla_{\xi}v_{0})  _{\Lambda_{j}}+\gamma
_{j}(  \nabla_{\xi}w_{j},\nabla_{\xi}v_{j})  _{Q_{j}}=(
\mathcal{F},v)  _{\Xi_{j}}\ \ \ \ \ \ \ \forall v\in\mathcal{H}_{j},
\label{2.37}%
\end{equation}
where $\nabla_{\xi}=\operatorname{grad},$ $(  \ ,\ )  _{\Xi}$ is
the natural scalar product in the Lebesgue space $L^{2}(  \Xi)  ,$
and $\mathcal{H}_{j}$ is the completion of $C_{c}^{\infty}(
\overline{\Xi_{j}})  $ (infinitely differentiable functions with compact
supports) in the norm%
\begin{equation}
\Vert w;\mathcal{H}_{j}\Vert =(  \Vert \nabla_{\xi
}w;L^{2}(  \Xi_{j})  \Vert ^{2}+\Vert w;L^{2}(
\Theta_{j})  \Vert ^{2})  ^{1/2} \label{2.38}%
\end{equation}
where $\Theta_{j}=\omega_{j}\times(  0,1)  .$ Notice that a constant belongs to $\mathcal{H}_{j}$ (see, e.g., \cite{na164, na220}).

%\begin{remark}
%\label{rem2.2}The norm $\Vert w;L^{2}(  \Theta_{j})
%\Vert $ in the finite domain $\Theta_{j}$ cannot be taken away from the
%definition (\ref{2.38}) because $\Vert \nabla_{\xi}w;L^{2}(  \Xi
%_{j})  \Vert $ implies only a semi-norm in the Hilbert space
%$\mathcal{H}_{j}.$ Indeed, one may verify directly that a constant function
%belongs to $\mathcal{H}_{j}$ since the Dirichlet integral
%\[
%\int_{\Xi_{j}}\vert \nabla_{\xi}X_{T}(  \xi)  \vert
%^{2}d\xi
%\]
%calculated for the function $X_{T}$ which approximates $1,$ has a compact
%support and is given by%
%\begin{equation}
%X_{T}(  \eta,\zeta)  =\{
%\begin{array}
%[c]{c}%
%\chi(  \zeta-N)  ,\ \ \ \xi\in Q_{j},\\
%\chi(  (  \ln T)  ^{-1}\ln\vert \eta\vert )
%,\ \ \ \xi\in\Lambda_{j},
%\end{array}
%.  \ \ \ \ \ \chi(  t)  =\{
%\begin{array}
%[c]{c}%
%1,\ \ \ t<1,\\
%0,\ \ \ t>2,
%\end{array}
%.  \label{2.39}%
%\end{equation}
%becomes infinitesimal as $T\rightarrow+\infty.$ $\blacksquare$
%\end{remark}

\begin{lemma}
\label{lem2.3} An equivalent norm in $\mathcal{H}_{j}$ can be chosen as
$(  \Vert \nabla_{\xi}w;L^{2}(  \Xi_{j})  \Vert
^{2}+\Vert \mathcal{R}^{-1}w;L^{2}(  \Xi_{j})  \Vert
^{2})  ^{1/2} $ where $\mathcal{R}$\ is a smooth positive function in $\overline{\Xi_{j}}$
such that%
\begin{equation}
\mathcal{R}(  \xi)   =\zeta\text{ \ for \ }\xi\in
Q_{j},\ \zeta>2, \qquad \mathcal{R}(  \xi)    =\rho\ln\rho\text{ \ for \ }\xi
\in\Lambda_{j},\ \rho=\vert \eta\vert >2.\label{2.41}
\end{equation}

\end{lemma}

\textbf{Proof.} First of all, we observe that the cylinder $\Theta_{j}%
=\omega_{j}\times(  0,1)  $ can be obviously replaced in (\ref{2.38}) by any
bigger domain $\Theta_{j}^{\prime}$ with a compact closure.
%This follows from
%the obvious inequality%
%\begin{equation}
%\Vert w;L^{2}(  \Theta_{j}^{\prime})  \Vert \leq
%C(  \Theta_{j},\Theta_{j}^{\prime})  (  \Vert \nabla
%_{\xi}w;L^{2}(  \Theta_{j}^{\prime})  \Vert +\Vert
%w;L^{2}(  \Theta_{j})  \Vert )  . \label{2.411}%
%\end{equation}
Next, we recall two one-dimensional Hardy inequalities%
\begin{align}
\int_{0}^{L}\zeta^{-2}\vert W(  \zeta)  \vert
^{2}d\zeta &  \leq4\int_{0}^{L}\left\vert \frac{\partial W}%
{\partial\zeta}(  \zeta)  \right\vert ^{2}d\zeta\ \ \ \ \forall
W\in C_{c}^{1}(  0,L]  , \quad L\in (0,+\infty],\label{2.42}\\
\int_{1}^{+\infty}\rho^{-1}\vert \ln\rho\vert ^{-2}\vert
W(  \rho)  \vert ^{2}d\rho &  \leq4\int_{1}^{+\infty}%
\rho\left\vert \frac{\partial W}{\partial\rho}(  \rho)  \right\vert
^{2}d\rho\ \ \ \ \forall W\in C_{c}^{1}(  1,+\infty)  .
\label{2.43}%
\end{align}
Note that (\ref{2.43}) is derived from (\ref{2.42}) by the change $\rho\mapsto
t=\ln\rho.$
%Moreover, for reader's convenience, we demonstrate a verification
%of a simple analog of (\ref{2.42}) in calculation (\ref{4.0}).

Finally, we integrate inequality (\ref{2.42}) with $W(  \zeta)
=(  1-\chi(  \zeta)  )  w(  \eta,\zeta)  $
over $\omega_{j}\ni\eta$ and inequality (\ref{2.43}) with $L=\infty $ and $W(
\rho)  =(  1-\chi(  \rho)  )  w(  \eta
,\zeta)  $ in the angular variable $\varphi\in(  0,2\pi)$,
where $\chi(t)=1$, if $t<1$, and $\chi(t)=0$, if $t>2$.
As a result we obtain the estimates
\begin{align}
\Vert \zeta^{-2}w;L^{2}(  \omega_{j}\times(  2,+\infty)
)  \Vert ^{2}\leq C(  \Vert \partial_{\zeta}%
w;L^{2}(  \omega_{j}\times(  1,+\infty))\Vert
^{2}&+\Vert w;L^{2}(  \omega_{j}\times(  1,2))
\Vert ^{2})  ,\label{2.44}\\
\Vert \rho^{-1}\vert \ln\rho\vert ^{-1}w;L^{2}((
\mathbb{R}^{2}\setminus \mathbb{B}_{2})  \times(
0,1) )  \Vert ^{2}  &  \leq C(\Vert
\partial_{\rho}w;L^{2}(  (  \mathbb{R}^{2}\setminus \mathbb{B}_{1})  \times(  0,1)  ) \Vert ^{2}+\nonumber
\\
&  +\Vert w;L^{2}((  \mathbb{B}_{2}  \setminus
\mathbb{B}_{1})  \times(  0,1))\Vert
^{2})  .\nonumber
\end{align}
where here and in the sequel $\mathbb{B}_{R}=\{\eta\in\mathbb{R}^{2}:\vert \eta\vert<R\}$ denotes the ball of radius $R$. We mention that cut-off functions were introduced in order to fulfil the
conditions $W(  0)  =0$ in (\ref{2.42}) and $W(  1)  =0$
in (\ref{2.43}). Furthermore, we have used the relations%
\begin{align*}
\partial_{\zeta}(  (  1-\chi(  \zeta)  )  w(
\eta,\zeta)  )   &  =(  1-\chi(  \zeta)  )
\partial_{\zeta}w(  \eta,\zeta)  -w(  \eta,\zeta)
\partial_{\zeta}\chi(  \zeta)  ,\\
1-\chi(  \zeta)   &  =1\text{ \ for \ }\zeta>2\text{ \ and
\ }\partial_{\zeta}\chi(  \zeta)  =0\text{ \ for \ }\zeta
\notin(  1,2)
\end{align*}
to get the first estimate and similar relations together with the formula
$d\eta=\rho d\rho d\varphi$ for the second estimate in (\ref{2.44}). Choosing $\Theta_{j}^{\prime}=\{  \xi\in\Xi_{j}:\rho<\max\{
2,\text{diam }\omega_{j}\}  ,\ \zeta\in(  0,2)  \}  $
 and according to (\ref{2.38}), (\ref{2.41}), we get $\Vert \mathcal{R}^{-1}w;L^{2}(  \Xi_{j})  \Vert ^{2}\leq
c\Vert w;\mathcal{H}_{j}\Vert ^{2}$ which suffices to conclude with the proof. $\blacksquare$

By a standard argument, the Riesz representation theorem leads to the following assertion.
\begin{proposition}
\label{prop2.4}Let the right-hand side $\mathcal{F}$ satisfy the conditions
$\mathcal{RF}\in L^{2}(  \Xi_{j})  $ and
\begin{equation}
\int\nolimits_{\Xi_{j}} \mathcal{F}(  \xi)  d\xi=0. \label{2.45}%
\end{equation}
Then the integral identity (\ref{2.37}) (problem (\ref{2.28})-(\ref{2.33}) in
the differential form) has a solution $w\in\mathcal{H}_{j}$ defined up to an
additive constant. Under the orthogonality condition $\int_{\Theta_{j}}w(  \xi)  d\xi=0$
the solution becomes unique and meets the estimate
$
\Vert w;\mathcal{H}_{j}\Vert ^{2}\leq c\Vert \mathcal{RF}%
;L^{2}(  \Xi_{j})  \Vert ^{2}.
$
\end{proposition}
%\textbf{Proof.} Let us consider the auxiliary problem%
%\begin{equation}
%(  \nabla_{\xi}w_{0},\nabla_{\xi}v_{0})  _{\Lambda_{j}}+\gamma
%_{j}(  \nabla_{\xi}w_{j},\nabla_{\xi}v_{j})  _{Q_{j}}+\mu(
%w,v)  _{\Theta_{j}}=(  \mathcal{F},v)  _{\Xi_{j}%
%}\ \ \ \ \ \ \ \forall v\in\mathcal{H}_{j}. \label{2.48}%
%\end{equation}
%One readily sees that, for $\mu>0,$ the left-hand side of (\ref{2.48}) becomes
%a scalar product in the Hilbert space $\mathcal{H}_{j}.$ Moreover, due to
%Lemma \ref{lem2.3} the right-hand side is a continuous functional in
%$v\in\mathcal{H}_{j}.$ Hence, the Riesz representation theorem delivers a
%unique solution of the problem (\ref{2.48}) with $\mu>0.$ The additional term
%$\mu(  w,v)  _{\Theta_{j}}$ generates a compact perturbation of the
%operator of the original problem (\ref{2.37}) because the domain $\Theta_{j}$
%is finite and therefore the embedding $\mathcal{H}_{j}\subset L^{2}(
%\Theta_{j})  $ is compact. We apply the Fredholm alternative and readily
%deduce that any solution $w\in\mathcal{H}_{j}$ of the homogeneous
%($\mathcal{F}=0$) problem (\ref{2.37}) is constant. Thus, (\ref{2.45}) is the
%only defect functional and the equality (\ref{2.46}) fixes the solution
%uniquely. $\blacksquare$
%
%\bigskip
%
In addition to a constant solution of the homogeneous problem (\ref{2.28}%
)-(\ref{2.33}) we shall need an unbounded solution $\mathbf{w}^{j}$ defined
uniquely by its asymptotic forms%
\begin{align}
\mathbf{w}^{j}(  \xi)   &  =-(2\pi)^{-1}\ln\rho + o(  1)  ,\ \ \ \xi\in\Lambda_{j},\ \ \rho=\vert
\eta\vert \rightarrow+\infty,\label{2.50}\\
\mathbf{w}^{j}(  \xi)   &  =\gamma_{j}^{-1}\vert \omega
_{j}\vert ^{-1}\zeta+\mathbf{q}_{j}+o(  1)  ,\ \ \ \xi\in
Q_{j},\ \ \zeta\rightarrow+\infty, \label{2.49}%
\end{align}
where $\mathbf{q}_{j}$ is a constant described below. A distinct way to find
out this solution is to search for the remainder $\widehat{\mathbf{w}}^{j}%
\in\mathcal{H}_{j}$ in the representation%
\begin{equation}
\mathbf{w}^{j}(  \xi)  =\gamma_j^{-1}\vert\omega_j\vert^{-1} X_{Q}(  \zeta)  -(2\pi)^{-1} X_{\Lambda}(  \eta)\ln\rho
+\widehat{\mathbf{w}}_{j}(  \xi)
\nonumber%
\end{equation}
where $X_{Q}(  \zeta)  =1-\chi(  \zeta)  $ and
$X_{\Lambda}(  \eta)  =1-\chi(  \rho/R_{j})  $ are
smooth cut-off functions,
\begin{align}
X_{Q}(  \zeta)   &  =0\text{ \ for \ }\zeta<1,\ \ X_{Q}(
\zeta)  =1\text{ \ for \ }\zeta>2,\label{QL}\\
X_{\Lambda}(  \eta)   &  =0\text{ \ for \ }\vert
\eta\vert <R_{j},\ \ X_{\Lambda}(  \eta)  =1\text{ \ for
\ }\vert \eta\vert >2R_{j},\nonumber
\end{align}
and radius $R_{j}>0$ is fixed such that $\overline{\omega_{j}}\subset\{
\eta:\vert \eta\vert <R_{j}\}  .$ The remainder must satisfy
problem (\ref{2.28})-(\ref{2.33}) \ with the smooth and compactly supported
right-hand sides
$
\mathcal{F}_{0}^{j}=-(2\pi)^{-1}[  \Delta_{\xi},X_{\Lambda}] \ln\rho,\ \mathcal{F}_{j}^{j}=\vert\omega_j\vert^{-1}[  \Delta_{\xi},X_{Q}]\zeta,
$
where $[  \Delta_{\xi},X]  $ is the commutator of the Laplace
operator and a cut-off function $X,$%
\begin{equation}
[  \Delta_{\xi},X]  w=2\nabla_{\xi}w\cdot\nabla_{\xi}X+w\Delta
_{\xi}X. \label{2.com}%
\end{equation}
Clearly, $\mathcal{RF}\in L^{2}(  \Xi_{j})  $ and, to get
$\widehat{\mathbf{w}}_{j},$ we only need to verify the validity of the
orthogonality condition (\ref{2.45}) in Proposition \ref{prop2.4}.
It is assured by the following calculation:
\begin{align}
&  \frac{1}{2\pi}\int_{\Lambda_{j}}[  \Delta_{\xi},X_{\Lambda}]
\ln\frac{1}{\rho}d\xi+\frac{1}{\vert \omega_{j}\vert }\int_{Q_{j}%
}[  \Delta_{\xi},X_{Q}]  \zeta d\xi\label{2.52}\\
&  =\frac{1}{2\pi}\lim_{T\rightarrow+\infty}\int_{0}^{1}\int_{\mathbb{B}%
_{T}\setminus\mathbb{B}_{R_{j}}}\Delta_{\xi}\left(  X_{\Lambda}(
\eta)  \ln\frac{1}{\rho}\right)  d\eta d\zeta+\frac{1}{\vert
\omega_{j}\vert }\lim_{T\rightarrow+\infty}\int_{\omega_{j}}\int_{1}%
^{T}\Delta_{\xi}(  X_{Q}(  \zeta)  \zeta)  d\zeta
d\eta\nonumber\\
&  =\frac{1}{2\pi}\lim_{T\rightarrow+\infty}\int_{\partial\mathbb{B}_{T}}%
\left.\frac{\partial}{\partial\rho}\ln\frac{1}{\rho}\right|_{\rho=T}~ds_{\eta}+\frac
{1}{\vert \omega_{j}\vert }\lim_{T\rightarrow+\infty}\int
_{\omega_{j}}\left.\frac{\partial\zeta}{\partial\zeta}\right|_{\zeta=T}~d\eta
=-1+1=0.\nonumber
\end{align}
The decompositions (\ref{2.50}) and (\ref{2.49}) are supported
by the Fourier method where $\mathbf{q}_{j}$ is a constant which is defined uniquely and depends on
$\gamma_{j}$ and $\omega_{j}$ (cf. Remark \ref{rem2.q}).

We finally formulate a representation of a solution to the homogeneous problem (\ref{2.28})-(\ref{2.33}) which, for example, follows from our consideration of the decoupled problems in the next section.

\begin{theorem}
\label{th2.5}Any solution $w\in H_{loc}^{1}(  \Xi_{j})  $ of the
homogeneous problem (\ref{2.28})-(\ref{2.33}) satisfying the estimates%
\begin{equation}
\vert w_{0}(  \xi)  \vert \leq c\ln(
1+\rho)  \text{\ in }\Lambda_{j},\ \rho>R_{j},\quad \vert w_{j}(  \xi)  \vert \leq c\,\zeta\text{ \ in
}Q_{j},\ \zeta\geq 2,\label{2.53}
\end{equation}
is a linear combination $c_{0}+c_{1}\mathbf{w}^{j}$ with some
coefficients $c_{i}$ and the constructed special solution $\mathbf{w}^{j}$.
\end{theorem}

\subsection{The limit problems in a perforated layer and in a
semi-cylinder\label{sect2.4}}

In the case $\alpha=1,$ equation (\ref{1.16}) contains the big factor $h^{-1}$
and, therefore, after the coordinate dilation $x\mapsto\xi^{j},$ see
(\ref{2.36}), the transmission conditions (\ref{1.15}), (\ref{1.16}) decouple
while junction (\ref{2.34}) splits into the perforated layer and the
semi-cylinder (\ref{2.35}), see fig. \ref{f3},b. We shall see in Section
\ref{sect3.1} that the transmission conditions give rise to the Dirichlet
condition on the hole surface $\upsilon_{j}^{1}\subset\partial\Lambda_{j}$ and
to the Neumann condition on the ring $\upsilon_{j}^{1}\subset\partial Q_{j}$
near the cylinder end. In this way, to describe the boundary layer phenomenon, we have to deal with two
problems, namely the mixed boundary value problem%
\begin{align}
-\Delta_{\xi}W_{0}(  \xi)   &  =0,\ \xi\in\Lambda_{j}%
,\ \ \partial_{\zeta}W_{0}(  \eta,0)  =\partial_{\zeta}W_{0}(
\eta,1)  =0,\ \ \ \eta\in\mathbb{R}^{2}\setminus\overline{\omega_{j}%
},\label{2.57}\\
W_{0}(  \xi)   &  =g_{0}(  \xi)  ,\ \ \xi\in\upsilon
_{j}^{1},\nonumber
\end{align}
in the perforated layer and the Neumann problem in the semi-cylinder%
\begin{align}
-\gamma_{j}\Delta_{\xi}W_{0}(  \xi)   &  =0,\ \xi\in Q_{j}%
,\ \ \gamma_{j}\partial_{\nu}W_{0}(  \eta,\zeta)  =g_{j}(
\eta,\zeta)  ,\ \ \ (  \eta,\zeta)  \in\partial\omega
_{j}\times\mathbb{R}_{+},\label{2.58}\\
\gamma_{j}\partial_{\zeta}W_{0}(  \eta,0)   &  =0,\ \ \eta\in
\omega_{j}. \label{2.59}%
\end{align}
Both the problems permit separation of variables and are rather standard in
the asymptotic analysis of rods and perforated plates, even isolated (cf.,
respectively, the monographs \cite{MaNaPl, KoMaMo, Panas, Nabook} and the
papers \cite{na164, na472}). We here present only some comprehensible pieces
of information on them which will be used for asymptotic structures in Section
\ref{sect3}.

First of all, a solution of the homogeneous problem (\ref{2.57}) with the
logarithmical growth at infinity, cf. (\ref{2.53}), does not depend on the variable $\zeta$ and takes the form of the
\textit{logarithmic potential} $\mathbf{W}_{j}(  \eta)  $ that is
a harmonic function in $\mathbb{R}^{2}\setminus\overline
{\omega_{j}}$ which vanishes at $\partial\omega_{j}$ and admits the
representation%
\begin{equation}
\mathbf{W}_{j}(  \eta)  =(2\pi)^{-1}(  -\ln\rho +\ln c_{\log}(\omega_{j}))  +O(\rho^{-1})  ,\ \ \rho\rightarrow+\infty, \label{2.60}%
\end{equation}
where $c_{\log}(  \omega_{j})  $ is the \textit{logarithmic
capacity} of the set $\overline{\omega_{j}}$ (see, e.g., \cite{PoSe, Land}).

\begin{remark}
\label{rem2.q} Owing to (\ref{2.49}) and (\ref{2.50}), the difference
$\mathbf{w}^{j}(  \xi)  -\mathbf{q}^{j}$ also can be represented in
form (\ref{2.60}) inside the layer $\Lambda_{j}$ as $-(  2\pi)  ^{-1}\ln\rho-\mathbf{q}^{j}+O(  \rho^{-1})
, \rho\rightarrow+\infty,$ so that the quantity $e^{-2\pi\mathbf{q}^{j}}$ could be called the logarithmic
capacity in the junction of the layer and semi-cylinder (\ref{2.34}).
\end{remark}

\begin{lemma}
\label{lem2.6}There holds the equality $%
%TCIMACRO{\dint _{\partial\omega_{j}}}%
%BeginExpansion
{\displaystyle\int_{\partial\omega_{j}}}
%EndExpansion
\partial_{\nu}\mathbf{W}_{j}(  \eta)  ds_{\eta}=1.$
\end{lemma}

\textbf{Proof.} We just repeat the calculation (\ref{2.52}):%
\[
\int_{\partial\omega_{j}}\partial_{\nu}\mathbf{W}_{j}(  \eta)
ds_{\eta}=-\lim_{T\rightarrow+\infty}\int_{\partial\mathbb{B}_{T}}%
\frac{\partial\mathbf{W}_{j}}{\partial\rho}(  \eta)  ds_{\eta
}=-\frac{1}{2\pi}\lim_{T\rightarrow+\infty}\int_{\partial\mathbb{B}_{T}}%
\frac{\partial}{\partial\rho}\ln\frac{1}{\rho}ds_{\eta}=1.\ \ \ \ \blacksquare
\]

As for the Neumann problem in the semi-cylinder $Q_{j},$ we will use the
following assertion which is based on the Fourier method.

\begin{lemma}
\label{lem2.7}Let the right-hand side $g_{j}\in L^{2}(  \partial
\omega_{j}\times\mathbb{R}_{+})  $ in problem (\ref{2.58}), (\ref{2.59})
possess a compact support. Then the problem has a solution $W_{j}$ which is
defined up to an additive constant and gets the form%
\begin{equation}
W_{j}(  \eta,\zeta)  =C_{j}\zeta+C_{0j}+O(  e^{-\delta\zeta
})  ,\ \ \zeta\rightarrow+\infty, \label{2.61}%
\end{equation}
where $C_{0j}$ is a constant, $\delta$ is some positive number and%
\begin{equation}
C_{j}=-\frac{1}{\gamma_{j}\vert \omega_{j}\vert }\int_{0}^{+\infty
}\int_{\partial\omega_{j}}g_{j}(  \eta,\zeta)  ds_{\eta}d\zeta.
\nonumber%
\end{equation}

\end{lemma}

%\textbf{Proof.} The existence of a solution in form (\ref{2.61}) can be
%demonstrated by the Fourier method (see, e.g., \cite{Smir} and also \cite[Ch.2
%and 5]{NaPl} for another approach based on the Kondratiev theory \cite{Ko}).
%Similarly to (\ref{2.52}), the Green formula in the truncated cylinder
%$\omega_{j}\times(  0,T)  $ yields
%\[
%\int_{0}^{T}\int_{\partial\omega_{j}}g_{j}(  \eta,\zeta)  ds_{\eta
%}d\zeta=-\gamma_{j}\int_{\omega_{j}}\frac{\partial W_{j}}{\partial\zeta
%}(  \eta,T)  d\eta=-\gamma_{j}C_{j}\vert \omega_{j}\vert
%+O(  e^{-\delta T})  .
%\]
%The limit passage $T\rightarrow+\infty$ completes the proof. $\blacksquare$

\begin{remark}
\label{rem2.8}To describe the boundary layer effect at the rod "soles"
$\omega_{j}^{h}(  l_{j})  $ with the Dirichlet conditions
(\ref{1.14}), one has to consider the mixed boundary value problem in the
semi-cylinder $Q_{j},$ too, which now is obtained as a result of the
coordinate dilation $x\mapsto(  h^{-1}(  y-P^{j})
,h^{-1}(  l_{j}-z)  )  ,$ cf. (\ref{2.36}). This problem
consists of equations (\ref{2.58}) and the Dirichlet condition at the cylinder end
\begin{equation}
W_{j}(  \eta,0)  =g_{j}^{l}(  \eta)  ,\ \ \eta\in
\omega_{j}. \label{2.63}%
\end{equation}
Although we do not involve solutions of (\ref{2.58}), (\ref{2.63}) into
asymptotic structures in Section \ref{sect3}, we again mention the monographs
\cite{MaNaPl, KoMaMo, Panas, Nabook} where this problem was investigated and
applied. $\blacksquare$
\end{remark}

\section{Constructing asymptotic expansions\label{sect3}}

Both cases $\alpha=1$ and $\alpha=0$ are investigated by means of the method of
matched asymptotic expansions, see, e.g, the monographs \cite{VanDyke,
Ilinbook} and \cite[Ch.2]{MaNaPl}, so that the asymptotic procedures for the
contrasting and homogeneous junctions defined according to (\ref{1.17}),
(\ref{1.18}) look quite similar although provide very different formulas for
the main asymptotic terms in the solution $u(  h,x)  $ of problem
(\ref{1.10})-(\ref{1.16}). An evident distinction of final formulas originates
in the different structure of the inner expansions composed from solution to
the limit problems (\ref{2.57})-(\ref{2.59}) and (\ref{2.28})-(\ref{2.33})
posed on infinite domains depicted in fig. \ref{f3}, a and b, respectively. In
other words, the asymptotic structure of solutions in the junction (\ref{1.4})
is crucially prescribed by boundary layer effects around the junction zones.

\subsection{The case $\alpha=1.$\label{sect3.1}}

Based on our preliminary consideration in Section \ref{sect2}, we assume the
asymptotic ans\"{a}tze (\ref{2.3}) in the rods $\Omega_{j}(  h)  $
and (\ref{2.10}) in the plate $\Omega_{\bullet}(  h)  .$ The main term $U_{j}^{0}$ in (\ref{2.3}) satisfies the ordinary differential
equation (\ref{2.6}) and the Dirichlet condition (\ref{2.7}), however a
condition at $z=0$ is still absent and needs to be derived. The main term $U_{0}^{0}$ in (\ref{2.10}) is a singular solution of the
Neumann problem (\ref{2.21}), (\ref{2.14}) in the punctured domain
(\ref{2.00}) while coefficients in its representation (\ref{2.25}) also remain
unknown. In order to find out appropriate values of the coefficients
$A_{0},A_{1},...,A_{J}$ as well as right-hand sides in
\begin{equation}
-\gamma_{j}\vert \omega_{j}\vert \partial_{z}U_{j}^{0}(
0)  =G_{j}^{0} \label{3.1}%
\end{equation}
we construct inner asymptotic expansions in the vicinity of the points
$P^{1},...,P^{J}.$ Namely, in a neighborhood $V_{j}$ of $P^{j}$ where the
stretched coordinates (\ref{2.36}) were defined, we search for the asymptotic
forms%
\begin{align}
u(  h,x)   &  =W_{j}^{0}(  \eta^{j},\zeta)  +hW_{j}%
^{1}(  \eta^{j},\zeta)  +...\ \ \ \text{ in }\Omega_{j}(
h)  \cap V_{j},\label{3.2}\\
u(  h,x)   &  =W_{0j}^{0}(  \eta^{j},\zeta)  +...\text{
\ \ \ \ in }\Omega_{0}(  h)  \cap V_{j}. \label{3.3}%
\end{align}
Aiming to determine the entries $W_{j}^{0}$ and $W_{j}^{1}$ in (\ref{3.2}), we apply the Taylor formula to $U_{j}^{0}$ and write
\begin{equation}
U_{j}^{0}(  z)  =U_{j}^{0}(  0)  +z\partial_{z}U_{j}%
^{0}(  0)  +O(  z^{2})  =U_{j}^{0}(  0)
+h\zeta\partial_{z}U_{j}^{0}(  0)  +O(  h^{2}\zeta^{2})
. \label{3.4}%
\end{equation}
The matching procedure (see, e.g., \cite{VanDyke, Ilinbook}, \cite[Ch.2]%
{MaNaPl}) requires the entries to inherit an asymptotic behavior from terms in
(\ref{3.4}). In this way we readily set
\begin{equation}
W_{j}^{0}(  \eta^{j},\zeta)  =U_{j}^{0}(  0)
\label{3.5}%
\end{equation}
and, furthermore, we fix the behavior of $W_{j}^{1}$ at infinity as follows:%
\begin{equation}
W_{j}^{1}(  \eta^{j},\zeta)  =\zeta\partial_{z}U_{j}^{0}(
0)  +o(  1)  ,\ \ \zeta\rightarrow+\infty. \label{3.6}%
\end{equation}

The constant function (\ref{3.5}) does not bring a discrepancy into the
transmission conditions (\ref{1.16}). Satisfying the other transmission
condition (\ref{1.15}), we subject $W_{0j}^{0}$ to the Dirichlet condition
\begin{equation}
W_{0j}^{0}(  \eta^{j},\zeta)  =U_{j}^{0}(  0)
,\ \ \ \ (  \eta^{j},\zeta)  \in\upsilon_{j}^{1}=\partial\omega
_{j}\times(  0,1)  . \nonumber%
\end{equation}
on the lateral boundary of the perforated layer $\Lambda_{j}$, see
(\ref{2.35}) and (\ref{1.9}). As in Section \ref{sect2.4}, the harmonic
function $W_{0j}^{0}$ satisfies the homogeneous Neumann condition at the bases
of the layer and, hence, it becomes independent of the variable $\zeta$. With some coefficient $A_j$, we
thus set%
\begin{equation}
W_{0j}^{0}(  \eta^{j},\zeta)  =U_{j}^{0}(  0)
+A_{j}\mathbf{W}_{j}(  \eta^{j}).  \label{3.A}%
\end{equation}

In view of (\ref{3.2}), (\ref{3.3}), and (\ref{3.5}), (\ref{3.A}), the
transmission condition (\ref{1.16}) converts into%
\[
h^{-1}\gamma_{j}\partial_{\nu}(  U_{j}^{0}(  0)  +hW_{j}%
^{1}(  \xi^{j})  )  =\partial_{\nu}(  U_{j}^{0}(
0)  +A_{j}\mathbf{W}_{j}(  \eta^{j})  )  ,\ \ \ \xi
^{j}\in\upsilon_{j}^{1},
\]
and, hence, yields the following Neumann condition on the ring
$\upsilon_{j}^{1}\subset\partial Q_{j}$ near the cylinder end:%
\begin{equation}
\gamma_{j}\partial_{\nu}W_{j}^{1}(  \eta^{j},\zeta)  =A_{j}%
\partial_{\nu}\mathbf{W}_{j}(  \eta^{j})  ,\ \ \ (  \eta
^{j},\zeta)  \in\upsilon_{j}^{1}. \nonumber
\end{equation}
In other words, the term $W_{j}^{1}$ in (\ref{3.2}) satisfies problem
(\ref{2.58}), (\ref{2.59}) with the right-hand side%
\[
g_{j}(  \eta,\zeta)  =A_{j}\partial_{\nu}\mathbf{W}_{j}(  \eta^{j}),\ \zeta\in(  0,1)  ,
\quad g_{j}(  \eta,\zeta)  =0, \ \zeta>1.
\]

Lemmas \ref{lem2.7} and \ref{lem2.6} ensure the existence of the unique
solution $W_{j}^{1}$ of the asymptotic form (\ref{2.61}) where
\begin{equation}
C_{0j}=0,\ \ C_{j}=\frac{A_{j}}{\gamma_{j}\vert \omega_{j}\vert
}\int_{0}^{1}\int_{\partial\omega_{j}}\partial_{\nu}\mathbf{W}_{j}(
\eta^{j})  ds_{\eta}d\zeta=-\frac{A_{j}}{\gamma_{j}\vert \omega
_{j}\vert }. \nonumber%
\end{equation}
Comparing this form with (\ref{3.6}), we conclude that condition (\ref{3.1})
provides%
\begin{equation}
G_{j}^{0}=A_{j}. \label{3.G}%
\end{equation}
The limit problem (\ref{2.6}), (\ref{2.7}), (\ref{3.1}) has the unique
solution
\begin{equation}
U_{j}^{0}(  z)  =U_{j}^{\#}(  z)  +A_{j}\gamma_{j}^{-1} \vert \omega_{j}\vert^{-1} (  l_{j}-z_{j})
\label{3.U}%
\end{equation}
where $U_{j}^{\#}$ is independent of $A_{j}$ and verifies%
\begin{align}
-\gamma_{j}\vert \omega_{j}\vert \partial_{z}^{2}U_{j}^{\#}(
z)   &  =\vert \omega_{j}\vert f_{j}^{0}(  z)
,\ \ z\in(  0,l_{j})  ,\ \ U_{j}^{\#}(  l_{j})
=0,\label{3.P}\\
-\gamma_{j}\vert \omega_{j}\vert \partial_{z}U_{j}^{\#}(
0)   &  =0. \label{3.PN}%
\end{align}
We obviously have
\begin{equation}
U_{j}^{\#}(  0)  =\frac{1}{\gamma_{j}}\int_{0}^{l_{j}}(
l_{j}-z)  f_{j}^{0}(  z)  dz. \label{3.UO}%
\end{equation}
Using decomposition (\ref{2.60}) in (\ref{3.A}) yields%
\begin{align}
W_{0j}^{0}(  \eta^{j},\zeta)   &  =U_{j}^{0}(  z)
-(2\pi)^{-1}A_{j}(  \ln\rho_{j}-\ln c_{\log}(\omega_{j}))  +O(\rho_{j}^{-1})
=\label{3.8}\\
&  =U_{j}^{0}(  0)  -(2\pi)^{-1}A_{j}(  \ln r_{j}-\ln h-\ln c_{\log}(  \omega_{j})  )  +O(h r_{j}^{-1})  .\nonumber
\end{align}
On the other hand, by (\ref{2.18}), the linear combination (\ref{2.25})
decomposes as follows:%
\begin{equation}
U_{0}^{0}(  y)  =U_{\bot}^{0}(  P^{j})  +A_{0}%
+(2\pi)^{-1}A_{j}\ln r_{j}^{-1}+{\textstyle\sum\nolimits_{k}}A_{k}G_{kj}+O(
r_{j})  . \label{3.9}%
\end{equation}
Note that $\ln r_{j}$ in (\ref{3.8}) got the same factor $-(
2\pi)  ^{-1}A_{j}$ as in (\ref{3.9}) due to our choice (\ref{3.A}).

The matching procedure requires that the sums of the terms detached in
(\ref{3.8}) and (\ref{3.9}) coincide with each other. Thus, in view of
(\ref{3.U}), we obtain the relations%
\begin{equation}
U_{j}^{\#}(  0)  +\gamma_{j}^{-1}\vert \omega_{j}\vert^{-1} A_{j}+(2\pi)^{-1}A_{j}(  \ln h+\ln c_{\log}(
\omega_{j})  )  =U_{\bot}^{0}(  P^{j})  +A_{0}%
+{\textstyle\sum\nolimits_{k}}A_{k}G_{kj} \label{3.10}%
\end{equation}
with $j=1,...,J$. In order to write them differently%
\begin{equation}
M(  \ln h)  \overrightarrow{A}=-A_{0}E+\overrightarrow{F}
\label{3.11}%
\end{equation}
we introduce the columns%
\begin{align}
\overrightarrow{A}  &  =(  A_{1},...,A_{J})  ^{\top}\in
\mathbb{R}^{J}\label{3.13}\\
E  &  =(  1,...,1)  ^{\top}\in\mathbb{R}^{J},\ \ \ \overrightarrow
{F}=\overrightarrow{U}^{\#}(  0)  -\overrightarrow{U}_{\bot}%
^{0}(  P)  ,\label{3.14}\\
\overrightarrow{U}^{\#}(  0)   &  =(  U_{1}^{\#}(
0)  ,...,U_{J}^{\#}(  0)  )  ^{\top}%
,\ \ \ \overrightarrow{U}_{\bot}^{0}(  P)  =(  U_{\bot}%
^{0}(  P^{1})  ,...,U_{\bot}^{0}(  P^{J})  )
^{\top},\nonumber
\end{align}
while $\top$ stands for transposition, and the matrix of size $J\times J$%
\begin{equation}
M(  \ln h)  =-(2\pi)^{-1}\ln h~\mathbb{I}+G-L
\label{3.12}%
\end{equation}
where $\mathbb{I}$ is the unit $J\times J$-matrix, $G$ is defined in
(\ref{2.19}) and $L$ is a diagonal matrix, namely,
\begin{equation}
L=diag\{  (2\pi)^{-1}\ln c_{\log}(  \omega_{1})  +\gamma_{1}^{-1}\vert \omega_{1}\vert^{-1},...,(2\pi)^{-1}\ln
c_{\log}(  \omega_{J})  +\gamma_{J}\vert\omega_{J}\vert^{-1}\}  . \nonumber%
\end{equation}

Matrix (\ref{3.12}) is symmetric, see (\ref{2.20}), and positive definite for
$h\in(0,h_{0}]$ when $h_{0}\in(0,1)$ is sufficiently small and, therefore,
$-\ln h=|\ln h|$. We then have%
\begin{equation}
\overrightarrow{A}=-A_{0}M(  \ln h)  ^{-1}E+M(  \ln h)
^{-1}\overrightarrow{F}. \label{3.15}%
\end{equation}
Recalling restriction (\ref{2.26}) in Proposition \ref{prop2.1}, we multiply
(\ref{3.15}) scalarly by the column $E$ and obtain%
\begin{equation}
A_{0}=m(\ln h)^{-1}(\langle f_{0}^{0}\rangle_0+E^{\top}M(\ln h)  ^{-1}\overrightarrow{F})
\label{3.16}%
\end{equation}
with the positive coefficient%
\begin{equation}
m(  \ln h)  =E^{\top}M(  \ln h)  ^{-1}E=O(
\vert \ln h\vert ^{-1})  . \label{3.m}%
\end{equation}
Inserting (\ref{3.16}) into (\ref{3.15}) finishes our construction of the
terms $U_{0}^{0}$ in (\ref{2.25}) and $U_{j}^{0}$ in (\ref{3.U}) of the
asymptotic ans\"{a}tze (\ref{2.3}), (\ref{2.10}). Due to (\ref{3.12}) and
(\ref{3.m}) all their ingredients with exception of $U_{\bot}^{0}$ and
$U_{j}^{\#}$ are rational functions in $\ln h$ while the estimates%
\begin{equation}
\vert A_{0}(  \ln h)  \vert \leq c\vert \ln
h\vert \vert \overrightarrow{F}\vert ,\ \ \ \vert
A_{j}(  \ln h)  \vert \leq c\vert \overrightarrow
{F}\vert ,\ \ \ j=1,...,J, \label{3.ln}%
\end{equation}
are valid where, according to the compact embeddings $H^{2}(  \omega
_{0})  \subset C(  \omega_{0})  $ and $H^{1}(
0,l_{j})  \subset C[  0,l_{j}]  ,$
\begin{equation}
\vert \overrightarrow{F}\vert \leq\vert \overrightarrow
{U}_{\bot}^{0}(  P)  \vert +\vert \overrightarrow
{U}^{\#}(  0)  \vert \leq c( \Vert U_{\bot}%
^{0};H^{2}(  \omega_{0})  \Vert +{\textstyle\sum\nolimits_{j}}\Vert
U_{j}^{\#};H^{1}(  0,l_{j})  \Vert)  . \label{3.FUU}%
\end{equation}
The inner expansions (\ref{3.2}) and (\ref{3.3}) are completed, too.

\subsection{The particular case $J=1$ and $\alpha=1$.\label{sect3.2}}

The mushroom shape of the junction, fig. \ref{f2},b, simplifies all the
above-obtained asymptotic terms because matrix (\ref{3.12}) at $J=1$ becomes
the scalar%
\begin{equation}
M(  \ln h)  =(2\pi)^{-1}(  \vert \ln h\vert -\ln
c_{\log}(  \omega_{1})  )  +G_{11}-\gamma_{1}^{-1}\vert\omega_{1}\vert^{-1} \nonumber%
\end{equation}
and, moreover, $E=1\in\mathbb{R}$ and $m(  \ln h)  =M(  \ln
h)  ^{-1}$ in (\ref{3.m}). In view of (\ref{2.BOT}) and (\ref{3.UO}),
now formulas (\ref{3.16}) and (\ref{3.15}) read as follows:%
\[
A_{0}(  \ln h)=M(  \ln h)  \langle f_{0}^{0}\rangle_0 +\frac{1}{\gamma_{1}}\int_{0}^{l_{1}}(
l_{1}-z)  f_{1}^{0}(  z)  dz-\int_{\omega_{0}}G_{1}(
y)  f_{0}^{0}(  y)  dy,\quad
A_{1}(  \ln h) =-\langle f_{0}^{0}\rangle_0.
\]
The latter apparently coincides with (\ref{2.26}). Only under the
orthogonality condition $\langle f_0^0\rangle_0=0$ the term $A_{0}(  \ln h)  $ in
(\ref{2.25}) and the solution $U_{0}^{0}$ itself stay bounded as
$h\rightarrow+0$ while $A_{1}(  \ln h)  G_{1}(y)$
disappears, too.

Let us outline the asymptotic structures of the solution $u(h,x)$ to problem
(\ref{1.10}), (\ref{1.11}), (\ref{1.d}), (\ref{1.13})-(\ref{1.16}) with the
Dirichlet condition on the lateral side of the plate $\Omega_{0}$. For any
$J$, the outer (\ref{2.3}), (\ref{2.10}) and inner (\ref{3.2}), (\ref{3.3})
expansions, of course, keep their forms but now the constant term $A_{0}$ in
the singular solution%
\begin{equation}
U_{0}^{0}(  y)  =U_{0}^{\#}(  y)  +{\textstyle\sum\nolimits_{j}} %
A_{j}G_{j}(  y)  \label{3.Dir}%
\end{equation}
is absent while $U_{0}^{\#}\in H^{2}(  \omega_{0})  $ is nothing
but the unique solution of%
\begin{equation}
-\Delta_{y}U_{0}^{\#}(  y)  =f_{0}^{0}(  y)
,\ y\in\omega_{0},\ \ \ \ \ \ \ \ U_{0}^{\#}(  y)  =0,\ y\in
\partial\omega_{0}, \label{3.DO}%
\end{equation}
and $G_{j}(y)=G(y,P^{j})$ denotes the standard Green function of the Dirichlet
Laplacian. Although the matching procedure resulting in the system
(\ref{3.11}) with $A_{0}=0$ remains the same as in Section \ref{sect3.1},
properties of the coefficient column (\ref{3.13}) in (\ref{3.Dir}) change
crucially. Examining the solution%
\begin{align}
\overrightarrow{A}  &  =(  A_{1},...,A_{J})  ^{\top}=M(  \ln
h)  ^{-1}(  \overrightarrow{U}^{\#}(  0)
-\overrightarrow{U}_{0}^{\#}(  P)  )  ,\label{3.AD}\\
\overrightarrow{U}_{0}^{\#}(  P)   &  =(  U_{0}^{\#}(
P^{1})  ,...,U_{0}^{\#}(  P^{J})  )  ^{\top}\nonumber
\end{align}
of the algebraic system, one sees that, instead of (\ref{3.ln}), there hold
the estimates $\vert A_{j}(  \ln h)  \vert \leq c\vert \ln
h\vert ^{-1}$. In other words, singular components in (\ref{3.Dir}) get small coefficients.
%About the case $J=1$, we only mention that%
%\[
%A_{0}=M(  \ln h)  ^{-1}(  \frac{1}{\gamma_{1}}\int_{0}^{l_{1}%
%}(  l_{1}-z)  f_{1}^{0}(  z)  dz-\int_{\omega_{0}}%
%G_{1}(  y)  f_{0}^{0}(  y)  dy)
%\]
%with $M(\ln h)$ as in (\ref{3.M1}).

\subsection{The case $\alpha=0.$\label{sect3.3}}

As mentioned in Section \ref{sect1.2}, a change of the exponent $\alpha$ in
relation (\ref{1.17}), comparing properties of elements of junction
(\ref{1.4}), may crucially modify the asymptotic ans\"{a}tze for the solution
$u(  h,x)  $ of problem (\ref{1.10})-(\ref{1.16}). At $\alpha=0,$
for example in the homogeneous ($\gamma_{j}=1$) junction, we replace
expansions (\ref{2.10}) and (\ref{2.3}), suitable for $\alpha=1,$ with the
following ones involving terms of order $h^{-1}:$%
\begin{align}
u_{0}(  h,x)   &  =h^{-1}a_{0}+U_{0}^{0}(  y)
+hU_{0}^{1}(  y,\zeta)  +h^{2}U_{0}^{2}(  y,\zeta)
+...,\label{3.34}\\
u_{j}(  h,x)   &  =h^{-1}U_{j}^{-1}(  z)  +U_{j}%
^{0}(  z)  +hU_{j}^{1}(  \eta^{j},z)  +h^{2}U_{j}%
^{2}(  \eta^{j},z)  +..., \label{3.35}%
\end{align}
where $a_{0}$ is a constant to be computed and, as we will confirm below,
$U_{j}^{-1}$ are linear functions,%
\begin{equation}
U_{j}^{-1}(  z)  =a_{0}(  1-l_{j}^{-1}z)  . \label{3.l}%
\end{equation}
Since a constant and the linear function (\ref{3.l}) verify the homogeneous
differential equations (\ref{1.10}) and (\ref{1.11}) as well as the boundary
conditions (\ref{1.12}) and (\ref{1.13}), (\ref{1.14}), the terms $U_{0}^{k}$
in (\ref{3.34}) and $U_{j}^{k}$ in (\ref{3.35}) with $k\geq0$ keep their forms
given in Sections \ref{sect2.2} and \ref{sect2.1}, respectively. Furthermore,
the main asymptotic terms in the above ans\"{a}tze satisfy the transmission
condition (\ref{1.15}) but leave a discrepancy $O(  1)  $ in the
condition (\ref{1.16}). Hence, the inner expansion which
substitutes for (\ref{3.2}) and (\ref{3.3}), looks as follows:%
\begin{equation}
u(  h,x)  =h^{-1}w_{j}^{-1}(  \eta^{j},\zeta)
+w_{j}^{0}(  \eta^{j},\zeta)  +...\ \ \ \text{in \ }\Xi(
h)  \cap V_{j} \label{3.36}%
\end{equation}
where $w_{j}^{-1}$ and $w_{j}^{0}$ are to be chosen as appropriate solutions
of the homogeneous problem (\ref{2.28})-(\ref{2.33}) studied in Section
\ref{sect2.3}. Evidently, due to (\ref{3.34})-(\ref{3.l}) we set%
\begin{equation}
w_{j}^{-1}(  \eta^{j},\zeta)  =a_{0}. \label{3.37}%
\end{equation}
Moreover, the matching procedure and formula (\ref{3.l}) requires that%
\begin{equation}
w_{j}^{0}(  \eta^{j},\zeta)  =-a_{0}l_{j}^{-1}\zeta+O(
1)  \text{ \ \ in }Q_{j},\ \ \ \zeta\rightarrow+\infty. \label{3.38}%
\end{equation}
All solutions of the homogeneous problem (\ref{2.28})-(\ref{2.33}) with the
asymptotic behavior (\ref{3.38}) have been described in Theorem \ref{th2.5}
and, therefore,%
\begin{equation}
w_{j}^{0}(  \eta^{j},\zeta)  =-a_{0}l_{j}^{-1}\vert \omega
_{j}\vert \gamma_{j}\mathbf{w}^{j}(  \eta^{j},\zeta)  +b_{j}
\label{3.40}%
\end{equation}
with a certain constant $b_{j}.$ Using decomposition (\ref{2.50}) of
$\mathbf{w}^{j}$ in the layer $\Lambda_{j},$ we now obtain that%
\begin{align}
h^{-1}w_{j}^{-1}(  \eta^{j},\zeta)  +w_{j}^{0}(  \eta
^{j},\zeta)   &  =h^{-1}a_{0}+a_{0}\gamma_{j}\vert \omega_{j}\vert(2\pi l_{j})^{-1}\ln\rho_{j}+b_{j}+O(\rho_{j}^{-1})=\label{3.41}\\
&  =h^{-1}a_{0}+a_{0}\gamma_{j}\vert \omega_{j}\vert(2\pi l_{j})^{-1}(\ln r_{j}-\ln h)+b_{j}+O(hr_{j}^{-1})  \ \ \ \ \text{in }\Lambda_{j}.\nonumber
\end{align}
Matching (\ref{3.41}) and (\ref{3.34}) subjects the function $U_{0}^{0}$ to
the conditions%
\begin{equation}
U_{0}^{0}(  y)  =a_{0}\gamma_{j}\vert \omega_{j}\vert(2\pi l_{j})^{-1}(\ln r_{j}-\ln h)
+b_{j}+O(  r_{j})  ,\ \ r_{j}\rightarrow0. \label{3.42}%
\end{equation}
In other words, representation (\ref{2.25}) of the singular solution
$U_{0}^{0}$ to the limit problem (\ref{2.21}), (\ref{2.14}) involves the
coefficients%
\begin{equation}
A_{j}=-a_{0}l_{j}^{-1}\vert \omega_{j}\vert \gamma_{j}.
\label{3.43}%
\end{equation}
Furthermore, restriction (\ref{2.26}) requires that%
\begin{equation}
a_{0}=(  l_{1}^{-1}\vert \omega_{1}\vert \gamma_{1}%
+...+l_{J}^{-1}\vert \omega_{J}\vert \gamma_{J})  ^{-1}%
\langle f_{0}^{0}\rangle_0. \label{3.44}%
\end{equation}
Thus, the main asymptotic terms in expansions (\ref{3.34}) and (\ref{3.35})
have been determined.

The second term $U_{0}^{0}$ of the outer expansion (\ref{3.34}) in the
perforated plate $\Omega_{\bullet}(  h)  $ is defined up to the
constant addendum $A_{0}$ in (\ref{2.25}). Comparing formulas (\ref{2.25}),
(\ref{2.18}) and (\ref{3.42}) yields%
\begin{equation}
b_{j}=U_{\bot}^{0}(  P^{j})  +A_{0}+{\textstyle\sum\nolimits_{k}}A_{k}%
G_{kj}(  y)
+a_{0}\gamma_{j}\vert \omega_{j}\vert(2\pi l_{j})^{-1}\ln h. \label{3.45}%
\end{equation}
Since, by (\ref{3.37}), (\ref{3.40}) and (\ref{2.49}), we have%
\begin{align*}
h^{-1}w_{j}^{-1}(  \eta^{j},\zeta)  +w_{j}^{0}(  \eta
^{j},\zeta)   &  =h^{-1}a_{0}-a_{0}l_{j}^{-1}\zeta+a_{0}\vert
\omega_{j}\vert \gamma_{j}\mathbf{q}_{j}+b_{j}+O(  e^{-\delta\zeta
})  =\\
&  =h^{-1}a_{0}(  1-l_{j}^{-1}z)  +b_{j}+a_{0}\vert \omega
_{j}\vert \gamma_{j}\mathbf{q}_{j}+O(  e^{-\delta z/h})
\text{ \ \ in }Q_{j},
\end{align*}
matching the inner expansion (\ref{3.36}) with the outer expansion
(\ref{3.35}) specified as%
\[
h^{-1}U_{j}^{-1}(  z)  +U_{j}^{0}(  z)  =h^{-1}%
a_{0}(  1-l_{j}^{-1}z)  +U_{j}^{0}(  0)  +O(
z)
\]
leads to the boundary condition%
\begin{equation}
U_{j}^{0}(  0)  =b_{j}+a_{0}\gamma_{j}\vert \omega
_{j}\vert \mathbf{q}_{j}. \label{3.46}%
\end{equation}
According to (\ref{3.45}), the solution of the limit problem (\ref{2.6}),
(\ref{2.7}) and (\ref{3.46}) in $(  0,l_{j})  $ looks
as follows:%
\begin{equation}
U_{j}^{0}(  z;\ln h)  =(  A_{0}+a_{0}\gamma_{j}\vert \omega_{j}\vert(2\pi l_{j})^{-1}\ln h)
(1-l_{j}^{-1}z)  +U_{j}^{\#}(z)  \label{3.47}%
\end{equation}
where $A_{0}$ is not fixed yet but $U_{j}^{\#}$ does not depend on $\ln h$ and
is unambiguously defined as a unique solution of problem (\ref{3.P}) equipped
with the Dirichlet condition%
\begin{equation}
U_{j}^{\#}(  0)  =U_{\bot}^{0}(  P^{j})  +{\textstyle\sum\nolimits_{k}}A_{k}G_{kj}(  y)  +a_{0}\gamma_{j}\vert \omega
_{j}\vert \mathbf{q}_{j} \nonumber%\label{3.PD}%
\end{equation}
instead of the Neumann ones (\ref{3.PN}). One may detect $A_{0}$ by
constructing next asymptotic terms but we avoid to encumber our paper with
computations of access, although we will get troubles with the unknown $A_{0}$
in the justification procedure. We only mention that $A_{0}=A_{0}(  \ln
h)  $ could be proved to depend linearly on $\ln h.$ In this way terms
(\ref{3.47}) of the asymptotic ans\"{a}tze (\ref{3.35}) on the rods become linear in $\ln h,$ too. Recall that
$a_{0}$ and the constructed terms $U_{j}^{-1}$, $U_{0}^{0}-A_{0}$
are independent of $\ln h.$

\subsection{The particular case $J=1$ and $\alpha=0.$\label{sect3.4}}

For the junction $\Xi(  h)  =\Omega_{0}(  h)  \cup
\Omega_{1}(  h)  $ of a plate and a single rod, the main asymptotic
terms in ans\"{a}tze (\ref{3.34}) and (\ref{3.35}) are determined by the
formulas $a_{0}=\gamma_{1}^{-1}\vert \omega_{j}\vert^{-1}l_{1}\langle f^0_0\rangle_0$
and (\ref{3.l}). Other formulas in the previous section only need the
specification $j=1$.

As in the case $\alpha=1,$ a simplification of asymptotic ans\"{a}tze occurs
in problem (\ref{1.10}), (\ref{1.11}), (\ref{1.d}), (\ref{1.13})-(\ref{1.16})
with the Dirichlet condition on the lateral side of the plate $\Omega
_{0}(  h)  .$ First of all, terms of order $h^{-1}$ disappear and
we return to the outer expansions (\ref{2.3}) and (\ref{2.10}). Moreover, the
main term $U_{0}^{0}$ in (\ref{2.10}) is no longer singular and implies a
solution in $H^{2}(  \omega_{0})  $ of the Dirichlet problem
(\ref{3.DO}). Thus, the values $U_{0}^{0}(  P^{j})  $ are defined
properly and the main terms $U_{j}^{0}$ in (\ref{2.3}) become nothing but
solutions of the Dirichlet problems in the interval $(  0,l_{j})  $
composed from (\ref{2.6}), (\ref{2.7}) and
\begin{equation}
U_{j}^{0}(  0)  =U_{0}^{0}(  P^{j})  . \label{3.DD}%
\end{equation}
Owing to such elementary asymptotic structures it is not even necessary to
deal with inner expansions in the justification scheme. However, we mention
that the inner expansion (\ref{3.36}), of course, loses the term $h^{-1}%
w_{j}^{-1}$ and gains the constant term $w_{j}^{0}(  \xi^{j})
=U^{0}(  P^{j})  .$

\section{Estimation of the asymptotic remainders\label{sect4}}

Since, in view of our consideration in Section \ref{sect3}, boundary layer
effects play an important role in the obtained asymptotic structures, error
estimates for the asymptotic expansions are seriously based on certain
weighted inequalities presented below. A general approach for derivation of
such inequalities on thin domains and junctions with different limit
dimensions can be found in the review paper \cite{na396}.

\subsection{Weighted estimates on thin elements of the junction\label{sect4.1}%
}

We integrate in $y\in\omega_{j}^{h}$ the Hardy inequality (\ref{2.42}) with $L=l_j, \zeta =l_j-z, W(\zeta)=u_j(y,l_j-\zeta)$,
%\begin{equation}
%\int_{0}^{l_{j}}\vert l_{j}-z\vert ^{-2}\vert u_{j}(
%y,z)  \vert ^{2}dz\leq4\int_{0}^{l_{j}}\vert \frac{\partial
%u_{j}}{\partial z}(  y,z)  \vert ^{2}dz \label{4.1}%
%\end{equation}
which is true due to the condition $u_{j}(  y,l_{j})  =0,$ cf.
(\ref{1.14}), and we obtain%
\begin{equation}
\Vert \vert l_{j}-z\vert ^{-1}u_{j};L^{2}(  \Omega
_{j}(  h)  )  \Vert \leq2\Vert \partial_{z}%
u_{j};L^{2}(  \Omega_{j}(  h)  )  \Vert .\label{4.1}
\end{equation}
%For reader's convenience, we mention that (\ref{4.1}) follows from
%Newton-Leibnitz formula, as follows:%
%\begin{align}
%\int_{0}^{l}t^{-2}\vert U(  t)  \vert ^{2}dt  &
%=2\int_{0}^{l}t^{-2}\int_{0}^{t}U(  \tau)  \frac{\partial
%U}{\partial\tau}(  \tau)  d\tau dt\leq2\int_{0}^{l}\vert
%U(  \tau)  \vert \vert \frac{\partial U}{\partial\tau
%}(  \tau)  \vert \int_{\tau}^{l}t^{-2}dtd\tau=\label{4.0}\\
%&  =2\int_{0}^{l}\vert U(  \tau)  \vert \vert
%\frac{\partial U}{\partial\tau}(  \tau)  \vert (
%\frac{1}{\tau}-\frac{1}{l})  d\tau\leq2\int_{0}^{l}\frac{1}{\tau
%}\vert U(  \tau)  \vert \vert \frac{\partial
%U}{\partial\tau}(  \tau)  \vert d\tau\leq\nonumber\\
%&  \leq2(  \int_{0}^{l}\tau^{-2}\vert U(  \tau)
%\vert ^{2}d\tau)  ^{1/2}(  \int_{0}^{l}\vert
%\frac{\partial U}{\partial\tau}(  \tau)  \vert ^{2}%
%d\tau)  ^{1/2}.\nonumber
%\end{align}
Furthermore, setting%
\begin{equation}
\overline{u}_{j}(  z)  =\frac{1}{\operatorname*{mes}_{2}\omega
_{j}^{h}}\int_{\omega_{j}^{h}}u_{j}(  y,z)  dy, \label{4.3}%
\end{equation}
the Poincar\'{e} inequality on cross-sections of the rod $\Omega_{j}(
h)  $ yields%
\begin{equation}
\Vert u_{j}-\overline{u}_{j};L^{2}(  \Omega_{j}(  h)
)  \Vert \leq c(  \omega_{j})  h\Vert \nabla
_{y}u_{j};L^{2}(  \Omega_{j}(  h)  )  \Vert
\label{4.4}%
\end{equation}
where $c(  \omega_{j})  ^{-2}$ is the first positive eigenvalue of
the Neumann Laplacian in the domain $\omega_{j}.$

For $u_{0}\in H^{1}(  \Omega_{0}(  h)  )  ,$ we set%
\begin{equation}
\overline{u}_{0}(  y)  =\frac{1}{h}\int_{0}^{h}u_{0}(
y,z)  dz,\ \ \ u_{0}(  x)  =u_{\bot}(  x)  +b_{0},
\quad b_{0}=\frac{1}{\operatorname*{mes}_{3}\Omega_{0}(  h)  }%
\int_{\Omega_{0}(  h)  }u_{0}(  x)  dx. \label{4.6}%
\end{equation}
Similarly to (\ref{4.4}), we have%
\begin{equation}
\Vert u_{0}-\overline{u}_{0};L^{2}(  \Omega_{0}(  h)
)  \Vert \leq\pi^{-1}h\Vert \partial_{z}u_{0};L^{2}(
\Omega_{0}(  h)  )  \Vert \label{4.7}%
\end{equation}
Using the change of variables $x\mapsto(  y,\zeta)  =(
y,h^{-1}z)  $ and therefore transforming $\Omega_{0}(  h)  $
into the domain $\Omega_{0}(  1)  =\omega_{0}\times(
0,1)  $ independent of $h$, we obtain
\begin{align}
\Vert u_{\bot};L^{2}(  \Omega_{0}(  h)  )
\Vert ^{2}  &  =\Vert u_{0}-b_{0};L^{2}(  \Omega_{0}(
h)  )  \Vert ^{2}=h\Vert u_{0}-b_{0};L^{2}(
\Omega_{0}(  1)  )  \Vert ^{2}\leq\label{4.8}\\
&  \leq c_{0}h\Vert \nabla_{(  y,\zeta)  }u_{0};L^{2}(
\Omega_{0}(  1)  )  \Vert ^{2}\leq\nonumber\\
&  \leq c_{0}h(  \Vert \nabla_{y}u_{0};L^{2}(  \Omega
_{0}(  1)  )  \Vert ^{2}+\Vert \partial_{\zeta
}u_{0};L^{2}(  \Omega_{0}(  1)  )  \Vert
^{2})  \leq\nonumber\\
&  \leq c_{0}h(  \Vert \nabla_{y}u_{0};L^{2}(  \Omega
_{0}(  1)  )  \Vert ^{2}+h^{-2}h_{0}^{2}\Vert
\partial_{\zeta}u_{0};L^{2}(  \Omega_{0}(  1)  )
\Vert ^{2})  =\nonumber\\
&  =c_{0}(  \Vert \nabla_{y}u_{0};L^{2}(  \Omega_{0}(
h)  )  \Vert ^{2}+h_{0}^{2}\Vert \partial_{z}%
u_{0};L^{2}(  \Omega_{0}(  h)  )  \Vert
^{2})  \leq\nonumber\\
&  \leq c_{0}\max\{  1,h_{0}^{2}\}  \Vert \nabla_{x}%
u_{0};L^{2}(  \Omega_{0}(  h)  )  \Vert
^{2},\nonumber
\end{align}
where $c_{0}$ is the first positive eigenvalue of the Neumann Laplacian in $\Omega_{0}(  1)  .$ Note that in (\ref{4.8}) we have used
the Poincar\'{e} inequality in $\Omega_{0}(  1)  $ based on the
orthogonality condition $(  u_{0}-b_{0},1)_{\Omega_{0}(1)}=0$ in (\ref{4.6}).

\subsection{Weighted estimates on the whole junction\label{sect4.2}}

In relation (\ref{4.1}) we give an inequality for the restriction $u_{j}=u|_{\Omega_{j}(  h)  }$. However, the relations (\ref{4.7}) and (\ref{4.8}) do not provide an inequality for the restriction $u_{0}=u|_{\Omega_{0}(
h)  }$ but only for its components $u_{0}-\overline{u}_{0}$ and
$u_{\bot}=u_{0}-b_{0}.$ Notice that, in contrast to the definitions in Section
\ref{sect1.1}, we treat $u_{0}$ as a function in the intact plate (\ref{1.1}),
not in the perforated plate (\ref{1.5}), so that here we have%
\begin{equation}
u_{0}(  x)  =u_{j}(  x)  ,\ \ \ x\in\theta_{j}%
^{h}=\omega_{j}^{h}\times(  0,h)  \label{4.9}%
\end{equation}
where $\theta_{j}^{h}=\Omega_{0}(  h)  \cap\Omega_{j}(
h)  $ is a small embedded piece of the rod  into the plate. Based on
(\ref{4.6}) and (\ref{4.9}), we write%
\begin{equation}
\vert b_{0}\vert =\frac{1}{\operatorname*{mes}_{3}\theta_{j}^{h}%
}\left\vert \int_{\theta_{j}^{h}}(  u_{j}(  x)  -u_{\bot
}(  x)  )  dx\right\vert \leq ch^{-3/2}(  \Vert
u_{j};L^{2}(  \theta_{j}^{h})  \Vert +\Vert u_{\bot
};L^{2}(  \theta_{j}^{h})  \Vert )  \label{4.10}%
\end{equation}
because the volume of $\theta_{j}^{h}$ is of order $h^{3}.$ We derive from (\ref{2.42}) that%
\begin{align}
h^{-1}\Vert u_{j};L^{2}(  \theta_{j}^{h})  \Vert ^{2}
&  =h^{-1}\int_{0}^{h}\int_{\omega_{j}^{h}}\vert u_{j}(
y,\mathfrak{z})  \vert ^{2}dyd\mathfrak{z}=h^{-1}\int_{0}^{h}%
\int_{\omega_{j}^{h}}\left\vert \int_{\mathfrak{z}}^{l_{j}}\frac{\partial
}{\partial z}(  \chi_{j}(  z)  u_{j}(  y,z)
)  dz\right\vert ^{2}dyd\mathfrak{z}\leq\nonumber\\
&  \leq c_{j}h^{-1}\int_{0}^{h}\int_{\omega_{j}^{h}}\int_{0}^{l_{j}}\left(
\left\vert \frac{\partial u_{j}}{\partial z}(  y,z)  \right\vert
^{2}+\vert u_{j}(  y,z)  \vert ^{2}\right)
dzdyd\mathfrak{z}=\nonumber\\
&  =c_{j}\int_{\Omega_{j}(  h)  }(  \vert \nabla_{x}%
u_{j}(  x)  \vert ^{2}+\vert u_{j}(  x)
\vert ^{2})  dx\leq C_{j}\Vert \nabla_{x}u_{j};L^{2}(
\Omega_{j}(  h)  )  \Vert ^{2}.\label{3.1415}
\end{align}
To examine the last norm in (\ref{4.10}), we apply the Hardy inequality
(\ref{2.43}) with logarithm to the product $\chi_{0j}u_{\bot},$ where
$\chi_{0j}\in C_{c}^{\infty}(  \omega_{0})  $ is a cut-off function
such that
\begin{equation}
\chi_{0j}(  y)  =1\text{ for }y\in\mathbb{B}_{R_{0}/2}(
P^{j})  \text{ \ and \ }\chi_{0j}(  y)  =0\text{ for }%
y\notin\mathbb{B}_{R_{0}}(  P^{j})  \label{4.27}%
\end{equation}
and radius $R_{0}>0$ is fixed to fulfil $\mathbb{B}_{R_{0}}(
P^{j})  \subset\omega_{0}$ and $\mathbb{B}_{R_{0}}(  P^{j})
\cap\mathbb{B}_{R_{0}}(  P^{j})  =\emptyset$ as $j\neq k.$

Owing to (\ref{4.8}) we thus obtain%
\begin{align}
\int_{0}^{h}\int_{\omega_{0}}r_{j}^{-2}(  1+\vert \ln r_{j}%
\vert )  ^{-2}\vert u_{\bot}(  y,z)  \vert
^{2}dydz  &  \leq c_{j}\int_{\Omega_{0}(  h)  }(  \vert
\nabla_{y}u_{\bot}(  x)  \vert ^{2}+\vert u_{\bot
}(  x)  \vert ^{2})  dx\leq\label{4.11}\\
&  \leq C_{j}\int_{\Omega_{0}(  h)  }\vert \nabla_{y}u_{\bot
}(  x)  \vert ^{2}dx\leq C_{j}\int_{\Omega_{0}}\vert
\nabla_{x}u_{0}(  x)  \vert ^{2}dx\nonumber
\end{align}
where $r_{j}=dist(  x,P^{j})  .$ Hence,%
\begin{align}
\int_{\theta_{j}^{h}}\vert u_{\bot}(  x)  \vert ^{2}dx
&  \leq c_{j}h^{2}(  1+\vert \ln h\vert )  ^{2}%
\int_{\theta_{j}^{h}}r_{j}^{-2}(  1+\vert \ln r_{j}\vert
)  ^{-2}\vert u_{\bot}(  x)  \vert ^{2}%
dx\leq\label{4.12}\\
&  \leq c_{j}h^{2}(  1+\vert \ln h\vert )
^{2}\Vert \nabla_{x}u_{0};L^{2}(  \Omega_{0}(  h)
)  \Vert ^{2}.\nonumber
\end{align}
From (\ref{4.10}) and (\ref{4.11}), (\ref{4.12}), we get an estimate of the
component $b_{0}$ in the representation (\ref{4.6}).

\begin{theorem}
\label{th4.1}Let $\alpha\geq0$ in (\ref{1.17}). There hold the weighted
inequalities%
\begin{gather}
\min\{  h^{-\alpha+1},(  1+\vert \ln h\vert )
^{-2}\}  \Vert r^{-1}(  1+\vert \ln r\vert )
^{-1}u_{0};L^{2}(  \Omega_{0}(  h)  )  \Vert
^{2}+\label{4.13}\\
+h^{-\alpha}\Vert (  l_{j}-z)  ^{-1}u_{j};L^{2}(
\Omega_{j}(  h)  )  \Vert ^{2}\leq c_{\Xi}a(
u,u;\Xi(  h)  ) \nonumber
\end{gather}
where $r=\min\{  1,r_{1},...,r_{J}\}  ,$ $a$ is the quadratic form
(\ref{1.20}) with the coefficients $\gamma_{j}(  h)  =O(
h^{-\alpha})  $ and the constant $c_{\Xi}$ is independent of the
parameter $h\in(  0,h_{0}]  $ and the function $u\in H_{0}%
^{1}(  \Xi(  h)  ,\Gamma(  h)  )  .$
\end{theorem}

\textbf{Proof.} It suffices to take into account the estimates (\ref{4.8}), (\ref{3.1415})
 together with the calculation%
\begin{align}
\Vert r_{j}^{-1}(  1+\vert \ln r_{j}\vert )
^{-1}b_{0};L^{2}(  \Omega_{0}(  h)  )  \Vert ^{2}
&  \leq ch\vert b_{0}\vert ^{2}\int_{0}^{R}r_{j}^{-2}(
1+\vert \ln r_{j}\vert )  ^{-2}r_{j}dr_{j}\leq Ch\vert
b_{0}\vert ^{2}\leq\label{4.14}\\
&  \leq Ch^{-2}(  \Vert u_{j};L^{2}(  \theta_{j}^{h})
\Vert ^{2}+\Vert u_{\bot};L^{2}(  \theta_{j}^{h})
\Vert ^{2})  \leq\nonumber\\
&  \leq Ch^{-2}(  h)  \Vert \nabla_{x}u_{j};L^{2}(
\Omega_{j}(  h)  )  \Vert ^{2}+h^{2}(
1+\vert \ln h\vert )  ^{2}\Vert \nabla_{x}u_{0}%
;L^{2}(  \Omega_{0}(  h)  )  \Vert ^{2}%
\leq\nonumber\\
&  \leq C(  h^{-\alpha+1}+(  1+\vert \ln h\vert )
^{2})  a(  u,u;\Xi(  h)  ) \nonumber
\end{align}
which is based on (\ref{4.10})-(\ref{4.12}). Notice that the last inequality
in (\ref{4.14}) is valid due to the assumption $\alpha\geq0$ which assures
that%
\begin{align*}
\Vert \nabla_{x}u_{0};L^{2}(  \Omega_{0}(  h)  )
\Vert ^{2}  &  =\Vert \nabla_{x}u_{0};L^{2}(  \Omega_{\bullet
}(  h)  )  \Vert ^{2}+{\textstyle\sum\nolimits_{j}} \Vert
\nabla_{x}u_{j};L^{2}(  \theta_{j}^{h})  \Vert ^{2}\leq\Vert \nabla_{x}u_{0};L^{2}(  \Omega_{\bullet}(
h)  )  \Vert ^{2}+\\
& +( h_{0}/h)^{\alpha
}{\textstyle\sum\nolimits_{j}} \Vert \nabla_{x}u_{j};L^{2}(  \theta_{j}^{h})
\Vert ^{2}\leq ca(  u,u;\Xi(  h)  )  .
\end{align*}

\begin{remark}
\label{rem4.2}Let us verify the asymptotic accuracy of the distribution of
weights on the left-hand side of (\ref{4.13}). Clearly, the exponent $-1$ of $(  l_{j}-z)  $ cannot be reduced.
Indeed, for any $\delta>0,$ the function $u_{j}^{\delta}(  x)
=(  l_{j}-z)  ^{(  2+\delta)  /4}$ belongs to $H_{0}%
^{1}(  \Omega_{j}(  h)  ;\omega_{j}^{h}(  l_{j})
)  $ and can be extended over $\Xi(  h)  $ but the integral%
\[
\int_{\Omega_{j}(  h)  }\vert l_{j}-z\vert ^{-2-\delta
}\vert u_{j}^{\delta}(  x)  \vert ^{2}dx
\]
diverges. In the case $\alpha=0$ a trial function to confirm the precision of the
inequality (\ref{4.13}) can be taken in the form $u_{0}(  x)  =1,\ \ u_{j}(  x)  =1-\chi_{j}(
z)$, where
\begin{equation}
\chi_{j}\in C_{0}^{\infty}[  0,l_{j})  ,\text{ }\chi_{j}(
z)  =1\text{ for }z<l_{j}/3\text{ \ and \ }\chi_{j}(
z)  =0\text{ for }z>2l_{j}/3. \label{4.l}%
\end{equation}
Then we have
\begin{align*}
\Vert r_{j}^{-1}(  1+r_{j})  ^{-1}u_{0};L^{2}(
\Omega_{0}(  h)  )  \Vert ^{2}  &  \geq c_{0}%
h,\ \Vert (  l_{j}-z)  ^{-1}u_{j};L^{2}(  \Omega
_{j}(  h)  )  \Vert ^{2}\geq c_{j}h^{2},\ c_{p}>0,\\
a(  u,u;\Xi(  h)  )   &  \leq c{\textstyle\sum\nolimits_{j}} \Vert
\partial_{z}\chi_{j};L^{2}(  \Omega_{j}(  h)  )
\Vert ^{2}\leq ch^{2}.
\end{align*}
We see that all terms in (\ref{4.13}) become $O(  h^{2})  .$

In the case $\alpha=1$ we assume for simplicity that the domain $\omega_{k}$
is the circle $\{  y:r_{k}<h\}  $ and set%
\[
u_{0}(  x)  =\chi_{k}^{0}(  y)  \ln\left\vert \frac{\ln
r_{k}}{\ln h}\right\vert ,\ \ u_{j}(  x)  =0.
\]
We then have
\begin{align*}
(  1+\vert \ln h\vert )  ^{-2}\Vert r_{k}%
^{-1}(  1+r_{k})  ^{-1}u_{0};L^{2}(  \Omega_{\bullet}(
h)  )  \Vert ^{2}  &  \geq ch(  1+\vert \ln
h\vert )  ^{-2}\int_{\vert \ln r_{\chi}\vert
}^{\vert \ln h\vert }(  1+\lambda)  ^{-2}\left\vert
\ln\frac{\lambda}{\vert \ln\lambda\vert }\right\vert ^{2}%
d\lambda\geq\\
&  \geq ch(  1+\vert \ln h\vert )  ^{-2}\vert
\ln\vert \ln h\vert \vert ^{2},\ \ c>0,
\end{align*}
\[
a(  u,u;\Xi(  h)  ) =(  \nabla_{x}u_{0}%
,\nabla_{x}u_{0})  _{\Omega_{\bullet}(  h)  }\leq ch\int
_{h}^{R_{\chi}}\left\vert \nabla_{y}\ln\left\vert \frac{\ln\lambda}{\ln
h}\right\vert \right\vert ^{2}rdr\leq ch\int_{h}^{R_{\chi}}\vert \ln r\vert ^{-2}\frac{dr}{r}\leq
Ch\vert \ln h\vert ^{-1}.
\]
Here, $R_{\chi}>0$ and $r_{\chi}$ are small and such that $\chi_{k}^{0}(
y)  =1$ for $r<r_{\chi}$\ and\ $\chi_{k}^{0}(  y)  =0$ for
$r>R_{\chi}.$ The above relations show that the inequality (\ref{4.13})
with $\alpha=1$ is sharp with respect to powers of the small parameter $h.$
Moreover, we conclude that the constant $c_{\Xi}$ cannot hold without a
logarithmical factor on the left. However, the authors do not know how to
confirm the optimality of the factor $(  1+\vert \ln h\vert
)  ^{-2}.$ $\blacksquare$
\end{remark}

%\begin{remark}
%\label{rem4.3}Inequality (\ref{4.13}) remains valid when $\alpha<0$. However,
%to prove it, one needs to find out a substitutor for the estimate (\ref{4.11})
%because, in view of the small coefficient $h^{-\alpha}$ in the last terms in
%(\ref{1.20}), and (\ref{1.17}), the relation%
%\[
%\Vert \nabla_{x}u_{\bot};L^{2}(  \Omega_{0}(  h)
%)  \Vert ^{2}=\Vert \nabla_{x}u_{0};L^{2}(  \Omega
%_{0}(  h)  )  \Vert ^{2}\leq ca(  u,u;\Xi(
%h)  )
%\]
%is no longer true with a constant $c$ independent of $h\in(
%0,h_{0}]  .$ Nevertheless, a function $u_{\bot}\in H^{1}(
%\Omega_{\bullet}(  h)  )  $ of mean zero in $\Omega_{\bullet
%}(  h)  $ meets the Poincar\'{e} inequality%
%\begin{equation}
%\Vert u_{\bot};L^{2}(  \Omega_{\bullet}(  h)  )
%\Vert ^{2}\leq C\Vert \nabla_{x}u_{\bot};L^{2}(
%\Omega_{\bullet}(  h)  )  \Vert ^{2}. \label{4.15}%
%\end{equation}
%Furthermore, one may use an approach in \cite[\S 2]{na396} to verify that the
%constant $C$ can be fixed independent of $h.$ The right-hand side of
%(\ref{4.15}) is equal to $c\Vert \nabla_{x}u_{0};L^{2}(
%\Omega_{\bullet}(  h)  )  \Vert ^{2}$ and does not
%exceed $ca(  u,u;\Xi(  h)  )  ,$ see (\ref{4.6}) and
%(\ref{1.20}). $\blacksquare$
%\end{remark}

\subsection{The requirements for the problem data\label{sect4.3}}

In this section we assume that the right-hand sides of equations (\ref{1.10})
and (\ref{1.11}) satisfy%
\begin{align}
f_{0}(  h,x)   &  =F_{0}(  y,h^{-1}z)  ,\ \ \ f_{j}%
(  h,x)  =h^{-\alpha}F_{j}(  h^{-1}(  y-P^{j})
,z)  ,\label{4.16}\\
F_{0}  &  \in L^{2}(  \omega_{0}\times(  0,1)  )
,\ \ \ F_{j}\in L^{2}(  \omega_{j}\times(  0,l_{j})  )
\label{4.F}%
\end{align}
and denote by $\mathcal{N}$ the sum of norms of functions (\ref{4.F}) in the
indicated spaces. We also set%
\begin{align}
f_{0}^{0}(  y)   &  =\int_{0}^{1}F_{0}(  y,\zeta)
d\zeta,\ \ \ f_{00}^{\bot}(  y,\zeta)  =F_{0}(  y,\zeta
)  -f_{0}^{0}(  y)  ,\label{4.P}\\
f_{j}^{0}(  z)   &  =\frac{1}{\vert \omega_{j}\vert
}\int_{\omega_{j}}F_{j}(  \eta^{j},z)  d\eta^{j},\ \ \ f_{j0}%
^{\bot}(  \eta^{j},z)  =F_{j}(  \eta^{j},z)  -f_{j}%
^{0}(  z) \nonumber
\end{align}
and observe that first, $f_{00}^{\bot}$ and $f_{j0}^{\bot}$ meet the
conditions (\ref{2.9}) and (\ref{2.2}), respectively, and
second,
\begin{equation}
\Vert f_{0}^{0};L^{2}(  \omega_{0})  \Vert +{\textstyle\sum\nolimits_{j}}\Vert f_{j}^{0};L^{2}(  0,l_{j})  \Vert \leq
c\mathcal{N}. \nonumber%
\end{equation}
The obtained representations%
\begin{equation}
f_{0}(  h,x)  =f_{0}^{0}(  y)  +f_{00}^{\bot}(
y,h^{-1}z)  ,\ \ \ f_{j}(  h,x)  =h^{-\alpha}(
f_{j}^{0}(  z)  +f_{j0}^{\bot}(  h^{-1}(  y-P^{j})
,z)  )  \nonumber%
\end{equation}
differ from representations (\ref{2.8}) and (\ref{2.1}) proposed in Section
\ref{sect2} in the absence of the small remainders $\widetilde{f}_{p}$ and the
big factor $h^{-1}$ on $f_{p0}^{\bot}.$ However, these simplified
representations are sufficient to demonstrate all technicalities in deriving
the error estimates and to achieve the goals of this section. We will return
to discuss the general forms of the right-hand sides in Section \ref{sect5.1}.

\textit{The case }$\alpha=1.$ Recalling materials of Sections \ref{sect2.1},
\ref{sect2.2} and \ref{sect3.1}, we observe that ingredients of the asymptotic
ans\"{a}tze (\ref{2.3}) and (\ref{2.10}) constructed from the functions
$f_{0}^{0}$ and $f_{j}^{0}$ in (\ref{4.P}) satisfy the estimate%
\begin{equation}
\Vert U_{\bot}^{0};H^{2}(  \omega_{0})  \Vert
+{\textstyle\sum\nolimits_{j}} \Vert U_{j}^{0};H^{2}(  0,l_{j})  \Vert
+\vert \ln h\vert ^{-1}\vert A_{0}\vert +{\textstyle\sum\nolimits_{j}}\vert A_{j}\vert \leq c\mathcal{N}. \label{4.20}%
\end{equation}
The first norm on the left of (\ref{4.20}) has appeared in (\ref{2.bot}) and
estimates for the norms of the solutions $U_{j}^{0}$ to problems (\ref{2.6}),
(\ref{2.7}), (\ref{3.1}) with the right-hand sides (\ref{3.G}) are evident in
view of estimates (\ref{3.ln}), (\ref{3.FUU}) which also are displayed in
(\ref{4.20}).

Inequality (\ref{4.20}) will be used in the next two sections to derive
estimates (\ref{4.55}) and (\ref{4.61}) of asymptotic remainders with the
bounds $ch\mathcal{N},$ where the factor $\mathcal{N}$ expresses the whole
dependence of these bounds on the right-hand sides (\ref{4.16}) in the
original problem in the junction $\Xi(  h)  $.

\textit{The case }$\alpha=0.$ Our justification scheme for asymptotics in the junction with $\alpha=0,$ e.g. an homogeneous junction
(cf. comment to (\ref{1.17}) and (\ref{1.18})), requires similar estimates of
ingredients of ans\"{a}tze (\ref{3.34}), (\ref{3.35}). As explained in Sections
\ref{sect3.3}, (\ref{3.44}), (\ref{2.bot}) and (\ref{3.43}), (\ref{3.47}) assure the estimate
\begin{equation}
\vert a_{0}\vert +\Vert U_{\bot}^{0};H^{2}(  \omega
_{0})  \Vert +{\textstyle\sum\nolimits_{j}} (  \vert A_{j}\vert
+(  1+\vert \ln h\vert )  ^{-1}\Vert U_{j}%
^{0};H^{2}(  0,l_{j})  \Vert )  \leq c\mathcal{N}.
\nonumber
\end{equation}

\subsection{The global approximation of the solution in the case $\alpha
=1.$\label{sect4.4}}

In order to glue the outer (\ref{3.2}), (\ref{3.3}) and inner (\ref{2.3}),
(\ref{2.10}) asymptotic expansions we introduce the following cut-off
functions:%
\begin{align}
X_{0}^{h}(  y)   &  =1-{\textstyle\sum\nolimits_{j}} \chi_{0j}^{h}(  y)
,\ \ \ \chi_{0j}^{h}(  y)  =\chi(  h^{-1}R_{j}^{-1}%
r_{j})  ,\label{4.24}\\
X_{j}^{h}(  z)   &  =1-\chi_{j}^{h}(  y)  ,\ \ \ \chi
_{j}^{h}(  y)  =\chi(  h^{-1}z)  \label{4.25}%
\end{align}
where $R_{j}$ was defined in (\ref{QL}) and $\chi\in C^\infty (\mathbb{R})$ is such that $\chi(t)=1$ for $t<1$ and $\chi(t)=0$ for $t>2$.
Notice that the function $X_{0}^{h}$ (the function $X_{j}^{h})$ is equal to
one everywhere in the plate $\Omega_{\bullet}(  h)  $ (in the rod
$\Omega_{j}(  h)  ),$ except for the vicinity of the holes
$\theta_{1}^{h},...,\theta_{J}^{h}$ (the rod end (\ref{1.8})). Clearly,%
\begin{equation}
\vert \nabla_{y}^{k}X_{0}^{h}(  y)  \vert \leq
ch^{-k},\ \ \ \vert \nabla_{z}^{k}X_{j}^{h}(  z)  \vert
\leq ch^{-k},\ \ \ k=1,2. \label{4.26}%
\end{equation}
We also need the cut-off functions $\chi_{0j}$ and $\chi_{j}$ determined in
(\ref{4.27}), (\ref{4.l}).

As an approximation of the solution $u(  h,x)  $ of problem
(\ref{1.10})-(\ref{1.16}) with $\alpha=1,$ we take%
\begin{align}
\mathbf{u}_{0}^{\prime}(  h,x)   &  =X_{0}^{h}(  y)
U_{0}^{0}(  y,\ln h)  +{\textstyle\sum\nolimits_{j}} \chi_{0j}(  y)
W_{0j}^{0}(  \eta^{j},\ln h)  -X_{0}^{h}(  y)
{\textstyle\sum\nolimits_{j}} \chi_{0j}(  y)  U_{j}^{0}(  0,\ln h)
+\label{4.28}\\
&  +A_{j}(  \ln h)  (2\pi)^{-1}(\ln(h/r_{j}+\ln
c_{\log}(  \omega_{j}))  ,\nonumber\\
\mathbf{u}_{j}^{\prime}(  h,x)   &  =X_{j}^{h}(  z)
U_{j}^{0}(  z,\ln h)  +\chi_{j}(  z)  (  U_{j}%
^{0}(  0,\ln h)  +hW_{j}^{1}(  \xi^{j},\ln h)  )
-\label{4.29}\\
&  -X_{j}^{h}(  z)  \chi_{j}(  z)  (  U_{j}%
^{0}(  0,\ln h)  +z\partial_{z}U_{j}^{0}(  0,\ln h)
)  .\nonumber
\end{align}
Let us comment on these formulas where the asymptotic terms constructed in
Section \ref{sect3.1} are used. We display explicitly their dependence on $\ln
h$ but we skip it in further calculations though. In (\ref{4.28}) and
(\ref{4.29}), the main terms (\ref{2.25}), (\ref{3.A}) and $U_{j}^{0},$
$W_{j}^{0}+hW_{j}^{1}$ of the outer and inner expansions, respectively, are
inserted entirely so that the matched terms%
\begin{align}
S_{0j}(  y,\ln h)   &  =U_{j}^{0}(  0,\ln h)
+A_{j}(  \ln h)  (2\pi)^{-1}(\ln(h/r_{j}+\ln
c_{\log}(  \omega_{j}))  ,\label{4.30}\\
S_{j}(  z,\ln h)   &  =U_{j}^{0}(  0,\ln h)
+z\partial_{z}U_{j}^{0}(  0,\ln h) \nonumber
\end{align}
do appear twice, however the subtrahends compensate for this reduplication.
Cut-off functions are distributed in (\ref{4.28}) and (\ref{4.29}) in such a
way that in the sequel it is worth to make good use of the following relations
with commutators, see (\ref{2.com}),%
\begin{align}
[  \Delta_{x},X_{0}^{h}\chi_{0j}]   &  =\chi_{0j}[  \Delta
_{x},X_{0}^{h}]  +X_{0}^{h}[  \Delta_{x},\chi_{0j}]
=-[  \Delta_{x},\chi_{0j}^{h}]  +[  \Delta_{x},\chi
_{0j}]  ,\label{4.31}\\
[  \Delta_{x},X_{j}^{h}\chi_{j}]   &  =-[  \Delta_{x},\chi
_{j}^{h}]  +[  \Delta_{x},\chi_{j}]  ,\nonumber
\end{align}
which are readily apparent from definitions (\ref{4.24}), (\ref{4.25}) and
(\ref{4.26}), (\ref{4.l}). After commuting, the expressions $-[
\Delta_{x},\chi_{0j}^{h}]  S_{0j},$ $-[  \Delta_{x},\chi_{j}%
^{h}]  S_{j}$ and $[  \Delta_{x},\chi_{0j}]  S_{0j},$
$[  \Delta_{x},\chi_{j}]  S_{j}$ will be added to discrepancies
produced by outer and inner expansions, respectively, and this rearrangement
will assist in diminishing residuals.

To fulfil our plan, we insert (\ref{4.28}) into equation (\ref{1.10}) and, in
view of (\ref{4.31}), obtain%
\begin{align}
\Delta_{x}\mathbf{u}_{0}^{\prime}  &  =X_{0}^{h}(  y)  \Delta
_{x}U_{0}^{0}+{\textstyle\sum\nolimits_{j}} \chi_{0j}\Delta_{x}W_{0j}^{0}+X_{0}^{h}(
y)  {\textstyle\sum\nolimits_{j}} \chi_{0j}\Delta_{x}S_{0j}-\label{4.32}\\
&  -{\textstyle\sum\nolimits_{j}} [  \Delta_{x},\chi_{0j}^{h}]  (  U_{0}%
^{0}-S_{0j})  +{\textstyle\sum\nolimits_{j}} [  \Delta_{x},\chi_{0j}]
(  W_{0j}^{0}-S_{0j})  .\nonumber
\end{align}
We denote by $I_{1}^{0},...,I_{5}^{0}$ terms on the right of (\ref{4.32}) and,
by (\ref{2.21}) and (\ref{2.57}), immediately conclude that
$I_{1}^{0}=-X_{0}^{h}f_{0}^{0},\ I_{2}^{0}=0 \text{ and } I_{3}^{0}=0. $
The other two terms require an estimation but, due to the rearrangement
explained above, the differences $U_{0}^{0}-S_{0j}$ and $W_{0j}^{0}-S_{0j}$
become small on supports of the commutator coefficients. Indeed, comparing
(\ref{4.10}) and (\ref{2.25}), (\ref{2.18}) shows that the difference
\begin{equation}
T_{0j}(  y)  =U_{0}^{0}(  y)  -S_{0j}(  y)
\label{4.Tj}%
\end{equation}
loses the logarithmic term and, therefore, falls into
$H^{2}(  \omega_{0})  .$ By (\ref{3.10}) and
(\ref{4.20}), we conclude that%
\begin{equation}
T_{0j}(  P^{j})  =0,\ \ \Vert T_{0j};H^{2}(
\mathbb{B}_{R_{0}}(  P^{j})  )  \Vert \leq c\vert
\ln h\vert \mathcal{N}. \label{4.T}%
\end{equation}
Coefficients of the first-order differential operator $[  \Delta_{x}%
,\chi_{0j}^{h}]  ,$ see (\ref{2.com}), are located in the annulus
$\Upsilon_{j}^{h}=\{  y:R_{j}h\leq r_{j}\leq2R_{j}h\}  $ where the
variable $r_{j}$ is equivalent to $h.$ We derive the weighted estimate
\begin{align}
(  1+\vert \ln h\vert )  ^{2}&\Vert r(
1+\vert \ln r\vert )  I_{4}^{0};L^{2}(  \Omega_{\bullet
}(  h)  )  \Vert ^{2}\leq
\label{4.34}\\
\leq & c(  1+\vert \ln h\vert )  ^{2}\sum\nolimits_{j} %
\int_{0}^{h}\int_{\Upsilon_{j}^{h}}r_{j}^{2}(  1+\vert \ln
r_{j}\vert )  ^{2}(  h^{-2}\vert \nabla_{y}T_{0j}(
y)  \vert ^{2}+h^{-4}\vert T_{0j}(  y)
\vert ^{2})  dydz\leq\nonumber\\
\leq &chh^{2}(  1+\vert \ln h\vert )  ^{6}\sum\nolimits_{j}\int_{\Upsilon_{j}^{h}}r_{j}^{-2}(  1+\vert \ln r_{j}\vert
)  ^{-2}(  \vert \nabla_{y}T_{0j}(  y)  \vert
^{2}+r_{j}^{-2}\vert T_{0j}(  y)  \vert ^{2})
dy\leq\nonumber\\
\leq &ch^{3}(  1+\vert \ln h\vert )^{6}\sum\nolimits_{j}\Vert T_{j};H^{2}(  \mathbb{B}_{R_{0}}(  P^{j})
)  \Vert ^{2}\leq ch^{3}(  1+\vert \ln h\vert
)  ^{6}\mathcal{N}^{2}.\nonumber
\end{align}
Here, the weight $r_{j}(  1+\vert \ln r_{j}\vert )  $
arises from the first norm on the left-hand side of (\ref{4.13}) and,
recalling the equivalence $r_{j}\sim h$ in $\Upsilon_{j}^{h},$ we have changed
$r_{j}(  1+\vert \ln r_{j}\vert )  $ for $h(
1+\vert \ln h\vert )  (  1+\vert \ln r_{j}%
\vert )  ^{-1}$ as well as replace by $cr_{j}^{-1}$ and
$cr_{j}^{-2}$ the big bounds $ch^{-1}$ and $ch^{-2}$ for the coefficients of
the commutator, see (\ref{4.26}) and (\ref{2.com}) again. In the end of
calculation (\ref{4.34}) we have applied formula (\ref{4.T}) together with
estimate%
\begin{align}
\Vert r_{j}^{-2}(  1+\vert \ln r_{j}\vert )
^{-1}T_{0j};L^{2}(  \mathbb{B}_{R_{0}}(  P^{j})  )
\Vert  &  \leq c\Vert r_{j}^{-1}(  1+\vert \ln
r_{j}\vert )  ^{-1}\nabla_{y}T_{0j};L^{2}(  \mathbb{B}%
_{R_{0}}(  P^{j})  )  \Vert \leq\label{4.H}\\
&  \leq C\Vert \nabla_{y}T_{0j};H^{1}(  \mathbb{B}_{R_{0}}(
P^{j})  )  \Vert .\nonumber
\end{align}
The latter requires the above-mentioned relation $T_{0j}(  P^{j})
=0$ and is inherited from the following one-dimensional inequalities of
Hardy's type%
\begin{align}
\int_{0}^{R_{0}}r^{-3}\left(  1+\left\vert \ln\frac{r}{R_{0}}\right\vert
\right)  ^{-2}\vert T(  r)  \vert ^{2}dr  &  \leq
c\int_{0}^{R_{0}}r^{-1}\left(  1+\left\vert \ln\frac{r}{R_{0}}\right\vert
\right)  ^{-2}\left\vert \frac{dT}{dr}(  r)  \right\vert
^{2}dr,\ \ T(  0)  =0,\label{4.35}\\
\int_{0}^{R_{0}}r^{-1}\left(  1+\left\vert \ln\frac{r}{R_{0}}\right\vert
\right)  ^{-2}\left\vert \mathcal{T}(  r)  \right\vert ^{2}dr  &
\leq C\int_{0}^{R_{0}}r\left(  \left\vert \frac{d\mathcal{T}}{dr}(
r)  \right\vert ^{2}+\vert \mathcal{T}(  r)  \vert
^{2}\right)  dr.\nonumber
\end{align}
Both the inequalities are derived in a standard way.
% The
%second one looks like the standard Hardy inequality with logarithm
%(\ref{2.43}) but we have added $1$ to $\vert \ln\frac{r}{R_{0}%
%}\vert $ to avoid a singularity of the weight at $r=R_{0}$ and this
%modification has lead to the summand $\vert \mathcal{T}(  r)
%\vert ^{2}$ on the right.
Note that $T(  r)  $ and
$\mathcal{T}(  r)  $ in (\ref{4.35}) substitute for $T_{0j}(
y)  $ and $\nabla_{y}T_{0j}(  y)  $ so that integrating in
the angular variable $\varphi\in[  0,2\pi)  $ is needed, cf.
(\ref{2.20}) and (\ref{2.44}). We emphasize that (\ref{4.34}) is the only
cumbersome calculation in the section and there exist other ways to treat the
discrepancy term $I_{4}^{0}$ but we prefer to use one tool throughout the
paper, namely weighted inequalities of Hardy's type.

By (\ref{3.A}), (\ref{3.8}) and (\ref{4.30}), we have%
\begin{align}
W_{0j}^{0}(  h^{-1}(  y-P^{j})  )  -S_{0j}(
y)   &  =O(  hr_{j}^{-1}\vert A_{j}\vert )
,\label{4.36}\\
\nabla_{y}W_{0j}^{0}(  h^{-1}(  y-P^{j})  )  -\nabla
_{y}S_{0j}(  y)   &  =O(  hr_{j}^{-2}\vert A_{j}%
\vert )  .\nonumber
\end{align}
Since coefficients of the differential operator $[  \Delta_{x},\chi
_{0j}]  $ vanish in the disk $\mathbb{B}_{R_{0}/2}(  P^{j})
,$ see (\ref{4.27}), we obtain%
\[
(  1+\vert \ln h\vert )  ^{2}\Vert r(
1+\vert \ln r\vert )  I_{5}^{0};L^{2}(  \Omega_{\bullet
}(  h)  )  \Vert ^{2}\leq c(  1+\vert \ln
h\vert )  ^{2}hh^{2}\mathcal{N}^{2},
\]
where the first $h$ stands due to the integration in $z\in(  0,h)
$ while $h^{2}$ and $\mathcal{N}^{2}$ come from (\ref{4.36}) and (\ref{4.20}).

Let us now consider the discrepancy%
\begin{align}
\Delta_{x}\mathbf{u}_{j}^{\prime}  &  =X_{j}^{h}\partial_{z}^{2}U_{j}^{0}%
+\chi_{j}(  W_{j}^{0}+hW_{j}^{1})  +X_{j}^{h}\chi_{j}\Delta
_{x}S_{j}-\label{4.38}\\
&  -[  \Delta_{x},\chi_{j}^{h}]  (  U_{j}^{0}-U_{j}^{0}(
0)  -z\partial_{z}U_{j}^{0}(  0)  )  +h[
\Delta_{x},\chi_{j}]  (  W_{j}^{1}-\zeta\partial_{z}U_{j}%
^{0}(  0)  )  .\nonumber
\end{align}
Denoting terms on the right by $I_{1}^{j},...,I_{5}^{j}$ we immediately derive
from (\ref{2.6}) and (\ref{2.58}) that
$I_{1}^{j}=-X_{j}^{h}\gamma_{j}^{-1}f_{j}^{0},\ \ I_{2}^{j}=0$ and $I_{3}^{j}=0$. Since $\widetilde{U}_{j}^{0}(  z)  =U_{j}^{0}(  z)
-U_{j}^{0}(  0)  -z\partial_{z}U_{j}^{0}(  0)  $ belongs
to $H^{2}(  0,l_{j})  $ and satisfies $\widetilde{U}_{j}^{0}(
0)  =z\partial_{z}\widetilde{U}_{j}^{0}(  0)  =0,$ the Hardy
inequalities provide%
\[
\int_{0}^{d}z^{-4}\vert \widetilde{U}_{j}^{0}(  z)
\vert ^{2}dz\leq\frac{4}{9}\int_{0}^{d}z^{-2}\vert \partial
_{z}\widetilde{U}_{j}^{0}(  z)  \vert ^{2}dz\leq\frac{16}%
{9}\int_{0}^{d}\vert \partial_{z}^{2}\widetilde{U}_{j}^{0}(
z)  \vert ^{2}dz
\]
with any $d>0.$ Recalling (\ref{4.25}), (\ref{4.26}) and (\ref{4.20}), we
write
\begin{gather}
h\Vert I_{4}^{j};L^{2}(  \Omega_{j}(  h)  )
\Vert ^{2}\leq chh^{2}\int_{h}^{2h}(  h^{-4}\vert
\widetilde{U}_{j}^{0}(  z)  \vert ^{2}+h^{-2}\vert
\partial_{z}\widetilde{U}_{j}^{0}(  z)  \vert ^{2})
dz\leq\label{4.40}\\
\leq ch^{3}\int_{h}^{2h}(  z^{-4}\vert \widetilde{U}_{j}^{0}(
z)  \vert ^{2}+z^{-2}\vert \partial_{z}\widetilde{U}_{j}%
^{0}(  z)  \vert ^{2})  dz\leq ch^{3}\mathcal{N}%
^{2}.\nonumber
\end{gather}
Notice that the first factor $h$ in (\ref{4.40}) is in accordance with
$h^{-\alpha}=h^{-1}$ in (\ref{4.13}) and the factor $h^{2}$ is caused by the
integration over the small cross-section $\omega_{j}^{h}\ni y.$ Furthermore,
the exponential decay of the difference $W_{j}^{1}(  \eta^{j}%
,\zeta)  -\zeta\partial_{z}U_{j}^{0}(  0)  $ as
$\zeta\rightarrow+\infty,$ see (\ref{3.6}) and (\ref{2.61}), and the location
of supports of coefficients in $[  \Delta_{x},\chi_{j}]  ,$ see
(\ref{4.l}), assure that
\begin{equation}
h\Vert I_{5}^{j};L^{2}(  \Omega_{j}(  h)  )
\Vert ^{2}\leq ce^{-\delta/h}\mathcal{N}^{2}. \label{4.41}%
\end{equation}
The residuals%
\begin{equation}
\mathbf{f}_{0}(  h,x)  =-\Delta_{x}\mathbf{u}_{0}(
h,x)  -X_{0}^{h}(  y)  f_{0}^{0}(  y)
,\ \ \ \mathbf{f}_{j}(  h,x)  =-h^{-1}\gamma_{j}\Delta
_{x}\mathbf{u}_{j}(  h,x)  -h^{-1}X_{j}^{h}(  z)
f_{j}^{0}(  z)  \nonumber
\end{equation}
in (\ref{1.10}) and (\ref{1.11}) are estimated. By definition, functions (\ref{4.28}) and (\ref{4.29}) meet the
homogeneous boundary conditions (\ref{1.10})-(\ref{1.14}) and the second
transmission condition (\ref{1.16}). However, the discrepancy%
\[
\mathbf{u}_{j}^{\prime}(  h,x)  -\mathbf{u}_{0}^{\prime}(
h,x)  =hW_{j}^{1}(  \xi^{j})  ,\ \ x\in\upsilon_{j}^{h},
\]
is left in the first transmission condition (\ref{1.15}) which can be
compensated by the function%
\begin{equation}
\mathbf{u}_{0j}^{\prime}(  h,x)  =hA_{j}W_{0j}^{1}(  \xi
^{j})  \label{4.43}%
\end{equation}
where $W_{0j}^{1}$ has a support in $(  \overline{\mathbb{B}_{R_{j}}%
}\setminus\omega_{j})  \times[  0,1]  $ and%
\begin{align}
\partial_{\nu}W_{0j}^{1}(  \xi^{j})   &  =0,\ \ \xi^{j}\in
\upsilon_{j}^{1},\ \ \ \ \ \ \ W_{0j}^{1}(  \xi^{j})
=0,\ \ \xi^{j}\in\partial\mathbb{B}_{R_{j}}\times(  0,1)
,\nonumber\\
\partial_{\zeta}W_{0j}^{1}(  \eta^{j},0)   &  =\partial_{\zeta
}W_{0j}^{1}(  \eta^{j},1)  =0,\ \ \eta^{j}\in\mathbb{B}_{R_{j}%
}\setminus\overline{\omega}_{j}.\nonumber
\end{align}
Furthermore, we have%
\begin{align}
(  1+\vert \ln h\vert )  ^{2}&\Vert r_{j}(
1+\vert \ln r_{j}\vert )  \nabla_{x}^{2}\mathbf{u}%
_{0j}^{\prime};L^{2}(  \Omega_{\bullet}(  h) )
\Vert ^{2}\leq\label{4.45}\\
\leq &ch^{2}\vert A_{j}\vert ^{2}h^{2}(  1+\vert \ln
h\vert )  ^{4}\Vert \Delta_{x}W_{0j}^{1};L^{2}(
\Omega_{\bullet}(  h)  )  \Vert ^{2}\leq ch^{3}(
1+\vert \ln h\vert )  ^{4}\mathcal{N}^{2}\nonumber
\end{align}
because the last norm gets order $h^{-1/2}.$

The function%
\begin{equation}
\mathbf{u}(  h,x)  =\left\{
\begin{array}
[c]{c}%
\mathbf{u}_{0}^{\prime}(  h,x)  +%
%TCIMACRO{\dsum _{j=1}^{J}}%
%BeginExpansion
{\displaystyle{\textstyle\sum\nolimits_{j}} }
%EndExpansion
\mathbf{u}_{0j}^{\prime}(  h,x)  ,\ \ \ x\in\Omega_{\bullet}(
h)  ,\\
\mathbf{u}_{j}^{\prime}(  h,x)  ,\ \ \ x\in\Omega_{j}(
h)
\end{array}
\right.  \label{4.46}%
\end{equation}
belongs to the space $H_{0}^{1}(  \Xi(  h)  ,\Gamma(
h)  )  $ and verifies the integral identity%
\begin{equation}
a(  \mathbf{u},v)  =(  f_{0}^{0},v_{0})  _{\Omega
_{\bullet}(  h)  }+h^{-1}{\textstyle\sum\nolimits_{j}} (  f_{j}^{0}%
,v_{j})  _{\Omega_{j}(  h)  }+(  \mathbf{f}^{\prime
},v)  _{\Xi(  h)  } \label{4.47}%
\end{equation}
where $\mathbf{f}_{0}^{\prime} =(  1-X_{0}^{h})  f_{0}^{0}, \quad
\mathbf{f}_{j}^{\prime}  =h^{-1}(  1-X_{j}^{h})  f_{j}%
^{0}+\mathbf{f}_{j}$ and
\begin{equation}
\vert (  \mathbf{f},v)  _{\Xi(  h)  }\vert
\leq ch^{3/2}(  (  1+\vert \ln h\vert )
^{3}+h^{-1/2})  \mathcal{N}a(  v,v;\Xi(  h)  )
^{1/2}. \label{4.49}%
\end{equation}
Let us comment. The terms $(  1-X_{0}^{h})
f_{0}^{0}$ and $(  1-X_{j}^{h})  f_{j}^{0}$ were estimated
according to (\ref{4.41}) and
\begin{align}
\vert (  (  1-X_{0}^{h})  f_{0}^{0},v_{0})
_{\Omega_{\bullet}(  h)  }\vert \leq &c\Vert r(  1+\vert \ln r\vert )  f_{0}%
^{0};L^{2}(  \text{supp}(  1-X_{0}^{h})  )  \Vert
\Vert r^{-1}(  1+\vert \ln r\vert )  ^{-1}%
v_{0};L^{2}(  \Omega_{\bullet}(  h)  )  \Vert
\leq\nonumber\\
\leq &ch(  1+\vert \ln h\vert )  ^{2}h^{1/2}%
\mathcal{N}a(  v,v;\Xi(  h)  )  ^{1/2},\label{4.50}\\
h^{-1}\vert (  (  1-X_{j}^{h}) & f_{j}^{0},v_{j})
_{\Omega_{j}(  h)  }\vert \leq ch^{-1}\Vert f_{j}%
^{0};L^{2}(  \text{supp}(  1-X_{j}^{h})  )  \Vert
\Vert v_{j};L^{2}(  \text{supp}(  1-X_{j}^{h})  )
\Vert \leq\label{4.51}\\
\leq &ch^{-1}h\mathcal{N}h^{1/2}\Vert \partial_{z}v_{j};L^{2}(
\Omega_{j}(  h)  )  \Vert \leq ch\mathcal{N}a(
v,v;\Xi(  h)  )  ^{1/2}.\nonumber
\end{align}
Here, we took into account that%
\[
\text{supp}(  1-X_{0}^{h})     =%
%TCIMACRO{\dbigcup _{j=1}^{J}}%
%BeginExpansion
{\displaystyle\bigcup_{j=1}^{J}}
%EndExpansion
\text{supp}\chi_{j}^{h}\subset%
%TCIMACRO{\dbigcup _{j=1}^{J}}%
%BeginExpansion
{\displaystyle\bigcup_{j=1}^{J}}
%EndExpansion
\overline{\mathbb{B}_{hR_{j}}(  P^{j})  }\times[  0,h],\quad
\text{supp}(  1-X_{j}^{h})   =\text{supp}\chi_{j}^{h}%
\subset\overline{\omega_{j}^{h}}\times[  0,2h]  .
\]
In (\ref{4.50}), we also made the change $r(  1+\vert \ln
r\vert )  \mapsto h(  1+\vert \ln h\vert )
$ and attached the factor $h^{1/2}$ because the norm in $L^{2}(
\Omega_{\bullet}(  h)  )  $ expects integration in
$z\in(  0,h)  .$ Calculation (\ref{4.51}) is based on (\ref{4.13})
and (\ref{1.20}), (\ref{1.17}) together with the obvious inequality%
\[
\int_{0}^{2h}\vert V(  z)  \vert ^{2}dz\leq2hl_{j}%
\int_{0}^{l_{j}}\vert \partial_{z}V(  z)  \vert
^{2}dz,\ \ V(  l_{j})  =0.
\]
%originating in the Newton-Leibnitz formula $V(  z)  =-%
%%TCIMACRO{\dint _{z}^{l_{j}}}%
%%BeginExpansion
%{\displaystyle\int_{z}^{l_{j}}}
%%EndExpansion
%\partial_{z}V(  \mathfrak{z})  ^{2}d\mathfrak{z}.$

We insert into the integral identities (\ref{1.19}) and (\ref{4.47}) the
difference%
\begin{equation}
v=u-\mathbf{u}\in H_{0}^{1}(  \Xi(  h)  ,\Gamma(
h)  )  \label{4.52}%
\end{equation}
between the exact and approximate solutions of problem (\ref{1.10}%
)-(\ref{1.16}). Subtracting one identity from the other yields%
\begin{equation}
a(  v,v;\Xi(  h)  )  =(  f_{00}^{\bot}%
,v_{0})  _{\Omega_{\bullet}(  h)  }+h^{-1}{\textstyle\sum\nolimits_{j}}(  f_{j0}^{\bot},v_{j})  _{\Omega_{j}(  h)
}+(  \mathbf{f}^{\prime},v)  _{\Xi(  h)  } \label{4.53}%
\end{equation}
The orthogonality conditions (\ref{2.9}) and (\ref{2.2}) allow us to replace
in the scalar products the functions $v_{0}$ and $v_{j}$ by $v_{0}%
-\overline{v}_{0}$ and $v_{j}-\overline{v}_{j}$ where $\overline{v}_{0}$ and
$\overline{v}_{j}$ are the mean-value functions (\ref{4.6}) and (\ref{4.3}),
respectively. Now the Poincar\'{e} inequalities (\ref{4.7}) and (\ref{4.4})
provide the estimates%
\begin{align}
\vert (  f_{00}^{\bot},v_{0})  _{\Omega_{\bullet}(
h)  }\vert  &  =\vert (  f_{00}^{\bot},v_{0}%
-\overline{v}_{0})  _{\Omega_{\bullet}(  h)  }\vert
\leq ch^{1/2}\mathcal{N}h\Vert \partial_{z}v_{0};L^{2}(
\Omega_{\bullet}(  h)  )  \Vert \leq\label{4.54}\\
&  \leq ch^{3/2}\mathcal{N}a(  v,v;\Xi(  h)  )
^{1/2},\nonumber\\
h^{-1}\vert (  f_{j0}^{\bot},v_{j})  _{\Omega_{j}(
h)  }\vert  &  =h^{-1}\vert (  f_{j0}^{\bot}%
,v_{j}-\overline{v}_{0j})  _{\Omega_{j}(  h)  }\vert
\leq ch^{-1}h\mathcal{N}h\Vert \nabla_{y}v_{j};L^{2}(  \Omega
_{j}(  h)  )  \Vert \leq\nonumber\\
&  \leq ch^{3/2}\mathcal{N}a(  v,v;\Xi(  h)  )
^{1/2}.\nonumber
\end{align}
In view of Theorem \ref{th4.1}, applying (\ref{4.49}) and (\ref{4.54}) to
(\ref{4.53}) leads to the following assertion.

\begin{proposition}
\label{prop4.N1}Under assumptions (\ref{4.16}) and (\ref{4.F}), the solution
$u$ of problem (\ref{1.10})-(\ref{1.16}) with $\alpha=1$ and its global
approximation (\ref{4.46}) constructed in Section \ref{sect3.1} are in the
relationship%
\begin{gather}
\Vert \nabla_{x}(  u_{0}-\mathbf{u}_{0})  ;L^{2}(
\Omega_{\bullet}(  h)  )  \Vert +(  1+\vert
\ln h\vert )  ^{-1}\Vert r^{-1}(  1+\vert \ln
r\vert )  ^{-1}(  u_{0}-\mathbf{u}_{0})  ;L^{2}(
\Omega_{\bullet}(  h)  )  \Vert +\label{4.55}\\
+h^{-1}{\textstyle\sum\nolimits_{j}} (  \Vert \nabla_{x}(  u_{j}-\mathbf{u}%
_{j})  ;L^{2}(  \Omega_{j}(  h)  )  \Vert
+\Vert (  l_{j}-z)  (  u_{j}-\mathbf{u}_{j})
;L^{2}(  \Omega_{j}(  h)  )  \Vert )  \leq
ch\mathcal{N},\nonumber
\end{gather}
where $\mathcal{N}$ is the sum of norms of functions (\ref{4.F}), $c$ is a
constant independent of $h\in(  0,h_{0}]  $ and $F_{0}%
,F_{1},...,F_{J}$ are given in (\ref{4.16}).
\end{proposition}

\subsection{The asymptotics in the case $\alpha=1.$\label{sect4.5}}

The complicated structures (\ref{4.28}), (\ref{4.29}) and (\ref{4.46}) were
introduced with a technical reason only. After proving estimate (\ref{4.45})
we easily simplify the final asymptotic structures.

In the rod $\Omega_{j}(  h)  $, expression (\ref{4.29}) differs
from the intact solution (\ref{3.U}) of the limit problem (\ref{2.6}),
(\ref{2.7}), (\ref{3.1}), (\ref{3.G}) by two terms%
\[
\chi_{j}^{h}(  z)  (  U_{j}^{0}(  z,\ln h)
-U_{j}^{0}(  0,\ln h)  -z\partial_{z}U_{j}^{0}(  0,\ln
h)  )  \text{ \ \ and \ \ }h\chi_{j}(  z)
\widetilde{W}_{j}^{1}(  \xi^{j},\ln h)
\]
where $\widetilde{W}_{j}^{1}(  \xi^{j},\ln h)  =W_{j}^{1}(  \xi
^{j},\ln h)  -\zeta\partial_{z}U_{j}^{0}(  0,\ln h)$
stands for exponentially decaying remainder in the asymptotic form
(\ref{3.6}), cf. Lemma \ref{lem2.7}. A direct calculation shows that the norm
$h^{-1/2}\Vert \cdot;H^{1}(  \Omega_{j}(  h)  )
\Vert $ of both the terms does not exceed $ch\mathcal{N},$ that is the
bound in (\ref{4.55}), and thus they can be neglected. In this way we derive
from Proposition \ref{prop4.N1} the estimate%
\begin{equation}
h^{-1/2}\Vert u_{j}-U_{j}^{0};H^{1}(  \Omega_{j}(  h)
)  \Vert \leq ch\mathcal{N}. \label{4.57}%
\end{equation}
It should be stressed that $h^{-1/2}\Vert U_{j}^{0};H^{1}(
\Omega_{j}(  h)  )  \Vert =O(  h^{1/2})  $
and, therefore, inequality (\ref{4.57}) indeed exhibits an asymptotics of
$u_{j}.$

In the perforated plate (\ref{1.5}) an asymptotic form of the solution $u_{0}$
is much more complicated. Function (\ref{4.28}) differs from the singular
solution (\ref{2.25}) of the limit problem (\ref{2.21}), (\ref{2.14}) by the
sum of the terms $\chi_{0j}^{h}(  y)  T_{0j}(  y,\ln h)$ and $\chi_{0j}(  y)  A_{j}(  \ln h)  \widetilde
{\mathbf{W}}^{j}(  \eta^{j})  ,\ \ \ j=1,...,J$, where $T_{0j}$ is defined in (\ref{4.Tj}) and $\widetilde{\mathbf{W}}^{j}$ is
the remainder in the asymptotic form (\ref{2.60}) of the logarithmic
potential. By Hardy's type inequality (\ref{4.H}), we, similarly to
(\ref{4.34}), derive the estimate%
\[
\Vert \nabla_{x}(  \chi_{0j}^{h}T_{0j})  ;L^{2}(
\Omega_{\bullet}(  h)  )  \Vert \leq ch^{3/2}(
1+\vert \ln h\vert )  ^{3}\mathcal{N}%
\]
which permits to omit $\chi_{0j}^{h}T_{0j}$ in the asymptotics of $u_{0},$ cf.
the bound in (\ref{4.55}). However, based on the decay rates $O(  (
1+\rho)  ^{-1})  $ and $O(  (  1+\rho)
^{-2})  $ of $\widetilde{\mathbf{W}}^{j}$ and $\nabla_{\eta}%
\widetilde{\mathbf{W}}^{j}$ respectively, we derive that%
\begin{align}
\Vert \nabla(  \chi_{0j}\widetilde{\mathbf{W}}^{j})
;L^{2}(  \Omega_{\bullet}(  h)  )  \Vert  &  \leq
ch\int_{0}^{R_{0}}\left(  \frac{1}{(  1+r/h)  ^{2}}+\frac{1}{h^{2}%
}\frac{1}{(  1+r/h)  ^{4}}\right)  rdr\leq\label{4.59}\\
&  \leq ch(  h^{2}(  1+\vert \ln h\vert )
+1)  \leq ch,\nonumber\\
\Vert \chi_{0j}\widetilde{\mathbf{W}}^{j};L^{2}(  \Omega_{\bullet
}(  h)  )  \Vert  &  \leq ch\int_{0}^{R_{0}}\frac
{r}{(  1+r/h)  ^{2}}dr\leq ch^{3}(  1+\vert \ln
h\vert )  .\nonumber
\end{align}
This means that the Dirichlet norm of the boundary layer term
\begin{equation}
A_{j}(  \ln h)  \widetilde{\mathbf{W}}^{j}(  \eta^{j})
=A_{j}(  \ln h)  (  \mathbf{W}^{j}(  \eta^{j})
+(2\pi)^{-1}(\ln\rho_{j}-\ln c_{\log}(  \omega
_{j})  )  )  \label{4.60}%
\end{equation}
as well as the first weighted norm on the left of (\ref{4.13}) get the same
order in $h$ as the singular solution $u_{0}(  h,x)  $ so that
$\chi_{0j}A_{j}\widetilde{\mathbf{W}}^{j}$ must be kept in the asymptotics.
Besides, functions (\ref{4.43}) with small supports which have been added in
(\ref{4.46}) can be neglected according to estimate (\ref{4.45}).

Let us formulate the obtained theorem on asymptotics.

\begin{theorem}
\label{th4.N3}Under conditions (\ref{4.16}) and (\ref{4.F}) the restriction
$u_{j}(  h,x)  $ on the rod $\Omega_{j}(  h)  $ of the
solution $u(  h,x)  $ of problem (\ref{1.10})-(\ref{1.16}) meets
the asymptotic formula (\ref{4.57}) where $U_{j}^{0}$ is given in (\ref{3.U}).
The restriction $u_{0}(  h,x)  =u(  h,x)  |_{\Omega
_{\bullet}(  h)  }$ on the perforated plate (\ref{1.5}) satisfies
the estimate%
\begin{align}
\Vert \nabla_{x}( & u_{0}-U_{0}^{0}-{\textstyle\sum\nolimits_{j}} \chi_{0j}%
A_{j}\widetilde{\mathbf{W}}^{j})  ;L^{2}(  \Omega_{\bullet}(
h)  )  \Vert +\label{4.61}\\
+&(  1+\vert \ln h\vert )  ^{-1}\Vert r^{-1}(
1+\vert \ln r\vert )  ^{-1}(  u_{0}-U_{0}^{0}-{\textstyle\sum\nolimits_{j}}\chi_{0j}A_{j}\widetilde{\mathbf{W}}^{j})  ;L^{2}(
\Omega_{\bullet}(  h)  )  \Vert \leq ch\mathcal{N}%
\nonumber
\end{align}
where $U_{0}^{0}$ is the linear combination (\ref{2.25}) with the coefficients
$A_{p}(  \ln h)  $ computed in (\ref{3.16}), (\ref{3.15}) and
$A_{j}(  \ln h)  \widetilde{\mathbf{W}}^{j}(  h^{-1}(
y-P^{j})  )  $ are the boundary layer terms (\ref{4.60}).
\end{theorem}

Since, by virtue of (\ref{3.12}),%
\[
M(  \ln h)  ^{-1}=2\pi\vert \ln h\vert ^{-1}%
\mathbb{I}+O(  \vert \ln h\vert ^{-2})  ,\ \ \ m(
\ln h)  =2\pi J\vert \ln h\vert ^{-1}+O(  \vert \ln
h\vert ^{-2})  ,
\]
formulas (\ref{3.16}), (\ref{3.15}) and (\ref{3.14}) indicate%
\begin{align}
A_{0}(  \ln h)   &  =A_{0}^{(  -1)  }\vert \ln
h\vert +A_{0}^{(  0)  }+O(  \vert \ln h\vert
^{-2})  ,\ \ \ A_{j}(  \ln h)  =A_{j}^{(  0)
}+O(  \vert \ln h\vert ^{-1})  ,\label{4.62}\\
A_{0}^{(  -1)  }  &  =\frac{1}{2\pi J}\int_{\omega_{0}}f_{0}%
^{0}(  y)  dy,\ \ \ A_{j}^{(  0)  }=-\frac{1}{J}%
\int_{\omega_{0}}f_{0}^{0}(  y)  dy.\nonumber
\end{align}

\begin{corollary}
\label{cor4.N4}Under conditions (\ref{4.16}) and (\ref{4.F}) the following
convergences occur%
\begin{align}
u(  h,P^{j}+h\eta^{j},z)   &  \rightarrow U_{j}^{\#}(
z)  +A_{j}^{(  0)  }\gamma_{j}^{-1}\vert
\omega_{j}\vert^{-1}(l_{j}-z)\text{ \ strongly in }H^{1}(  \omega_{j}%
\times(  0,l_{j})  )  ,\label{4.63}\\
\vert \ln h\vert ^{-1}u(  h,y,h\zeta)   &  \rightarrow
A_{0}^{(  -1)  }\text{ \ strongly in }H^{1}(  \omega_{0}%
\times(  0,1)  )  ,\label{4.64}\\
u(  h,y,h\zeta)   &  \rightarrow U_{\bot}^{0}+A_{0}^{(
0)  }+{\textstyle\sum\nolimits_{j}} A_{j}^{(  0)  }G_{j}(  y)
\text{ \ strongly in }H^{1}(  \omega_{0}\times(  0,1)
)  , \label{4.65}%
\end{align}
where $U_{j}^{\#}\in H^{2}(  0,l_{j})  $ is a solution of problem
(\ref{3.P}), (\ref{3.PN}), $U_{\bot}^{0}\in H^{2}(  \omega_{0})  $
is described in Proposition \ref{prop2.1} and $A_{p}^{(  q)  }$ are
shown in (\ref{4.62}). The convergence rate in (\ref{4.63}) and (\ref{4.65})
is of order $\vert \ln h\vert ^{-1}$ and in (\ref{4.64}) of order
$\vert \ln h\vert ^{-1/2}.$
\end{corollary}

\textbf{Proof.} First of all, we observe that $\Vert G_{j};H^{1}(  \omega_{\bullet})  \Vert =O(\vert \ln h\vert ^{1/2})$ according to (\ref{2.18}).
Moreover, $\Vert U_{j}^{0};H^{1}(  \omega_{j}^{h}\times(  0,h)
)  \Vert \leq ch\Vert U_{j}^{0};H^{1}(  0,l_{j})
\Vert$. These inequalities together with (\ref{4.59}) and (\ref{4.62}) show that after
multiplication of $u|_{\Omega_{0}(  h)  }$ with $\vert \ln
h\vert ^{-1}$ the Sobolev norms of all asymptotic terms with exception
of $\vert \ln h\vert ^{-1}A_{0}(  \ln h)  \rightarrow
A_{0}^{(  -1)  }$ become $O(  \vert \ln h\vert
^{-1/2})  $ and, therefore, (\ref{4.64}) is proved. Other formulas
follow from a simple analysis of asymptotic terms involved into the error
estimates (\ref{4.57}) and (\ref{4.61}). If (\ref{2.15}) occurs, the singular
solution (\ref{2.25}) with $A_{p}=A_{p}^{(  0)  }$ as in
(\ref{4.65}) lives outside $H^{1}(  \omega_{0})  $ and thus
convergence (\ref{4.65}) cannot hold true in $H^{1}(  \omega_{0}%
\times(  0,1)  )  .$ $\blacksquare$

\subsection{The asymptotics in the case $\alpha=0.$\label{sect4.6}}

The justification scheme stays the same as in Section \ref{sect4.4}, however,
by virtue of the small factor $h=\min\{  h^{-\alpha+1},(
1+\vert \ln h\vert )  ^{-2}\}  $ on the first norm in
the weighted anisotropic inequality (\ref{4.13}) at $\alpha=0,$ final error
estimates differ from ones in Theorem \ref{th4.N3} with $\alpha=1.$ We further
outline calculations which lead to asymptotic formulas for $u_{p}$ in Theorem
\ref{th4.M3} serving for the case $\alpha=0.$

Using cut-off functions in (\ref{4.25}), (\ref{4.26}) and copying asymptotic
structures in (\ref{4.28}), (\ref{4.29}), we set%
\begin{align}
\mathbf{u}_{0}(  h,x)   &  =X_{0}^{h}(  y)  (
h^{-1}a_{0}+U_{0}^{0}(  y,\ln h)  )  +{\textstyle\sum\nolimits_{j}} \chi
_{0j}(  y)  (  h^{-1}a_{0}+w_{j}^{0}(  \eta^{j}%
,\zeta)  )  -\label{4.66}\\
&  -X_{0}^{h}(  y)  {\textstyle\sum\nolimits_{j}} \chi_{0j}(  y)
S_{0j}(  h,y,\zeta)  ,\nonumber\\
\mathbf{u}_{j}(  h,x)   &  =X_{j}^{h}(  z)  (
h^{-1}a_{0}(  1-z/l_{j})  +U_{j}^{0}(  y,\ln h)
)  +\chi_{j}(  z)  (  h^{-1}a_{0}+w_{j}^{0}(
\eta^{j},\zeta)  )  -\label{4.67}\\
&  -X_{j}^{h}(  z)  \chi_{j}(  z)  S_{j}(
h,\zeta) \nonumber
\end{align}
where entries of the asymptotic expansions (\ref{3.34}), (\ref{3.35}) and
(\ref{3.36}) constructed in Section \ref{sect3.3} are involved together with
the following terms subject to matching:%
\[
S_{0j}^{0}(h,y)=h^{-1}a_{0}-a_{0}\gamma_{j}\vert
\omega_{j}\vert (2\pi l_{j})^{-1}\ln(h/r_{j})+b_{j},\quad
S_{j}(  h,\zeta)=h^{-1}a_{0}-a_{0}l_{j}^{-1}(
\zeta+\gamma_{j}\vert \omega_{j}\vert \mathbf{q}_{j})
+b_{j},
\]
The constants $b_{j}$ in (\ref{3.45}) and solutions (\ref{3.47})
of the limit problem (\ref{2.6}), (\ref{2.7}), (\ref{3.46}) have been
determined up to the additive constant $A_{0}$ which now is chosen arbitrarily
but later will be fixed properly.

Since term (\ref{3.40}) of the inner expansion (\ref{3.36}) verifies the
homogeneous limit problem (\ref{2.28})-(\ref{2.33}), functions (\ref{4.66})
and (\ref{4.67}) meet both conditions (\ref{1.15}) and
(\ref{1.16}). The boundary conditions (\ref{1.12})-(\ref{1.14}) are satisfied
as well. We thus need only to examine discrepancies in the equations (\ref{1.10}) and (\ref{1.11}). Since the constant $h^{-1}a_{0}$ and
the linear functions (\ref{3.l}) are eliminated by the Laplace operator, we
have%
\begin{gather*}
\Delta_{x}\mathbf{u}_{0}=X_{0}^{h}\Delta_{x}U_{0}^{0}+{\textstyle\sum\nolimits_{j}} \chi
_{0j}\Delta_{x}w_{j}^{0}+X_{0}^{h}{\textstyle\sum\nolimits_{j}} \chi_{0j}\Delta_{x}S_{0j}%
-{\textstyle\sum\nolimits_{j}} [  \Delta_{x},\chi_{0j}^{h}]  (  h^{-1}%
a_{0}+U_{0}^{0}-S_{0j})  +\\
+{\textstyle\sum\nolimits_{j}} [  \Delta_{x},\chi_{0j}]  (  h^{-1}a_{0}%
+w_{j}^{0}-S_{0j})  =:I_{1}^{0}+I_{2}^{0}+I_{3}^{0}+I_{4}^{0}+I_{5}^{0}.
\end{gather*}
Clearly, $I_{1}^{0}=-X_{0}^{h}f_{0}^{0}$, and $I_{2}^{0}=I_{3}^{0}=0$. The difference $T_{0j}=h^{-1}a_{0}+U_{0}^{0}-S_{0j}$ satisfies the formulas
$T_{0j}(  P^{j})  =0$ and $\Vert T_{0j};H^{2}(  \mathbb{B}_{R_{0}}(  P^{j})  )  \Vert\leq c\mathcal{N}$,
cf. (\ref{4.Tj}) and (\ref{4.T}). Hence, the Hardy inequality (\ref{4.H})
ensures that
\begin{align*}
h^{-1}&\Vert r(  1+\vert \ln r\vert )  I_{4}^{0};L^{2}(  \Omega_{\bullet}(  h)  )  \Vert
^{2}\leq \\
\leq &ch^{2}(  1+\vert \ln h\vert )  ^{4}{\textstyle\sum\nolimits_{j}}\int_{\Upsilon_{j}^{h}}r_{j}^{-2}(  1+\vert \ln r_{j}\vert
)  ^{-2}(  \vert \nabla_{y}T_{0j}(  y)  \vert
^{2}+r_{j}^{-2}\vert T_{0j}(  y)  \vert ^{2})
dy\leq ch^{2}(  1+\vert \ln h\vert )  ^{4}\mathcal{N}^{2}.
\end{align*}
The differences $h^{-1}a_{0}+w_{j}^{0}-S_{0j}$ get properties (\ref{4.36})
and, therefore, $h^{-1}\Vert r(  1+\vert \ln r\vert )  I_{5}^{0};L^{2}(  \Omega_{\bullet}(  h)  )  \Vert
^{2}\leq ch^{2}\mathcal{N}^{2}$. We now consider the expression%
\begin{align*}
\Delta_{x}\mathbf{u}_{j}=&X_{0}^{h}\partial_{z}^{2}U_{j}^{0}+\chi_{j}\Delta
_{x}w_{j}^{0}+X_{j}^{h}\chi_{j}\Delta_{x}S_{j}-[  \Delta_{x},\chi_{j}%
^{h}]  (  U_{j}^{0}-U_{j}^{0}(  0)  -z\partial_{z}%
U_{j}^{0}(  0)  )  +\\
&-a_{0}l_{j}^{-1}\gamma_{j}\vert \omega_{j}\vert [  \Delta
_{x},\chi_{j}]  (  \mathbf{w}^{j}+(  2\pi)  ^{-1}\ln
\rho_{j}-\mathbf{q}_{j})  =:I_{1}^{j}+I_{2}^{j}+I_{3}^{j}+I_{4}%
^{j}+I_{5}^{j}%
\end{align*}
and similarly to (\ref{4.38})-(\ref{4.41}), we obtain $I_{1}^{j}=-X_{j}^{h}\gamma_{j}^{-1}f_{j}^{0}$, $I_{2}^{j}=I_{3}^{j}=0$ and
\begin{align*}
\Vert I_{4}^{j};L^{2}(  \Omega_{j}(  h)  )
\Vert ^{2}  &  \leq ch^{2}\int_{h}^{2h}(  h^{-2}\vert
\partial_{z}U_{j}^{0}(  z)  -\partial_{z}U_{j}^{0}(  0)
\vert ^{2}+h^{-4}\vert U_{j}^{0}(  z)  -U_{j}^{0}(
0)  -z\partial_{z}U_{j}^{0}(  0)  \vert ^{2})
dz\leq\\
&  \leq ch^{2}\int_{0}^{2h}\vert \partial_{z}^{2}U_{j}^{0}(
z)  \vert ^{2}dz\leq ch^{2}\mathcal{N}^{2}, \qquad\quad
\Vert I_{5}^{j};L^{2}(  \Omega_{j}(  h)  )
\Vert ^{2}  \leq ce^{-\delta/h}\mathcal{N}^{2}.
\end{align*}
Setting
\[
\mathbf{f}_{0}(  h,x)=-\Delta_{x}\mathbf{u}_{0}(
h,x)  -X_{0}^{h}(  y)  f_{0}^{0}(  y)  , \quad
\mathbf{f}_{j}(  h,x)=-\gamma_{j}\Delta_{x}\mathbf{u}%
_{j}(  h,x)  -X_{j}^{h}(  z)  f_{j}^{0}(  z),
\]
we repeat an argument from Section \ref{sect4.3} and conclude that the
difference (\ref{4.52}) between the exact and approximate solutions of problem
(\ref{1.10})-(\ref{1.16}) satisfies the formula%
\[
a(  v,v;\Xi(  h)  )  =(  f_{00}^{\bot}%
,v_{0})  _{\Omega_{\bullet}(  h)  }-(  \mathbf{f}%
_{0},v_{0})  _{\Omega_{\bullet}(  h)  }+{\textstyle\sum\nolimits_{j}} (
(  f_{j0}^{\bot},v_{j})  _{\Omega_{j}(  h)  }-(
\mathbf{f}_{j},v_{0})  _{\Omega_{j}(  h)  })  .
\]
Collecting the above estimates, repeating calculation (\ref{4.54}) and
applying inequalities (\ref{4.13}), (\ref{4.7}), (\ref{4.4}) yield%
\begin{align}
\vert (  \mathbf{f}_{0},v_{0})  _{\Omega_{\bullet}(
h)  }\vert  &  \leq h^{-1/2}\Vert r(  1+\vert \ln
r\vert )  \mathbf{f}_{0};L^{2}(  \Omega_{\bullet}(
h)  )  \Vert h^{1/2}\Vert r^{-1}(  1+\vert
\ln r\vert )  ^{-1}v_{0};L^{2}(  \Omega_{\bullet}(
h)  )  \Vert \leq\label{est}\\
&  \leq ch(  1+\vert \ln h\vert )  ^{2}\mathcal{N}%
a(  v,v;\Xi(  h)  )  ^{1/2},\nonumber\\
\vert (  \mathbf{f}_{j},v_{j})  _{\Omega_{j}(  h)
}\vert  & \leq ch\mathcal{N}a(  v,v;\Xi(  h)  )
^{1/2},\
\vert (  f_{00}^{\bot},v_{0})  _{\Omega_{\bullet}(
h)  }\vert    \leq ch\mathcal{N}a(  v,v;\Xi(  h)
)  ^{1/2},\
\vert (  f_{j0}^{\bot},v_{j})  _{\Omega_{j}(  h)
}\vert    \leq ch\mathcal{N}a(  v,v;\Xi(  h)  )
^{1/2}.\nonumber
\end{align}
Thus, the following assertion is proved.

\begin{proposition}
\label{prop4.M1}Under assumptions (\ref{4.16}) and (\ref{4.F}), (\ref{2.15})
the solution $u$ of problem (\ref{1.10})-(\ref{1.16}) with $\alpha=0$ and its
global approximation (\ref{4.66}), (\ref{4.67}) constructed in Section
\ref{sect3.3} are in relationship%
\begin{align}
\Vert \nabla_{x}&(  u_{0}-\mathbf{u}_{0})  ;L^{2}(
\Omega_{\bullet}(  h)  )  \Vert +h^{1/2}\Vert
r^{-1}(  1+\vert \ln r\vert )  ^{-1}(
u_{0}-\mathbf{u}_{0})  ;L^{2}(  \Omega_{\bullet}(  h)
)  \Vert+\label{4.70}\\
+&{\textstyle\sum\nolimits_{j}} (  \Vert \nabla_{x}(  u_{j}-\mathbf{u}%
_{j})  ;L^{2}(  \Omega_{j}(  h)  )  \Vert
+\Vert (  l_{j}-z)  ^{-1}(  u_{j}-\mathbf{u}_{j})
;L^{2}(  \Omega_{j}(  h)  )  \Vert )  \leq
ch(  1+\vert \ln h\vert )  ^{2}\mathcal{N},\nonumber
\end{align}
where $\mathcal{N}$ is the sum of norms of functions (\ref{4.F}) and $c$ is a
constant independent of $h\in(  0,h_{0}]  $ and $F_{0}%
,F_{1},...,F_{J}$ in (\ref{4.16}).
\end{proposition}

Let us analyze the error estimate (\ref{4.70}) and detect the valid asymptotic
formulas for $u_{0}(  h,x)  $ and $u_{j}(  h,x)  .$ A
direct calculation shows that in the case (\ref{2.15}) the main terms
$h^{-1}U_{0}^{-1}=h^{-1}a_{0}$ and $h^{-1}U_{j}^{-1}=h^{-1}a_{0}(
1-z/l_{j})  $ in the asymptotic ans\"{a}tze (\ref{3.34}), (\ref{3.35})
acquire the weighted norms%
\begin{equation}
h^{-1/2}\Vert r^{-1}(  1+\vert \ln r\vert )
^{-1}(  hU_{0}^{-1})  ;L^{2}(  \Omega_{\bullet}(
h)  )  \Vert    =O(  1)  ,\quad\Vert (  l_{j}-z)  ^{-1}(  h^{-1}U_{j}^{-1})
;L^{2}(  \Omega_{j}(  h)  )  \Vert   =O(
1)  .\label{4.71}
\end{equation}
Furthermore,%
\begin{equation}
\Vert \nabla_{x}(  hU_{0}^{-1})  ;L^{2}(  \Omega
_{\bullet}(  h)  )  \Vert =0,\ \ \ \Vert
\nabla_{x}(  h^{-1}U_{j}^{-1})  ;L^{2}(  \Omega_{j}(
h)  )  \Vert =\vert a_{0}\vert l_{j}^{-1}\vert
\omega_{j}\vert ^{1/2}. \label{4.72}%
\end{equation}
The Dirichlet (\ref{4.72}) and weighted (\ref{4.71}) norms of the boundary
layer%
\begin{equation}
\widetilde{w}_{j}^{0}(  \eta^{j},\zeta)  =-a_{0}\frac{\gamma
_{j}\vert \omega_{j}\vert }{l_{j}}\left\{
\begin{array}
[c]{c}%
\mathbf{w}^{j}(  \eta^{j},\zeta)  +(  2\pi)  ^{-1}%
\ln\rho_{j},\ \ \ (  \eta^{j},\zeta)  \in\Lambda_{j},\\
\mathbf{w}^{j}(  \eta^{j},\zeta)  -\gamma_{j}^{-1}\vert
\omega_{j}\vert ^{-1}\zeta-\mathbf{q}_{j},\ \ \ (  \eta^{j}%
,\zeta)  \in Q_{j},
\end{array}
\right.  \label{4.73}%
\end{equation}
become of order $h^{1/2}.$ Notice that, as concluded in Section \ref{sect2.3},
function (\ref{4.73}) has the exponential decay in the semi-cylinder $Q_{j}$
and a power-law decay in the layer $\Lambda_{j}.$

Unfortunately, all the above mentioned norms of the secondary terms $U_{0}%
^{0}$ and $U_{j}^{0}$ in the ans\"{a}tze are smaller than the bound in
(\ref{4.70}) and, hence, these terms cannot figure in the next assertion.

\begin{theorem}
\label{th4.M3}Under conditions (\ref{4.16}) and (\ref{4.F}), (\ref{2.15}) the
restrictions $u_{0}=u|_{\Omega_{\bullet}(  h)  }$ and
$u_{j}=u|_{\Omega_{j}(  h)  }$ of the solution $u$ of problem
(\ref{1.10})-(\ref{1.16}) admit the estimates%
\begin{align*}
&\Vert r^{-1}(  1+\vert \ln r\vert )  ^{-1}(
u_{0}-h^{-1}a_{0}-{\textstyle\sum\nolimits_{j}} \chi_{0j}\widetilde{w}_{j}^{0})
;L^{2}(  \Omega_{\bullet}(  h)  )  \Vert \leq
ch^{1/2}(  1+\vert \ln h\vert )  ^{2}\mathcal{N}%
,\\
&\Vert \nabla_{x}u_{j}-\nabla_{x}(  h^{-1}U_{j}^{-1}+\chi
_{j}\widetilde{w}_{j}^{0})  ;L^{2}(  \Omega_{j}(  h)
)  \Vert+\\
&\qquad\qquad+\Vert (  l_{j}-z)  (  u_{j}-h^{-1}U_{j}^{-1}-\chi
_{j}\widetilde{w}_{j}^{0})  ;L^{2}(  \Omega_{j}(  h)
)  \Vert \leq ch(  1+\vert \ln h\vert )
\mathcal{N},
\end{align*}
where $a_{0}$ is the constant (\ref{3.44}), $U_{j}^{-1}$ is the linear
function (\ref{3.l}) and the boundary layer term $\widetilde{w}_{j}^{0}$ is
defined in (\ref{4.73}).
\end{theorem}

We are in position to derive an assertion on convergence.

\begin{corollary}
\label{cor4.M4}Under conditions (\ref{4.16}) and (\ref{4.F}), (\ref{2.15}) the
following convergences occur:%
\begin{align}
hu(  h,P^{j}+h\eta^{j},z)   &  \rightarrow a_{0}(  1-l_{j}^{-1}z)  \text{ \ strongly in }H^{1}(  \omega_{j}\times(
0,l_{j})  )  ,\label{4.75}\\
hu(  h,y,h\zeta)   &  \rightarrow a_{0}\text{ \ strongly in }%
H^{1}(  \omega_{0}\times(  0,1)  )  .\nonumber
\end{align}
The constant $a_{0}$ is given in (\ref{3.44}) and the convergence rate in
(\ref{4.75}) is $O(  h^{1/2})  .$
\end{corollary}

The asymptotic procedure designed in Section \ref{sect3.3} allows continuation
so that lower-order terms in the outer (\ref{3.34}), (\ref{3.35}) and inner
(\ref{3.36}) expansions can be elucidated, in particular, the constant $A_{0}$
in (\ref{2.25}) and (\ref{3.47}) which fully molds the terms $U_{0}^{0}$ and
$U_{j}^{0}$. However, the presence of the small factor $h$ on the left-hand
side of the a priori estimate (\ref{4.13}) requires for a sufficient reduction
of discrepancies, namely bounds in estimates (\ref{est}) must become
$ch^{1+\delta}\mathcal{N}a(  v,v;\Xi(  h)  )  ^{1/2}.$
To this end, explicit formulas for $U_{0}^{1},$ $U_{j}^{1}$ and $w_{j}^{1}$
are needed as well as boundary layers near the soles $\omega_{j}^{h}(
l_{j})  $ of the rods (cf. Remark \ref{rem2.8}). We avoid such
complications in the present paper and leave the asymptotic terms (\ref{2.15})
and (\ref{3.47}) unproved in the case $\alpha=0.$ It should be emphasized that
the boundary layer (\ref{4.73}) makes impossible to obtain $U_{0}^{0}$ and
$U_{j}^{0}$ as a result of limit passages like $U_{0}^{-1}$ and $U_{j}^{-1}$
in Corollary \ref{cor4.M4}.

\section{Conclusive remarks\label{sect5}}

\subsection{On the structure of right-hand sides\label{sect5.1}}

In section \ref{sect4} we dealt with the simplified right-hand sides
(\ref{4.16}), (\ref{4.F}) of equations (\ref{1.10}) and (\ref{1.11}) only in
order to formulate Corollaries \ref{cor4.N4} and \ref{cor4.M4} with
convergence results. The asymptotic procedure expounded in Section \ref{sect3}
allows us to consider the asymptotic forms (\ref{2.8}), (\ref{2.9}) of
$f_{0}(  h,x)  $ and (\ref{2.1}), (\ref{2.2}) of $f_{j}(
h,x)  .$ For instance, if the norms%
\begin{align}
\mathcal{N}_{0}&=\left(  \int_{0}^{1}\int_{\omega_{0}}(  \vert
\nabla_{y}f_{0}^{\bot}(  y,\zeta)  \vert ^{2}+\vert
f_{0}^{\bot}(  y,\zeta)  \vert ^{2})  dyd\zeta\right)
^{1/2},\label{6.1}\\
\mathcal{N}_{j}&=\left(  \int_{\omega_{j}}\int_{0}^{l_{j}}(
\vert \partial_{z}f_{j}^{\bot}(  \eta,z)  \vert
^{2}+\vert f_{j}^{\bot}(  \eta,z)  \vert ^{2})
dzd\eta\right)  ^{1/2}\nonumber
\end{align}
are finite and, moreover,
\begin{align}
f_{0}^{\bot}(  y,\zeta)   &  =0\text{ in the vicinity of the points
}y=P^{j},\label{6.2}\\
f_{j}^{\bot}(  \eta,z)   &  =0\text{ in the vicinity of the points
}z=0\text{ and }z=l_{j},\nonumber
\end{align}
then Theorems \ref{th4.N3} and \ref{th4.M3} remain valid with the following
modifications. First, the asymptotic expansions of $u_{0}$ and $u_{j}$ must be
augmented with the terms $hU_{0}^{1}(  y,\zeta)  $ and $hU_{j}%
^{1}(  \eta^{j},z)  $ as in (\ref{2.10}) and (\ref{2.3}),
respectively. Second, in the bounds in estimates (\ref{4.55}), (\ref{4.61})
and (\ref{4.70}), (\ref{4.71}) the factor $\mathcal{N}$\ must be changed for
$\mathcal{N+N}_{0}+\mathcal{N}_{1}+...+\mathcal{N}_{J}$.
It should be stressed that the supplementary smoothness (\ref{6.1}) of
$f_{0}^{\bot}$ in $y$ and of $f_{j}^{\bot}$ in $z$ is needed to achieve the
inclusions $U_{0}^{1}\in H^{1}(  \omega_{0}\times(  0,1)
)  $ and $U_{j}^{1}\in H^{1}(  \omega_{j}\times(
0,l_{j})  )  $ while requirements (\ref{6.2}) are introduced in
order to avoid a modification of boundary layers described in Section
\ref{sect3}. We also emphasize that, in view of the dependence on the fast
variables, the Dirichlet norms of the new terms $hU_{p}^{1}$ get the same
order in $h$ as the old terms $U_{p}^{0}$ and, therefore, an asymptotics of
the solution $u$ of problem (\ref{1.10})-(\ref{1.16}) with the right-hand
sides (\ref{2.8}), (\ref{2.1}) cannot be written without the terms $hU_{p}%
^{1}$ which hamper in deducing a convergence result.

As usual, a smallness of the asymptotic remainders $\widetilde{f}_{0}$ in
(\ref{2.8}) and $\widetilde{f}_{j}$ in (\ref{2.1}) is expressed in such a way
that according to the weighted anisotropic inequality (\ref{4.13}) the scalar
products $(  \widetilde{f}_{0},v_{0})  _{\Omega_{\bullet}(
h)  }$ and $(  \widetilde{f}_{j},v_{j})  _{\Omega_{j}(
h)  }$ get the same order in $h$ as the bound in an error estimate in
Propositions \ref{prop4.N1} and \ref{prop4.M1}.

We also refer to the monograph \cite[Ch. 2 and 5]{MaNaPl} where basic
principles of multi-scaled representations for right-hand sides are laid. In
particular, components which are written in the fast variables (\ref{2.36})
and decay as $\vert \xi^{j}\vert \rightarrow+\infty,$ may be
introduced into (\ref{2.8}) and (\ref{2.1}) to be taken into account in
inhomogeneous problems (\ref{2.28})-(\ref{2.33}) or
(\ref{2.57}) and (\ref{2.58}). No changes in the asymptotic procedure are then needed.

\subsection{On the shape of the junction\label{sect5.2}}

\begin{figure}
[ptb]
\begin{center}
\includegraphics[
height=0.7109in,
width=3.5215in
]%
{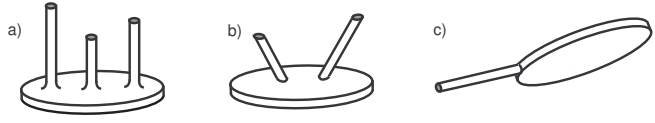}%
\caption{The smoothed junction zones (a), the skewed (b) and planar (c)
junctions.}%
\label{f4}%
\end{center}
\end{figure}

Asymptotic procedures of dimension reduction are known to support much more
general geometry of rods and plates, cf. \cite{Ciarlet, Sanchez, MaNaPl,
Panas, Nabook} and others. In particular, the rods $\Omega_{j}(
h)  $ may have varying cross-sections and distorted ends while the plate
$\Omega_{0}(  h)  $ may be with smoothly variable thickness and
indented edge. Local perturbations of junction zones, fig.\ref{f4},a, are
available (see \cite[Ch. 5 and 11]{MaNaPl}). These modifications do
not affect the performed asymptotic analysis but make the justification
schemes much more involved, although the above-used technique of cut-off
functions with overlapping supports still works.%

%TCIMACRO{\FRAME{ftbpFU}{2.3627in}{1.3387in}{0pt}{\Qcb{The smoothed junction
%zones (a), the skewed (b) and planar (c) junctions.}}{\Qlb{f4}}%
%{junction4new.eps}{\special{ language "Scientific Word";  type "GRAPHIC";
%width 2.3627in;  height 1.3387in;  depth 0pt;  original-width 6.6945in;
%original-height 3.7818in;  cropleft "0";  croptop "1";  cropright "1";
%cropbottom "0";  filename 'Junction4new.eps';file-properties "XNPEU";}} }%
%BeginExpansion
%EndExpansion

Even the cylindrical elements $\Omega_{0}(  h)  =\omega_{0}%
\times(  0,h)  $ and $\Omega_{j}(  h)  =\omega_{j}%
^{h}\times(  0,l_{j})  $ can be joined in a way different from the
junction $\Xi(  h)  $ in formula (\ref{1.4}) and fig. \ref{f1},a,
compare a skewed and planar junctions drawn in fig. \ref{f4},b and c. The
construction of the asymptotics in the rods and plate follows Sections
\ref{sect2.1} and \ref{sect2.2} directly but certain changes in our analysis
of the main terms in the inner expansions occur. For example, in the case of
skewed junction in fig. \ref{f4},b, problem (\ref{2.57}) is posed in a layer
with an inclined shaft where separation of variables is impossible and the
special solution $\mathbf{W}^{j}$ with decomposition (\ref{2.60}) becomes
three-dimensional. For the planar junction in fig. \ref{f4},c, the limit
problem of type (\ref{2.57}) is set in the half-layer $\Lambda_{+}=\{
(  \eta_{1},\eta_{2},\zeta)  :\eta_{1}<0,\ \zeta\in(
0,1)  \}  $ and the limit problem of type (\ref{2.28}%
)-(\ref{2.33}) in the union of $\Lambda_{+}$ and a semi-cylinder with the
$\eta_{1}$-axis. In order to make necessary conclusions on the asymptotic
behavior of solutions of these problems, we refer to the Kondratiev theory
\cite{Ko} (see also \cite[Ch.2 and 5]{NaPl}) and results in \cite{na164,
na243} about the cylinder and layer-like outlets to infinity. All other steps
in the asymptotic procedure remain without any alteration.

All the previous results admit slight and self-understood modifications
in the case of arbitrary elliptic second-order differential operators with smooth coefficients in $\Omega_0(1)$ and $\Omega_j(h)$, cf. \cite[Ch.1]{Nabook}.

\subsection{The Dirichlet condition on the lateral side of the
plate.\label{sect5.4}}

Let us outline certain primary results on the asymptotic structures of the
solution $u(  h,x)  $ of problem (\ref{1.10}), (\ref{1.11}),
(\ref{1.d}), (\ref{1.13})-(\ref{1.16}) with exponent (\ref{1.18}). As
explained in Sections \ref{sect3.2} and \ref{sect3.4}, the asymptotic and
justification procedures become much more plain and simple. Note that the role
of the anisotropic inequality (\ref{4.13}) is now passed over to the standard
one%
\begin{equation}
\Vert r^{-1}(  1+\vert \ln r\vert )  ^{-1}%
v_{0};L^{2}(  \Omega_{\bullet}(  h)  )  \Vert
+h^{-\alpha}{\textstyle\sum\nolimits_{j}} \Vert (  l_{j}-z)  ^{-1}v_{j}%
;L^{2}(  \Omega_{j}(  h)  )  \Vert \leq c_{\Xi
}a(  v,v;\Xi(  h)  )  . \label{6.10b}%
\end{equation}
We do not comment on proofs of the following assertions which are but a
simplified version of contents of Sections \ref{sect3} and \ref{sect4}.

\begin{theorem}
\label{th5.7}Under conditions (\ref{4.16}) and (\ref{4.F}), the solution
$u(  h,x)  $ of problem (\ref{1.10}), (\ref{1.11}), (\ref{1.d}),
(\ref{1.13})-(\ref{1.16}) with $\alpha=1$ satisfies%
\begin{equation}
\Vert \nabla_{x}(  u_{0}-U_{0}^{0}-{\textstyle\sum\nolimits_{j}} \chi_{0j}%
A_{j}\widetilde{\mathbf{W}}^{j})  ;L^{2}(  \Omega_{\bullet}(
h)  )  \Vert +h^{-1}{\textstyle\sum\nolimits_{j}} \Vert \nabla_{x}%
u_{j}-\nabla_{x}U_{j}^{0};L^{2}(  \Omega_{j}(  h)  )
\Vert \leq ch\mathcal{N}, \label{6.11}%
\end{equation}
where $U_{0}^{0}(  y,\ln h)  $ and $A_{j}(  \ln h)  $
are taken from (\ref{3.Dir}) and (\ref{3.AD}), $A_{j}(  \ln h)
\widetilde{\mathbf{W}}^{j}(  \eta^{j},z)  $ is the boundary layer
(\ref{4.60}) and $U_{j}^{0}(  z,\ln h)  $ is the solution of the
mixed boundary value problem (\ref{2.6}), (\ref{2.7}), (\ref{3.1}),
(\ref{3.G}). The constant $c$ in (\ref{6.11}) does not depend on $h\in(
0,h_{0}]  $ and the sum $\mathcal{N}$\ of norms of functions (\ref{4.F}).
\end{theorem}

Estimates of weighted Lebesgue norms are inherited from (\ref{6.11}) and
(\ref{6.10b}). Since due to (\ref{3.AD}) the coefficients $A_{j}(  \ln
h)  $ on the Green functions $G_{j}$ in the linear combination
(\ref{3.Dir}) are small, Theorem \ref{th5.7} ensures that%
\begin{align}
u(  h,y,h\zeta)   &  \rightarrow U_{0}^{\#}(  y)  \text{
\ strongly in }L^{2}(  \omega_{0}\times(  0,1)  )
,\label{6.12}\\
u(  h,P^{j}+h\eta^{j},z)   &  \rightarrow U_{j}^{\#}(
z)  \text{ \ strongly in }L^{2}(  \omega_{j}\times(
0,l_{j})  )  \label{6.13}%
\end{align}
where $U_{0}^{\#}\in H^{2}(  \omega_{0})  $ is a solution of the
Dirichlet problem (\ref{3.DO}) and $U_{j}^{\#}\in H^{2}(  0,l_{j})
$ is a solution of problem (\ref{2.6}), (\ref{2.7}) with the homogeneous
Neumann condition $\partial_{z}U_{j}^{\#}(  0)  =0.$ The
convergence rate in (\ref{6.12}) and (\ref{6.13}) is of order $\vert \ln
h\vert ^{-1}.$ A direct calculation considering the logarithmic
singularities of the Green functions demonstrates that the Dirichlet norm on
the intact rescaled plate $\omega_{0}\times(  0,1)  $ is also
infinitesimal, namely%
\[
\Vert \nabla_{y}u;L^{2}(  \omega_{0}\times(  0,1)
)  \Vert ^{2}+\Vert \partial_{\zeta}u;L^{2}(  \omega
_{0}\times(  0,1)  )  \Vert ^{2}\leq c\vert \ln
h\vert ^{-1}.
\]
Thus, the strong convergence (\ref{6.12}) occurs in $H^{1}(  \omega
_{0}\times(  0,1)  )  ,$ too, however with rate $\vert
\ln h\vert ^{-1/2}$ only.

\begin{theorem}
\label{th5.8}Under conditions (\ref{4.16}) and (\ref{4.F}), the solution of
problem (\ref{1.10}), (\ref{1.11}), (\ref{1.d}), (\ref{1.13})-(\ref{1.16})
with $\alpha=0$ satisfies%
\begin{equation}
\Vert \nabla_{x}(  u_{0}-U_{0}^{\#})  ;L^{2}(
\Omega_{\bullet}(  h)  )  \Vert +{\textstyle\sum\nolimits_{j}}\Vert \nabla_{x}(  u_{j}-U_{j}^{\#})  ;L^{2}(
\Omega_{j}(  h)  )  \Vert \leq ch^{3/2}\mathcal{N},
\label{6.14}%
\end{equation}
where $U_{0}^{\#}\in H^{2}(  \omega_{0})  $ and $U_{j}^{\#}\in
H^{2}(  0,l_{j})  $ are solutions of the problems
(\ref{3.DO}) and (\ref{2.6}), (\ref{2.7}) (\ref{3.DD}).
\end{theorem}

The error estimate (\ref{6.14}) also exhibits the strong convergences
(\ref{6.12}) in $H^{1}(  \omega_{0}\times(  0,1)  )  $
and (\ref{6.13}) in $H^{1}(  \omega_{j}\times(  0,l_{j})
)  $ which have been observed in \cite{G1} but without detecting the
convergence rates. We again emphasize the obvious difference in our results in
Theorems \ref{th4.N3} and \ref{th5.7} for $\alpha=1$ and Theorems \ref{th4.M3}
and \ref{th5.8} for $\alpha=0$ which serve, respectively, for the cases of the
Neumann and Dirichlet conditions on the lateral side $\upsilon_{0}(
h)  $ of the intact plate (\ref{1.1}).

%\bigskip
%
%\textbf{Acknowledgements}
%
%\medskip
%
%The third author have been supported by the Gruppo Nazionale per l' Analisi Matematica, la Probabilit\`{a} e le loro Applicazioni (GNAMPA) of the Istituto Nazionale di Alta Matematica (INdAM). The second author is member of GNAMPA of the INDAM.

\end{document}